\documentclass[a4paper]{amsart}

\usepackage{amsmath}
\usepackage{amssymb}
\usepackage{latexsym}
\usepackage[mathscr]{eucal}
\usepackage{mathrsfs}
\usepackage{graphics}
\usepackage{newtxtext}
\usepackage{newtxmath}
\usepackage{ascmac,color}
\usepackage{framed}
\usepackage{ulem}
\usepackage{cancel}
\usepackage{arydshln}
\usepackage{bm}
\usepackage{longtable}
\allowdisplaybreaks
\usepackage{multicol}
\usepackage{type1cm}
\usepackage{booktabs}
\usepackage{amscd}

\theoremstyle{plain}
\newtheorem{theorem}{\bf Theorem}[subsection]
\newtheorem{lemma}[theorem]{\bf Lemma}
\newtheorem{proposition}[theorem]{\bf Proposition}

\newtheorem*{claim}{\bf Claim}

\newcommand{\relmiddle}[1]{\mathrel{}\middle#1\mathrel{}}

\makeatletter
\def\thefootnote{\ifnum\c@footnote>\z@\leavevmode\lower.5ex%
      \hbox{$^{\@arabic\c@footnote)}$}\fi}
\makeatother

\theoremstyle{definition}

\theoremstyle{remark}

\def\Aut{\mbox{\rm {Aut}}}
\def\det{\mbox{\rm {det}}}

\def\diag{\mbox{\rm {diag}}}

\def\Hom{\mbox{\rm {Hom}}}

\def\Iso{\mbox{\rm {Iso}}}

\def\Ker{\mbox{\rm {Ker}}}

\def\tr{\mbox{\rm {tr}}}

\def\ov{\overline}

\def\dfrac#1#2{\displaystyle \frac{#1}{#2}}

\def\C{\mbox{\boldmath $C$}}

\def\H{\mbox{\boldmath $H$}}

\def\R{\mbox{\boldmath $R$}}

\def\Z{\mbox{\boldmath $Z$}}

\def\sR{\mbox{\boldmath $\scriptstyle{R}$}}
\def\sC{\mbox{\boldmath $\scriptstyle{C}$}}
\def\sH{\mbox{\boldmath $\scriptstyle{H}$}}

\def\0{\mbox{\boldmath {0}}}
\def\1{\mbox{\boldmath {1}}}
\def\2{\mbox{\boldmath {2}}}
\def\3{\mbox{\boldmath {3}}}
\def\4{\mbox{\boldmath {4}}}
\def\5{\mbox{\boldmath {5}}}
\def\6{\mbox{\boldmath {6}}}
\def\7{\mbox{\boldmath {7}}}
\def\8{\mbox{\boldmath {8}}}
\def\9{\mbox{\boldmath {9}}}

\def\a{\mbox{\boldmath $a$}}

\def\c{\mbox{\boldmath $c$}}

\def\f{\mbox{\boldmath $f$}}

\def\i{\mbox{\boldmath $i$}}
\def\j{\mbox{\boldmath $j$}}

\def\x{\mbox{\boldmath $x$}}

\def\tr{\mbox{\rm {tr}}}

\def\det{\mbox{\rm {det}}}
\def\diag{\mbox{\rm {diag}}}
\def\Iso{\mbox{\rm {Iso}}}
\def\Hom{\mbox{\rm {Hom}}}
\def\Ker{\mbox{\rm {Ker}}}
\def\ov{\overline}

\def\sR{\mbox{\boldmath $\scriptstyle{R}$}}
\def\sC{\mbox{\boldmath $\scriptstyle{C}$}}
\def\sH{\mbox{\boldmath $\scriptstyle{H}$}}
\def\dfrac#1#2{\displaystyle \frac{#1}{#2}}

\def\C{\mbox{\boldmath $C$}}

\def\H{\mbox{\boldmath $H$}}

\def\R{\mbox{\boldmath $R$}}

\def\Z{\mbox{\boldmath $Z$}}
\def\sR{\mbox{\boldmath $\scriptstyle{R}$}}

\def\0{\mbox{\boldmath {0}}}
\def\1{\mbox{\boldmath {1}}}
\def\2{\mbox{\boldmath {2}}}
\def\3{\mbox{\boldmath {3}}}
\def\4{\mbox{\boldmath {4}}}
\def\5{\mbox{\boldmath {5}}}
\def\6{\mbox{\boldmath {6}}}
\def\7{\mbox{\boldmath {7}}}
\def\8{\mbox{\boldmath {8}}}
\def\9{\mbox{\boldmath {9}}}
\def\a{\mbox{\boldmath $a$}}

\def\c{\mbox{\boldmath $c$}}

\def\f{\mbox{\boldmath $f$}}

\def\i{\mbox{\boldmath $i$}}
\def\j{\mbox{\boldmath $j$}}

\def\x{\mbox{\boldmath $x$}}

\definecolor{cyan}{cmyk}{1,0,0,0}
\definecolor{lightcyan}{cmyk}{0.5,0,0,0}
\definecolor{pastelcyan}{cmyk}{0.25,0,0,0}
\definecolor{magenta}{cmyk}{0,1,0,0}
\definecolor{yellow}{cmyk}{0,0,1,0}
\definecolor{lightyellow}{cmyk}{0,0,0.5,0}
\definecolor{pastelyellow}{cmyk}{0,0,0.25,0}
\definecolor{black}{cmyk}{0,0,0,1}
\definecolor{darkgray}{cmyk}{0,0,0,0.75}
\definecolor{gray}{cmyk}{0,0,0,0.5}
\definecolor{lightgray}{cmyk}{0,0,0,0.25}
\definecolor{white}{cmyk}{0,0,0,0}
\definecolor{red}{cmyk}{0,1,1,0}

\definecolor{orange}{cmyk}{0,0.5,1,0}
\definecolor{redorange}{cmyk}{0,0.77,0.87,0.35}
\definecolor{scarlet}{cmyk}{0,1,0.5,0}
\definecolor{brown}{cmyk}{0.5,0.75,1,0}
\definecolor{camel}{cmyk}{0.25,0.375,0.5,0}
\definecolor{cream}{cmyk}{0,0.2,0.3,0}
\definecolor{green}{cmyk}{1,0,1,0}
\definecolor{llightgreen}{cmyk}{0.55,0,0.45,0}
\definecolor{lightgreen}{cmyk}{0.5,0,0.5,0}
\definecolor{pastelgreen}{cmyk}{0.25,0,0.25,0}
\definecolor{mossgreen}{cmyk}{0.64,0.4,1,0}
\definecolor{yellowgreen}{cmyk}{0.5,0,1,0}
\definecolor{skyblue}{cmyk}{0.4,0.16,0,0}
\definecolor{royal}{cmyk}{1.0,0.5,0,0}
\definecolor{navyblue}{cmyk}{0.9,0.75,0.5,0}
\definecolor{navy}{cmyk}{0.9,0.8,0.5,0}
\definecolor{lightgray}{cmyk}{0,0,0,0.02}

\definecolor{lightnavy}{cmyk}{0.4,0.3,0.2,0}
\definecolor{blue}{cmyk}{1,1,0,0}
\definecolor{lightblue}{cmyk}{0.5,0.5,0,0}
\definecolor{lavender}{cmyk}{0.25,0.25,0,0}
\definecolor{violet}{cmyk}{0.75,1,0.25,0}
\definecolor{purple}{cmyk}{0.5,1,0.5,0}
\definecolor{pink}{cmyk}{0,0.5,0,0}
\definecolor{pastelpink}{cmyk}{0,0.25,0,0}
\definecolor{Midnightblue}{cmyk}{0.98,0.13,0,0.43}
\definecolor{skyblue2}{cmyk}{0.1,0.1,0,0}
\definecolor{Cerulean}{cmyk}{0.9,0.08,0,0}
\definecolor{ProcessBlue}{cmyk}{0.7,0,0,0}

\def\skyblue{\color{skyblue}}

\def\skyblue2{\color{skyblue2}}

\usepackage{epic,eepic}

\allowdisplaybreaks

\begin{document}

\title[On realizations of Lie groups $ (E_{6,\sR})^C, (E_{6,\sC})^C, (E_{6,\sH})^C $]
{On realizations of the complex Lie groups $ (E_{6,\sR})^C, (E_{6,\sC})^C, (E_{6,\sH})^C $ \\
and those real forms}

\author[Toshikazu Miyashita]{Toshikazu Miyashita}


\begin{abstract}
There exist six Lie groups of type $ E_6 $, and to be specific, $ {E_6}^C , E_6, E_{6(6)}, E_{6(-2)}, E_{6(-14)}, \allowbreak E_{6(-26)}$. In order to define these groups, we use usually the Cayley algebra $ \mathfrak{C} $ and the split Cayley algebra $ \mathfrak{C}' $. In the present article, we consider the Lie groups which are defined by replacing $ \mathfrak{C}^C, \mathfrak{C} $ and $ \mathfrak{C}' $ with the fields of real numbers $ \R $, complex numbers $ \C $, split complex numbers $ \C'$, quaternions $ \H $ and split quaternions $ \H' $. For instance, the group $ (E_{6,\sR})^C $ is given as a group defined by replacing $ \mathfrak{C} $ with $ \R $ in $ {E_6}^C $ and the group $ E_{6(-26),\sH} $ is given as a group defined by replacing $ \mathfrak{C} $ with $ \H $ in $ E_{6(-26)} $. We call { \it realization} to determine the structure of the group.
\end{abstract}

\subjclass[2010]{}

\keywords{exceptional Lie groups}

\address{1365-3 Bessho onsen      \endgraf
	     Ueda City                 \endgraf
	     Nagano Pref. 386-1431     \endgraf
	     Japan}
\email{anarchybin@gmail.com}

\maketitle

\setcounter{section}{0}

\section{Introduction}
\begin{center}
\begin{tabular}{p{5em}lll}
\hline
&&&
\\[-2mm]
  & $ K=\R $ & \hspace{20mm}$ \C $ & \hspace{5mm}$ \H $
\\[1mm]
\midrule[0.05mm]

\\[-2mm]
$ (E_{6,K})^C $ & $ SL(3,C) $ & $ (SU(3,\C^C) \times SU(3,\C^C))/\Z_3 \rtimes \Z_2 $ & $ SU(6,\C^C)/\Z_2 $
\\[1mm]
$ E_{6,K} $ & $ SU(3) $ & $ (SU(3) \times SU(3))/\Z_3 \rtimes \Z_2 $ & $ SU(6)/\Z_2 $
\\[1mm]
$ E_{6(2),K} $ & $ SU(3) $ & $ (SU(3) \times SU(3))/\Z_3 \rtimes \Z_2 $ & $ SU(6)/\Z_2 $
\\[1mm]
$ E_{6(-14),K} $ & $ SU(1,2) $ & $ (SU(1,2) \times SU(1,2))/\Z_3 \rtimes \Z_2 $ & $ SU(2,4)/\Z_2 $
\\[1mm]
$ E_{6(-26),K} $ & $ SL(3,\R) $ & $ SU(3,\C^C)/\Z_3 \rtimes \Z_2 $ & $ SU^*(6)/\Z_2 $
\\[1mm]
\hline
\end{tabular}
\end{center}
\vspace{3mm}

\begin{center}
\begin{tabular}{p{5em}l@{\hspace{8mm}}l@{\hspace{8mm}}l}
\hline
&&&
\\[-2mm]
  & \;\;$ K=\R $ & \hspace{18mm}$ \C' $ & \hspace{5mm} $ \H' $
\\[1mm]
\midrule[0.05mm]
&&&
\\[-2mm]
$ E_{6(6),K} $ & $ SL(3,\R) $ & $ (SU(3,\C') \times SU(3,\C')) \rtimes \Z_2  $ & $ SL(6,\R) \rtimes \Z_2 $
\\[1mm]
\hline
\end{tabular}
\end{center}

\section{Preliminaries}

\section{The complex Lie group $ (E_{6,\sR})^C $ and its real forms}

We define the group $ (E_{6,\sR})^C $ by
\begin{align*}
(E_{6,\sR})^C:&=\left\lbrace \alpha \in \Iso_{C}(\mathfrak{J}(3,\R^C))\relmiddle{|} \det(\alpha X)=\det\,X \right\rbrace
\\
&=\left\lbrace \alpha \in \Iso_{C}(\mathfrak{J}(3,\R^C))\relmiddle{|} (\alpha X,\alpha Y, \alpha Z)=(X,Y,Z) \right\rbrace
\\
&=\left\lbrace \alpha \in \Iso_{C}(\mathfrak{J}(3,\R^C))\relmiddle{|} \alpha X \times \alpha Y={}^t\!\alpha^{-1}(X \times Y) \right\rbrace.
\end{align*}
and its real forms are defined as follows:
\begin{align*}
E_{6,\sR}&:=\left\lbrace \alpha \in \Iso_{C}(\mathfrak{J}(3,\R^C))\relmiddle{|} \det(\alpha X)=\det\,X, \langle \alpha X, \alpha Y \rangle=\langle X,Y \rangle \right\rbrace,
\\
E_{6(-14),\sR}&:=\left\lbrace \alpha \in \Iso_{C}(\mathfrak{J}(3,\R^C))\relmiddle{|} \det(\alpha X)=\det\,X, \langle \alpha X, \alpha Y \rangle_\sigma=\langle X,Y \rangle_\sigma \right\rbrace,
\\
E_{6(-26),\sR}&:=\left\lbrace \alpha \in \Iso_{\sR}(\mathfrak{J}(3,\R))\relmiddle{|} \det(\alpha X)=\det\,X \right\rbrace,
\end{align*}
where $ \langle X,Y \rangle=(\tau X,Y), \langle X,Y \rangle_\sigma=(\tau\sigma X,Y) $,
and since $ E_{6(6),\sR}=E_{6(-26),\sR} $ and $ E_{6(2),\sR}=E_{6,\sR}$, these groups are omitted.

The structure of the group $ (E_{6,\sR})^C $ has been already determined as follows.

\begin{theorem}{\rm \cite[Theorem 5.0.5]{miya0}}\label{theorem 3.1}
The group $ (E_{6,\sR})^C $ is isomorphic to the group $ SL(3,C) ${\rm :} $ (E_{6,\sR})^C \allowbreak \cong SL(3,C) $.
\end{theorem}
\begin{proof}
We define a mapping $ f_{6,C}: SL(3,C) \to (E_{6,\sR})^C $ by
\begin{align*}
f_{6,C}(A)X=AX{}^t\!A, \;\; X \in \mathfrak{J}(3,\R^C).
\end{align*}
This mapping induces the isomorphism $ (E_{6,\sR})^C \allowbreak \cong SL(3,C) $.
\end{proof}

Subsequently, the structure of the group $ E_{6,\sR} $ also has been already determined as follows. The proof is omitted, we describe only the result obtained.

\begin{theorem}{\rm \cite[Theorem 5.0.9]{miya0}}\label{theorem 3.2}
The group $ E_{6,\sR} $ is isomorphic to the group $ SU(3) ${\rm :} $ E_{6,\sR} \cong SU(3) $.
\end{theorem}
Note that the author has proved under the definition of $ SU(3) $ which is defined by $ \{A \!\in M(3,C)| (\tau{}^t\!A)A\allowbreak =E,\det\,A=1\} $ in  Theorem \ref{theorem 3.2}.

Now, after preparation, we will determine the structure of the group $ E_{6(-14), \sR} $.

We define an involutive automorphism $ \lambda $ of $ (E_{6,\sR})^C $ by
\begin{align*}
\lambda(\alpha)={}^t\!\alpha^{-1},\;\;\alpha \in (E_{6,\sR})^C,
\end{align*}
where the transpose $ {}^t\!\alpha $ of $ \alpha $ is defined by $ ({}^t\!\alpha X,Y)=(X,\alpha Y) $. Moreover we define an $ \R $-linear transformation $ \sigma $ of $ \mathfrak{J}(3,\R) $ by
\begin{align*}
\sigma X=\begin{pmatrix}
\xi_1 & -x_3 & -\ov{x}_2 \\
-\ov{x}_3 & \xi_2 & x_1  \\
-x_2 & \ov{x}_1 & \xi_3
\end{pmatrix},\,\, X \in \mathfrak{J}(3,\R).
\end{align*}
Then we have $ \sigma \in E_{6,\sR} \subset (E_{6,\sR})^C $ and $ \sigma^2=1 $.
This transformation $ \sigma $ is naturally extended to the $ C $-linear transformation of $ \mathfrak{J}(3,\R^C) $. Let $ \tau $ be the complex conjugation in $ \mathfrak{J}(3,\R^C) $, then we consider an involutive automorphism $ \tilde{\tau\lambda\sigma} $ of $ (E_{6,\sR})^C $: $ \tilde{\tau\lambda\sigma}(\alpha)=(\tau\sigma)\lambda(\alpha)(\sigma\tau), \alpha \in (E_{6,\sR})^C $. Note that as in the proof of \cite[Lemma 3.2.1]{iy7}, we can prove $ \tilde{\tau\lambda\sigma}(\alpha) \in (E_{6,\sR})^C $.

\noindent Hence we can define a subgroup $ ((E_{6,\sR})^C)^{\tau\lambda\sigma} $ of $ (E_{6,\sR})^C $ by
\begin{align*}
((E_{6,\sR})^C)^{\tau\lambda\sigma}:=\left\lbrace \alpha \in (E_{6,\sR})^C \relmiddle{|} \tilde{\tau\lambda\sigma}(\alpha)=\alpha \right\rbrace.
\end{align*}

Then we have the following proposition.

\begin{proposition}\label{proposition 3.3}
The group  $ ((E_{6,\sR})^C)^{\tau\lambda\sigma} $ coincides with the group $ E_{6(-14),\sR} ${\rm :} $ ((E_{6,\sR})^C)^{\tau\lambda\sigma}= E_{6(-14),\sR} $.
\end{proposition}
\begin{proof}
Let $ \alpha \in ((E_{6,\sR})^C)^{\tau\lambda\sigma} $. Then it follows that
\begin{align*}
\langle \alpha X, \alpha Y \rangle_\sigma &=(\tau\sigma \alpha X, \alpha Y)=({}^t\!\alpha \tau\sigma \alpha X,Y)=(\tau\sigma \alpha^{-1}\alpha X,Y)=(\tau\sigma X,Y)=\langle X,Y \rangle_\sigma,
\end{align*}
that is, $ \langle \alpha X, \alpha Y \rangle_\sigma=\langle X,Y \rangle_\sigma $.
Hence we have $ \alpha \in E_{6(-14),\sR} $. Conversely, let $ \beta \in E_{6(-14),\sR} $. It follows from $ \langle \beta X, \beta Y \rangle_\sigma=\langle X,Y \rangle_\sigma $ that
\begin{align*}
(\tau\sigma X,Y)=(\tau\sigma\beta X,\beta Y)=({}^t\!\beta\tau\sigma\beta X,Y)\;\; \text{for any}\;\; Y \in \mathfrak{J}(3,\R^C).
\end{align*}
Hence we have $ \tau\sigma={}^t\!\beta\tau\sigma\beta $, that is, $ (\tau\sigma){}^t\!\beta^{-1}(\sigma\tau)=\beta $. Hence we have $ \beta \in ((E_{6,\sR})^C)^{\tau\lambda\sigma} $. With above, the proof of this proposition is completed.
\end{proof}

We prove the lemma needed in the proof of theorem below.

\begin{lemma}\label{lemma 3.4}
The mapping $ f_{6,C}:SL(3,C) \to (E_{6,\sR})^C $ of Theorem {\rm 3.1} satisfies
\begin{align*}
(1)\;\; {}^t\!{f_{6,C}(A)}^{-1}=f_{6,C}({}^t\! A^{-1}).\quad (2)\;\; \sigma f_{6,C}(A)\sigma =f_{6,C}(I_1AI_1)\quad (3)\;\; \tau f_{6,C}(A)\tau =f_{6,C}(\tau A),
\end{align*}
where $ I_1:=\diag(-1,1,1) \in M(3,\R) $.
\end{lemma}
\begin{proof}
(1) It follows that
\begin{align*}
({}^t\!f_{6,C}(A)X,Y)=(X,f_{6,C}(A)Y )=(X,AY{}^t\!A)=({}^t\!AXA,Y)
=(f_{6,C}({}^t\!A)X,Y)
\end{align*}
that is, $ {}^t\!f_{6,C}(A)=f_{6,C}({}^t\!A) $. Hence we have $ {}^t\!{f_{6,C}(A)}^{-1}=f_{6,C}({}^t\! A^{-1}) $.
\vspace{1mm}

It is easy to verify that (2) and (3) hold.

\end{proof}

We determine the structure of the group $ E_{6(-14),\sR} $.

\begin{theorem}\label{theorem 3.5}
The group $ E_{6(-14),\sR} $ is isomorphic to the group $ SU(1,2) ${\rm :} $ E_{6(-14),\sR} \cong SU(1,2) $.
\end{theorem}
\begin{proof}
Let the group $ SU(1,2)=\{A \in M(3,C)\,|\, AI_1(\tau{}^t\!A)=I_1, \det\,A=1 \} $ and the group $ E_{6(-14),\sR} $ as the group $ ((E_{6,\sR})^C)^{\tau\lambda\sigma} $ (Proposition \ref{proposition 3.3}). Then we define a mapping $ \varphi_{6(-14),\sR}: SU(1,2) \to ((E_{6,\sR})^C)^{\tau\lambda\sigma} $ by
\begin{align*}
\varphi_{6(-14),\sR}(A)X=AX{}^t\!A,\;\; X \in \mathfrak{J}(3,\R^C).
\end{align*}
Note that this mapping is the restriction of the mapping $ f_{6,C} $ (Theorem \ref{theorem 3.1}).
First, we will prove that $ \varphi_{6(-14),\sR} $ is well-defined. Since $ SU(1,2) \subset SL(3,C) $, it is easy to see $ \varphi_{6(-14),\sR}(A) \in (E_{6,\sR})^C $. Moreover, from Lemma \ref{lemma 3.4} we have
\begin{align*}
(\tau\sigma)\lambda(\varphi_{6(-14),\sR}(A))(\sigma\tau)=\varphi_{6(-14),\sR}(I_1(\tau{}^t\! A^{-1})I_1).
\end{align*}
Hence, since we have $ I_1(\tau{}^t\! A^{-1})I_1=A $ from $ A \in SU(1,2) $, we obtain $ (\tau\sigma)\lambda(\varphi_{6(-14),\sR}(A))(\sigma\tau)=\varphi_{6(-14),\sR}(A) $,
that is, $ \varphi_{6(-14),\sR}(A) \in ((E_{6,\sR})^C)^{\tau\lambda\sigma}$, so that $ \varphi_{6(-14),\sR} $ is well-defined. Subsequently, we will prove that $ \varphi_{6(-14),\sR} $ is a homomorphism, however since $ \varphi_{6(-14),\sR} $ is the restriction of the mapping $ f_{6,C} $, it is clear.

Next, we will prove that $ \varphi_{6(-14),\sR} $ is surjective. Let $ \alpha \in ((E_{6,\sR})^C)^{\tau\lambda\sigma} \subset (E_{6,\sR})^C $. Then there exists $ P \in SL(3,C) $ such that $ \alpha=f_{6,C}(P) $ (Theorem \ref{theorem 3.1}). Since $ \alpha $ satisfies the condition $ (\tau\sigma)\lambda(\alpha)(\sigma\tau)=\alpha $, that is, $ (\tau\sigma)\lambda(f_{6,C}(P))(\sigma\tau)=f_{6,C}(P) $, it follows from Lemma \ref{lemma 3.4} that
\begin{align*}
(\tau\sigma)\lambda(f_{6,C}(P))(\sigma\tau)=f_{6,C}(I_1(\tau{}^t\! P^{-1})I_1).
\end{align*}
Hence we have the following
\begin{align*}
I_1(\tau{}^t\! P^{-1})I_1=P,
\end{align*}
so that $ P $ satisfies $ PI_1(\tau{}^t\!P)=I_1 $ and $ \det\,P=1 $, that is, $ P \in SU(1,2) $. Thus there exists $ A \in SU(1,2) $ such that $ \alpha=f_{6,C}(A)=\varphi_{6(-14),\sR}(A) $. With above, the proof of surjective is completed.

Finally, we will determine $ \Ker\,\varphi_{6(-14),\sR} $. since $ \varphi_{6(-14),\sR} $ is the restriction of the mapping $ f_{6,C} $, it is easy to obtain $ \Ker\,\varphi_{6(-14),\sR}=\{E\} $.

Therefore we have the required isomorphism
\begin{align*}
 E_{6(-14),\sR} \cong SU(1,2).
\end{align*}
\end{proof}

Again, let $ \tau $ be complex conjugation in $ \mathfrak{J}(3,\R^C) $, then $ \tau $ induces an involutive automorphism $ \tilde{\tau} $ of $ (E_{6,\sR})^C $: $ \tilde{\tau}(\alpha)=\tau\alpha\tau, \alpha \in (E_{6,\sR})^C $. Indeed, by using the formulas $ (\tau X, \tau Y)=\tau(X,Y),\tau X \times \tau Y=\tau(X \times Y), X,Y \in \mathfrak{J}(3,\R^C) $, it follows from $ (\alpha X,\alpha Y,\alpha Z)=(X,Y,Z) $ that
\begin{align*}
(\tau\alpha\tau X,\tau\alpha\tau Y, \tau\alpha\tau Z )&=(\tau\alpha\tau X, \tau\alpha\tau Y \times \tau\alpha\tau Z)=(\tau\alpha\tau X, \tau((\alpha\tau Y) \times (\alpha\tau Z)))
\\
&=(\tau\alpha\tau X, \tau{}^t\!\alpha^{-1}(\tau Y \times \tau Z))=\tau(\alpha\tau X, {}^t\!\alpha^{-1}(\tau Y \times \tau Z))
\\
&=\tau(\alpha\tau X, \alpha\tau Y \times \alpha\tau Z)=\tau(\alpha\tau X, \alpha\tau Y ,\alpha\tau Z)
\\
&=\tau(\tau X, \tau Y ,\tau Z)=\tau(\tau X,\tau(Y \times Z))
\\
&=\tau^2(X,Y \times Z)
\\
&=(X,Y,Z),
\end{align*}
so we can confirm that $ \tau $ induces an involutive automorphism $ \tilde{\tau} $ of $ (E_{6,\sR})^C $.

\noindent Hence we can define a subgroup $ ((E_{6,\sR})^C)^\tau $ of $ (E_{6,\sR})^C $ by
\begin{align*}
((E_{6,\sR})^C)^\tau:=\left\lbrace \alpha \in (E_{6,\sR})^C \relmiddle{|} \tilde{\tau}(\alpha)=\alpha \right\rbrace.
\end{align*}

Then we have the following proposition.

\begin{proposition}\label{proposition 3.6}
The group $ ((E_{6,\sR})^C)^\tau $ coincides with the group $ E_{6(-26),\sR} ${\rm :} $ ((E_{6,\sR})^C)^\tau= E_{6(-26),\sR} $.
\end{proposition}
\begin{proof}
Let $ \alpha \in ((E_{6,\sR})^C)^\tau $. Then, for $ X \in \mathfrak{J}(3,\R) $, it follows from $ \tau X=X $ that
\begin{align*}
\alpha X=\alpha(\tau X)=\tau(\alpha X).
\end{align*}
Hence we have $ \alpha X \in \mathfrak{J}(3,\R) $, so that $ \alpha $ induces an $ \R $-linear isomorphism of $ \mathfrak{J}(3,\R) $. Thus we see $ \alpha \in E_{6(-26),\sR} $. Conversely, let $ \beta \in E_{6(-26),\sR} $. Then we define an action to $ X \in \mathfrak{J}(3,\R^C) $ of $ \beta $ by
\begin{align*}
\beta X=\beta(X_1+iX_2)=\beta X_1+i\beta X_2,\;\;X:=X_1+iX_2,X_i \in \mathfrak{J}(3,\R).
\end{align*}
Hence it follows that
\begin{align*}
\tau\beta X&=\tau\beta(X_1+iX_2)=\tau(\beta X_1+i\beta X_2)=\beta X_1-i\beta X_2
\\
&=\beta(X_1-iX_2)=\beta\tau(X_1+iX_2)
\\
&=\beta\tau X,
\end{align*}
that is, $ \tau\beta=\beta\tau $. Hence we see $ \beta \in ((E_{6,\sR})^C)^\tau $.

With above, the proof of this proposition is completed.
\end{proof}

We determine the structure of the group $ E_{6(-26),\sR} $.

\begin{theorem}\label{theorem 3.7}
The group $ E_{6(-26),\sR} $ is isomorphic to the group $ SL(3,\R) ${\rm :} $ E_{6(-26),\sR} \cong SL(3,\R) $.
\end{theorem}
\begin{proof}
Let the group $ E_{6(-26),\sR} $ as the group $ ((E_{6,\sR})^C)^\tau $ (Proposition \ref{proposition 3.6}). Then we define a mapping $ \varphi_{6(-26),\sR}: SL(3,\R) \to ((E_{6,\sR})^C)^\tau $ by
\begin{align*}
\varphi_{6(-26),\sR}(A)X=AX\,{}^t\!A,\;\;X \in \mathfrak{J}(3,\R^C).
\end{align*}
Note that this mapping is the restriction of the mapping $ f_{6,C} $ (Theorem \ref{theorem 3.1}).
First, we will prove that $ \varphi_{6(-26),\sR} $ is well-defined. Since $ SL(3,\R) \subset SL(3,C) $, it is easy to see $ \varphi_{6(-26),\sR}(A) \in ((E_{6,\sR})^C)^\tau $. Moreover, from Lemma \ref{lemma 3.4} (3) we have
\begin{align*}
\tau(\varphi_{6(-26),\sR}(A))\tau=\varphi_{6(-26),\sR}(\tau A).
\end{align*}
Hence, since we have $ \tau A=A $ from $ A \in SL(3,\R) $, we obtain $ \tau(\varphi_{6(-26),\sR}(A))\tau=\varphi_{6(-26),\sR}(A) $, so that $ \varphi_{6(-26),\sR} $ is well-defined.
Subsequently, we will prove that $ \varphi_{6(-26),\sR} $ is a homomorphism. It follows that
\begin{align*}
\varphi_{6(-26),\sR}(AB)X=(AB)X\,{}^t\!(AB)=A(BX\,{}^t\!B){}^t\!A=\varphi_{6(-26),\sR}(A)\varphi_{6(-26),\sR}(B)X,
\end{align*}
that is, $ \varphi_{6(-26),\sR}(AB)=\varphi_{6(-26),\sR}(A)\varphi_{6(-26),\sR}(B) $.

Next, we will prove that $\varphi_{6(-26),\sR}$ is surjective. Let $ \alpha \in E_{6(-26),\sR}=((E_{6,\sR})^C)^\tau \subset (E_{6,\sR})^C$ (Proposition \ref{proposition 3.6}). Then there exists $ P \in (E_{6,\sR})^C $ such that $ \alpha=f_{6,C}(P) $. Since $ \alpha $ satisfies the condition $ \tau\alpha\tau=\alpha $, that is, $ \tau f_{6,C}(P)\tau=f_{6,C}(P) $, it follows from Lemma \ref{lemma 3.4} (3) that
\begin{align*}
\tau f_{6,C}(P)\tau=f_{6,C}(\tau P).
\end{align*}
Hence we have the following
\begin{align*}
\tau P =P,
\end{align*}
that is, $ P \in SL(3,\R) $. Thus there exists $ A \in SL(3,\R) $ such that $ \alpha=f_{6,C}(A)=\varphi_{6(-26),\sR}(A) $. With above, the proof of surjective is completed.

Finally, we will determine $ \Ker\,\varphi_{6(-26),\sR} $. since $ \varphi_{6(-26),\sR} $ is the restriction of the mapping $ f_{6,C} $, it is easy to obtain $ \Ker\,\varphi_{6(-26),\sR}=\{E\} $.

Therefore we have the required isomorphism
\begin{align*}
 E_{6(-26),\sR} \cong SL(3,\R).
\end{align*}
\end{proof}

\section{The complex Lie group $ (E_{6,\sC})^C $ and its real forms}

We define the group $ (E_{6,\sC})^C $ by
\begin{align*}
(E_{6,\sC})^C:&=\left\lbrace \alpha \in \Iso_{C}(\mathfrak{J}(3,\C^C))\relmiddle{|} \det(\alpha X)=\det\,X \right\rbrace
\\
&=\left\lbrace \alpha \in \Iso_{C}(\mathfrak{J}(3,\C^C))\relmiddle{|} (\alpha X,\alpha Y, \alpha Z)=(X,Y,Z) \right\rbrace
\\
&=\left\lbrace \alpha \in \Iso_{C}(\mathfrak{J}(3,\C^C))\relmiddle{|} \alpha X \times \alpha Y={}^t\!\alpha^{-1}(X \times Y) \right\rbrace.
\end{align*}
and its real forms are defined as follows:
\begin{align*}
E_{6,\sC}&:=\left\lbrace \alpha \in \Iso_{C}(\mathfrak{J}(3,\C^C))\relmiddle{|} \det(\alpha X)=\det\,X, \langle \alpha X, \alpha Y \rangle=\langle X,Y \rangle \right\rbrace,
\\
E_{6(6),\sC'}&:=\left\lbrace \alpha \in \Iso_{\sR}(\mathfrak{J}(3,\C'))\relmiddle{|} \det(\alpha X)=\det\,X \right\rbrace,
\\
E_{6(-14),\sC}&:=\left\lbrace \alpha \in \Iso_{C}(\mathfrak{J}(3,\C^C))\relmiddle{|} \det(\alpha X)=\det\,X, \langle \alpha X, \alpha Y \rangle_\sigma=\langle X,Y \rangle_\sigma \right\rbrace,
\\
E_{6(-26),\sC}&:=\left\lbrace \alpha \in \Iso_{\sR}(\mathfrak{J}(3,\C))\relmiddle{|} \det(\alpha X)=\det\,X \right\rbrace,
\end{align*}
where $ \langle X,Y \rangle=(\tau X,Y), \langle X,Y \rangle_\sigma=(\tau\sigma X,Y) $,
and since $ E_{6(2),\sC}=E_{6,\sC}$, the definition of $ E_{6(2),\sC} $ is omitted.
\vspace{1mm}

Let $ \alpha \in (E_{6,\sC})^C  $. Then, as in \cite[Lemma 3.2.1]{iy7}, we have $ {}^t\!\alpha^{-1} \in (E_{6,\sC})^C $, where the transpose $ {}^t\!\alpha $ of $ \alpha $ is defined by $ ({}^t\!\alpha X,Y)=(X,\alpha Y), X,Y \in \mathfrak{J}(3,\C^C) $. Hence we can define an involutive automorphism $ \lambda $ of $ (E_{6,\sC})^C $ by
\begin{align*}
\lambda(\alpha):={}^t\!\alpha^{-1}, \alpha \in (E_{6,\sC})^C.
\end{align*}
Here, we define an $ \R $-linear transformation $ \gamma_{\scalebox{0.8}{\sC}} $ of $ \C $ by
\begin{align*}
\gamma_{\scalebox{0.8}{\sC}}(x)=\ov{x}.
\end{align*}
Needless to say, $ \gamma_{\scalebox{0.8}{\sC}} $ is also the complex conjugation of $ \C $. Then $ \gamma_{\scalebox{0.8}{\sC}} $ is naturally extended to the $ C $-linear transformation of $ \C^C $, so that $ \gamma_{\scalebox{0.8}{\sC}} $ can be extended to the $ C $-linear transformation of $ \mathfrak{J}(3,\C^C) $ by
\begin{align*}
\gamma_{\scalebox{0.8}{\sC}}X=\begin{pmatrix}
\xi_1 & \gamma_{\scalebox{0.8}{\sC}}x_3 & \ov{\gamma_{\scalebox{0.8}{\sC}}x}_2 \\
\ov{\gamma_{\scalebox{0.8}{\sC}}x}_3 & \xi_2 & \gamma_{\scalebox{0.8}{\sC}}x_1 \\
\gamma_{\scalebox{0.8}{\sC}}x_2 & \ov{\gamma_{\scalebox{0.8}{\sC}}x}_1 & \xi_3
\end{pmatrix},\;\; X \in \mathfrak{J}(3,\C^C)
\end{align*}
with the properties of $ \gamma_{\scalebox{0.8}{\sC}} \in G_{2,\sC} \subset (G_{2,\sC})^C \subset (F_{4,\sC})^C \subset (E_{6,\sC})^C $ and $ {\gamma_{\scalebox{0.8}{\sC}}}^2=1 $.
Moreover, let $ \tau $ be the complex conjugation in $ \mathfrak{J}(3,\C^C) $ and $ \sigma $ be the $ C $-linear transformation of $ \mathfrak{J}(3,\C^C) $.
Then $ (E_{6,\sC})^C $ has involutive automorphism $ \tilde{\tau\lambda}, ,\tilde{\tau\gamma_{\scalebox{0.8}{\sC}}}, \tilde{\tau\lambda\sigma} $ and $ \tilde{\tau} $. Indeed, for $ \alpha \in (E_{6,\sC})^C $, we have to show
\begin{align*}
\tilde{\tau\lambda}(\alpha):=\tau\lambda(\alpha)\tau,
\tilde{\tau\gamma_{\scalebox{0.8}{\sC}}}(\alpha):=(\tau\gamma_{\scalebox{0.8}{\sC}})\alpha(\gamma_{\scalebox{0.8}{\sC}}\tau),
\tilde{\tau\lambda\sigma}(\alpha):=(\tau\sigma)\lambda(\alpha)(\sigma\tau), \tilde{\tau}(\alpha):=\tau\alpha\tau \in (E_{6,\sC})^C.
\end{align*}
First, we can immediately confirm $ \tau\alpha\tau \in \Iso_C(\mathfrak{J}(3,\C^C)) $ for $ \alpha \in (E_{6,\sC})^C $, so that the transpose $ {}^t(\tau\alpha\tau) $ of $ \tau\alpha\tau $ is defined by $ ({}^t(\tau\alpha\tau)X,Y)=(X,(\tau\alpha\tau) Y) $. Hence, note that $ (\tau X,\tau Y)=\tau(X,Y) $ and $ \tau^2=1 $, it follows from
\begin{align*}
(\tau{}^t\!\alpha\tau X, Y)&=(\tau{}^t\!\alpha\tau X, \tau^2 Y)=\tau({}^t\!\alpha\tau X, \tau Y)=\tau(\tau X,\alpha \tau Y )
\\
&=\tau(\tau X,\tau^2\alpha \tau Y )=\tau^2(X,\tau\alpha \tau Y )
=(X,\tau\alpha \tau Y ), X, Y \in \mathfrak{J}(3,\C^C)
\end{align*}
that $ {}^t(\tau\alpha\tau)X=\tau{}^t\!\alpha\tau X $, that is, $ {}^t(\tau\alpha\tau)=\tau{}^t\!\alpha\tau $, so that we have $ {}^t(\tau\alpha\tau)^{-1}=\tau{}^t\!\alpha^{-1}\tau $.

Thus, using $ {}^t(\tau\alpha\tau)^{-1}=\tau{}^t\!\alpha^{-1}\tau $ above, it follows from $ \tau X \times \tau Y=\tau(X \times Y) $ that
\begin{align*}
\tilde{\tau\lambda}(\alpha)X \times \tilde{\tau\lambda}(\alpha)Y&=(\tau\lambda(\alpha)\tau)X \times (\tau\lambda(\alpha)\tau)Y=(\tau{}^t\!\alpha^{-1}\tau)X \times (\tau{}^t\!\alpha^{-1}\tau)Y
\\
&=\tau(({}^t\!\alpha^{-1}\tau)X \times ({}^t\!\alpha^{-1}\tau)Y)=\tau\alpha(\tau X \times \tau Y)
\\
&=\tau\alpha\tau(X \times Y)={}^t\!(\tau{}^t\!\alpha^{-1}\tau)^{-1}(X \times Y)
\\
&={}^t\!(\tilde{\tau\lambda}(\alpha))^{-1}(X \times Y).
\end{align*}
Hence we see $ \tilde{\tau\lambda}(\alpha) \in (E_{6,\sC})^C $.

Next, by the similar computation as above, it follows from $ \gamma_{\scalebox{0.8}{\sC}} X \times \gamma_{\scalebox{0.8}{\sC}} Y=\gamma_{\scalebox{0.8}{\sC}}(X \times Y) $ that
\begin{align*}
\tilde{\tau\gamma_{\scalebox{0.8}{\sC}}}(\alpha)X \times \tilde{\tau\gamma_{\scalebox{0.8}{\sC}}}(\alpha)Y
&=(\tau\gamma_{\scalebox{0.8}{\sC}})\alpha(\gamma_{\scalebox{0.8}{\sC}}\tau)X \times (\tau\gamma_{\scalebox{0.8}{\sC}})\alpha(\gamma_{\scalebox{0.8}{\sC}}\tau)Y
=\tau\gamma_{\scalebox{0.8}{\sC}}(\alpha(\gamma_{\scalebox{0.8}{\sC}}\tau)X \times \alpha(\gamma_{\scalebox{0.8}{\sC}}\tau)Y)
\\
&=\tau\gamma_{\scalebox{0.8}{\sC}}\,{}^t\!\alpha^{-1}((\gamma_{\scalebox{0.8}{\sC}}\tau)X \times (\gamma_{\scalebox{0.8}{\sC}}\tau)Y)
=(\tau\gamma_{\scalebox{0.8}{\sC}}){}^t\!\alpha^{-1}(\gamma_{\scalebox{0.8}{\sC}}\tau)(X \times Y)
\\
&={}^t((\tau\gamma_{\scalebox{0.8}{\sC}})\alpha(\gamma_{\scalebox{0.8}{\sC}}\tau))^{-1}(X \times Y)
\\
&={}^t(\tilde{\tau\gamma_{\scalebox{0.8}{\sC}}}(\alpha))^{-1}(X \times Y).
\end{align*}
Hence we see $ \tilde{\tau\gamma_{\scalebox{0.8}{\sC}}}(\alpha) \in (E_{6,\sC})^C$.

Moreover, it follows from $ \sigma X \times \sigma Y=\sigma(X \times Y) $ that
\begin{align*}
\tilde{\tau\lambda\sigma}X \times \tilde{\tau\lambda\sigma}Y&=((\tau\sigma){}^t\!\alpha^{-1}(\sigma\tau))X \times ((\tau\sigma){}^t\!\alpha^{-1}(\sigma\tau))Y
\\
&=\tau((\sigma{}^t\!\alpha^{-1}(\sigma\tau))X \times (\sigma{}^t\!\alpha^{-1}(\sigma\tau))Y)
\\
&=\tau\sigma(({}^t\!\alpha^{-1}(\sigma\tau))X \times ({}^t\!\alpha^{-1}(\sigma\tau))Y)
\\
&=(\tau\sigma)\alpha((\sigma\tau)X \times ((\sigma\tau)Y)
\\
&=(\tau\sigma)\alpha(\sigma\tau)(X \times Y)
\\
&={}^t((\tau\sigma){}^t\!\alpha^{-1}(\sigma\tau))^{-1}(X \times Y)
\\
&={}^t\!(\tilde{\tau\lambda\sigma}(\alpha))^{-1}(X \times Y).
\end{align*}
Hence we also see $ \tilde{\tau\lambda\sigma}(\alpha) \in (E_{6,\sC})^C$.
Finally, since it is easy to verify that $ \tilde{\tau}(\alpha) \in (E_{6,\sC})^C $, its proof is omitted.

\subsection{The group $ (E_{6,\sC})^C $}

First, we prove the lemmas used in order to prove theorem below.

\begin{lemma}\label{lemma 4.1.1}
Any element $ X \in \mathfrak{J}(3,\C^C) $ such that $ X^2=X, \tr(X)=1 $ can be transformed to any $ E_i,\allowbreak i=1,2,3 $ by some $ B \in U(3,\C^C) ${\rm :}
$ B^*XB=E_i $.
\end{lemma}
\begin{proof}
It is well-known that any element $ X \in \mathfrak{J}(3,\C) $ can be transformed to diagonal form by a certain $ B \in U(3) $. Hence, since $ U(3,\C^C) $ contains the subgroup $ U(3) $, we may assume $ X \in \mathfrak{J}(3,\C^C) $ as
\begin{align*}
X=\begin{pmatrix}
\xi_1 & ix_3 & i\ov{x}_2 \\
i\ov{x}_3 & \xi_2 & ix_1  \\
ix_2 & i\ov{x}_1 & \xi_3
\end{pmatrix},\;\;\begin{array}{l}
\xi_i \in C, \xi_1+\xi_2+\xi_3=1,\\
x_i \in \C.
\end{array}
\end{align*}
Subsequently, the computation of $ X^2 $ is obtained as follows:
\begin{align*}
X^2=\begin{pmatrix}
{\xi_1}^2-x_2\ov{x}_2-x_3\ov{x}_3 & -\ov{x}_2\ov{x}_1+i(\xi_1+\xi_2)x_3 &  *    \\
* & {\xi_2}^2-x_3\ov{x}_3-x_1\ov{x}_1 & -\ov{x}_3\ov{x}_2+i(\xi_2+\xi_3)x_1 \\
-\ov{x}_1\ov{x}_3+i(\xi_3+\xi_1)x_2 &  *  & {\xi_3}^2-x_1\ov{x}_1-x_2\ov{x}_2
\end{pmatrix} \in \mathfrak{J}(3,\C^C).
\end{align*}
Then we compare the diagonals of both of $ X^2=X $, so that we have that each $ \xi_i $ is real numbers. Indeed, for instance, since $ {\xi_1}^2-x_2\ov{x}_2-x_3\ov{x}_3=\xi_1 $, that is, $  {\xi_1}^2-\xi_1-x_2\ov{x}_2-x_3\ov{x}_3=0 $, we have $ \xi_1=(1/2)(1\pm\sqrt{1+(x_2\ov{x}_2+x_3\ov{x}_3)}) \in \R$, so are $ \xi_2,\xi_3 $. Hence by comparing $ F_k(x_k),k=1,2,3 $-part, we have
\begin{align*}
x_1x_2=x_2x_3=x_3x_1=0, \;\; \xi_1x_1=\xi_2x_2=\xi_3x_3=0 \;\cdots \,(*).
\end{align*}

In the case where $ x_1=x_2=x_3=0 $. Since $ X $ is diagonal form, we have $ \xi_i=0 $ or $ \xi_i=1 $ from $ X^2=X $, and together with $ \tr(X)=1 $, we see $ X=E_1, X=E_2 $ or $ X=E_3 $. Set $ C_2:=\begin{pmatrix}
0 & 1 & 0 \\
1 & 0 & 0 \\
0 & 0 & 1
\end{pmatrix}, C_3:=\begin{pmatrix}
1 & 0 & 0 \\
0  & 0 & 1 \\
0  & 1 & 0
\end{pmatrix} $, so we easily see $ C_k \in O(3) \subset U(3) \subset U(3,\C^C),\vspace{1mm} k=2,3 $. Then it follows from
\begin{align*}
{}^t\!C_2E_1C_2=E_2,\;\; {}^t\!(C_3C_2)E_1(C_3C_2)=E_3, \;\; {}^t\!C_3E_2C_3=E_3
\end{align*}
that this lemma is valid.

In the case where $ x_1\not=0 $. Then we have $ x_2=x_3=0, \xi_1=0, \xi_2+\xi_3=1 $ from the formulas $ (*) $ above. Hence $ X $ is of the form $ \begin{pmatrix}
0 & 0 & 0 \\
0 & \xi_2 & ix_1 \\
0 & i\ov{x}_1 & \xi_3
\end{pmatrix}$ with $ x_1\ov{x}_1=-\xi_2\xi_3 $. Here, note that $ \xi_2\xi_3 <0 $, if $ \xi_2 >0 $ and $ \xi_3 <0 $, $ X $ can be transformed to $ E_2 $ by $ B_1:=\begin{pmatrix}
1 & 0 & 0 \\
0 & x_1/\sqrt{-\xi_3} & -ix_1/\sqrt{\xi_2} \\
0 & i\sqrt{-\xi_3} & \sqrt{\xi_2}
\end{pmatrix} \in U(3,\C^C)$. Indeed, first it follows  from
\begin{align*}
B_1{B_1}^*=\begin{pmatrix}
1 & 0 & 0 \\
0 & x_1/\sqrt{-\xi_3} & -ix_1/\sqrt{\xi_2} \\
0 & i\sqrt{-\xi_3} & \sqrt{\xi_2}
\end{pmatrix}
\begin{pmatrix}
1 & 0 & 0 \\
0 & \ov{x}_1/\sqrt{-\xi_3} & i\sqrt{-\xi_3} \\
0 & -i\ov{x}_1/\sqrt{\xi_2} & \sqrt{\xi_2}
\end{pmatrix}
=E
\end{align*}
that $ B_1 \in U(3,\C^C) $. By straightforward computation, we have $ {B_1}^*XB_1=E_2 $.
If $ \xi_2 <0 $ and $ \xi_3 >0 $, $ X $ can be also transformed to $ E_2 $ by $ B_2:=\begin{pmatrix}
1 & 0 & 0 \\
0 & i\sqrt{-\xi_2} & \sqrt{\xi_3} \\
0 & \ov{x}_1/\sqrt{-\xi_2} & -i\ov{x}_1/\sqrt{\xi_3}
\end{pmatrix} \in U(3,\C^C) $. Indeed, it follows from
\begin{align*}
B_2{B_2}^*=\begin{pmatrix}
1 & 0 & 0 \\
0 & i\sqrt{-\xi_2} & \sqrt{\xi_3} \\
0 & \ov{x}_1/\sqrt{-\xi_2} & -i\ov{x}_1/\sqrt{\xi_3}
\end{pmatrix}\begin{pmatrix}
1 & 0 & 0 \\
0 & i\sqrt{-\xi_2} & x_1/\sqrt{-\xi_2} \\
0 & \sqrt{\xi_3} & -ix_1/\sqrt{\xi_3}
\end{pmatrix}=E
\end{align*}
that $ B_2 \in U(3,\C^C) $. As in the case above, we have $ {B_2}^*XB_2=E_2 $.

 In the case where $ x_2\not=0 $. As in the case where $ x_1\not=0 $, $ X $ is of the form $ \begin{pmatrix}
\xi_1 & 0 & i\ov{x}_2 \\
0 & 0 & 0 \\
ix_2 & 0 & \xi_3
\end{pmatrix} $ with $ x_2\ov{x}_2=-\xi_3\xi_1, \xi_3+\xi_1=1 $. Let $ C_2 $. Then we have
\begin{align*}
C_2X{}^t\!C_2=\begin{pmatrix}
0 & 0 & 0  \\
0 & \xi_1 & i\ov{x}_2 \\
0 & ix_2 & \xi_3
\end{pmatrix}.
\end{align*}
Hence this case is reduced to the case where $ x_1\not=0 $.

In the case where $ x_3\not=0 $. As in the case where $ x_1\not=0 $, $ X $ is of the form $ \begin{pmatrix}
\xi_1 & ix_3 & 0 \\
i\ov{x}_3 & \xi_2 & 0 \\
0 & 0 & 0
\end{pmatrix} $ with $ x_3\ov{x}_3=-\xi_1\xi_2, \xi_1+\xi_2=1 $. Let $ C_3 $ Then we have
\begin{align*}
C_3X{}^t\!C_3=\begin{pmatrix}
\xi_1 & 0 & ix_3  \\
0 & 0 & 0 \\
i\ov{x}_3 & 0 & \xi_2
\end{pmatrix}.
\end{align*}
Hence this case is also reduced to the case where $ x_2\not=0 $.

Finally $ E_1, E_2 $ and $ E_3 $ can be transformed to one another by $ C_2,C_3 $. With above, the proof of this lemma is completed.
\end{proof}

Here, we consider the group $ (F_{4,\sC})^C $ as complexification of the group $ F_{4,\sC} $:
\begin{align*}
(F_{4,\sC})^C:&=\left\lbrace \alpha \in \Iso_C(\mathfrak{J}(3,\C^C))\relmiddle{|} \alpha(X \circ Y)=\alpha X \circ \alpha Y\right\rbrace
\\
&=\left\lbrace \alpha \in \Iso_C(\mathfrak{J}(3,\C^C))\relmiddle{|} \alpha(X \times Y)=\alpha X \times \alpha Y\right\rbrace
\\
&=\left\lbrace \alpha \in \Iso_{C}(\mathfrak{J}(3,\C^C))\relmiddle{|} \det(\alpha X)=\det\,X, \alpha E=E \right\rbrace(=((E_{6,\sC})^C)_E),
\end{align*}
where as for $ F_{4,\sC} $, see \cite[Theorem 5]{iy9} in detail.

\begin{lemma}\label{lemma 4.1.2}
For $ \alpha \in (F_{4,\sC})^C $, $ \alpha $ satisfies $ \alpha E=E $.
\end{lemma}
\begin{proof}
Let $ \alpha \in (F_{4,\sC})^C $. Apply on the both of $ E \circ X=X, X \in \mathfrak{J}(3,\C^C) $, then we have $ \alpha E \circ \alpha X=\alpha X $. Here, set $ X:=\alpha^{-1} E $, so that $ \alpha E \circ E=E $, that is, $ \alpha E=E $.
\end{proof}

Using the lemma above, we prove the following theorem needed later by an argument similar to that in the proof of \cite[Theorem 5]{iy9}.

We consider a discrete group $ \Z_2:=\{1,\varepsilon\} $, where $ \varepsilon $ is the complex conjugation of $ \C $: $ \varepsilon x=\ov{x}, x \in \C $, and $ \varepsilon $ is naturally extended to the mapping $ \varepsilon: \C^C \to \C^C $.
Then this group acts on the group $ SU(3,\C^C) $ by
\begin{align*}
1A=A, \quad \varepsilon A=\ov{A}
\end{align*}
and let $ SU(3,\C^C) \rtimes \Z_2 $ be the semi-direct product of $ SU(3,\C^C) $ and $ \Z_2 $ with the multiplication
\begin{align*}
(A_1,1)(A_2,1)&=(A_1A_2,1),\;(A_1,1)(A_2,\varepsilon)=(A_1A_2,\varepsilon),
\\
(A_1,\varepsilon)(A_2,1)&=(A_1\ov{A}_2, \varepsilon), \; (A_1,\varepsilon)(A_2,\varepsilon)=(A_1\ov{A}_2, 1).
\end{align*}

\begin{theorem}\label{theorem 4.1.3}
The group $ (F_{4,\sC})^C $ is isomorphic to the semi-direct product of the group $ SL(3,C)/\Z_3 $ and the discrete group $ \Z_2, \Z_3=\{E,\omega E, \omega^2E \}, \Z_2=\{1,\varepsilon \}${\rm :} $ (F_{4,\sC})^C \cong SL(3,C)/\Z_3 \rtimes \Z_2 $, where $ \omega \in C, \allowbreak \omega^3=1 $.
\end{theorem}
\begin{proof}
Let $ SL(3,C) $ as the group $ SU(3,\C^C) $. We define a mapping $ f_{4,\sC^C}: SU(3,\C^C) \rtimes \Z_2 \to (F_{4,\sC})^C $ by
\begin{align*}
f_{4,\sC^C}(A,1)X&=AX A^*,
\\
f_{4,\sC^C}(A,\varepsilon)X&=A\ov{X} A^*,\;\; X \in \mathfrak{J}(3,\C^C).
\end{align*}
First, we will prove that $ f_{4,\sC^C} $ is well-defined. It follows that
\begin{align*}
f_{4,\sC^C}(A,1)X \circ f_{4,\sC^C}(A,1)Y&=(AX A^*) \circ (A Y A^*)
\\
&=\dfrac{1}{2}((AX A^*)(AY A^*)+(AY A^*)(AX A^*))
\\
&=\dfrac{1}{2}((A XY A^*)+(A YX A^*))
\\
&=A(\dfrac{1}{2}(XY+YX))A^*=A(X \circ Y)A^*
\\
&=f_{4,\sC^C}(A,1)(X \circ Y),\;\; X,Y \in \mathfrak{J}(3,\C^C),
\end{align*}
so that $ f_{4,\sC^C}(A,1) \in (F_{4,\sC})^C$, and we have $ f_{4,\sC^C}(A,\varepsilon) \in (F_{4,\sC})^C$ in exactly the same way. 
It is clear that $ f_{4,\sC^C} $ is a homomorphism. Indeed, it follows that
\begin{align*}
f_{4,\sC^C}(A_1,1)f_{4,\sC^C}(A_2,1)X&=f_{4,\sC^C}(A_1,1)(A_2X{A_2}^*)=A_1(A_2X{A_2}^*){A_1}^*
\\
&=(A_1A_2)X(A_1A_2)^*=f_{4,\sC^C}(A_1A_2,1)X
\\
&=f_{4,\sC^C}((A_1,1)(A_2,1))X,
\\[1mm]
f_{4,\sC^C}(A_1,1)f_{4,\sC^C}(A_2,\varepsilon)X&=f_{4,\sC^C}(A_1,1)(A_2\ov{X}{A_2}^*)=A_1(A_2\ov{X}{A_2}^*){A_1}^*
\\
&=(A_1A_2)\ov{X}(A_1A_2)^*=f_{4,\sC^C}(A_1A_2,\varepsilon)X
\\
&=f_{4,\sC^C}((A_1,1)(A_2,\varepsilon))X,
\\[1mm]
f_{4,\sC^C}(A_1,\varepsilon)f_{4,\sC^C}(A_2,1)X&=f_{4,\sC^C}(A_1,\varepsilon)(A_2X{A_2}^*)=A_1(\ov{A_2X{A_2}^*}){A_1}^*
\\
&=(A_1\ov{A}_2)\ov{X}(A_1\ov{A}_2)^*=f_{4,\sC^C}(A_1\ov{A}_2,\varepsilon)X
\\
&=f_{4,\sC^C}((A_1,\varepsilon)(A_2,1))X,
\\[1mm]
f_{4,\sC^C}(A_1,\varepsilon)f_{4,\sC^C}(A_2,\varepsilon)X&=f_{4,\sC^C}(A_1,\varepsilon)(A_2\ov{X}{A_2}^*)=A_1(\ov{A_2\ov{X}{A_2}^*}){A_1}^*
\\
&=(A_1\ov{A}_2)X(A_1\ov{A}_2)^*=f_{4,\sC^C}(A_1\ov{A}_2,1)X
\\
&=f_{4,\sC^C}((A_1,\varepsilon)(A_2,\varepsilon))X.
\end{align*}

Next we will prove that $ f_{4,\sC^C} $ is surjective. Let $ \alpha \in (F_{4,\sC})^C $. We consider the elements $ \alpha E_i \in \mathfrak{J}(3,\C^C), i=1,2,3 $. Then, since the formulas $ E_i \circ E_i=E_i $ and $ \tr(\alpha E_i)=\tr\,E_i $ hold, $ \alpha E_i $ satisfies the conditions
\begin{align*}
(\alpha E_i)^*=\alpha E_i,\;\; (\alpha E_i)^2=\alpha E_i,\;\;\tr(\alpha E_i)=1.
\end{align*}
Hence, from Lemma \ref{lemma 4.1.1}, there exists $ C_i \in U(3,\C^C) $ such that $ C_i(\alpha E_i){C_i}^*=E_i, i=1,2,3 $, so that let
$ \c_i:=(c_{i1}\; c_{i2}\; c_{i3}) $ as the $ i $-th row of $ C_i,i=1,2,3 $.
We construct a matrix $ C:=
\begin{pmatrix}
\c_1 \\
\c_2 \\
\c_3
\end{pmatrix}$, then we have  $ \alpha E_i=C^*E_iC , i=1,2,3 $ and $ C \in U(3,\C^C) $. Indeed, it is clear that $ \alpha E_i=C^*E_iC , i=1,2,3 $ hold, and using Lemma \ref{lemma 4.1.2}, it follows from
\begin{align*}
C^*C&=C^*EC=C^*(E_1+E_2+E_3)C=C^*E_1C+C^*E_2C+C^*E_3C
\\
&=\alpha E_1+\alpha E_2+\alpha E_3=\alpha E
\\
&=E
\end{align*}
that $ C \in U(3,\C^C) $. Hereafter, we may assume $ C \in SU(3,\C^C) $, if necessary, replace $ \a_1 $ with $ \a_1/\det\,C $. Now, set $ \beta:={f_{4,\sC^C}(C,1)}^{-1}\alpha $, then $ \beta $ satisfies $ \beta \in (F_{4,\sC})^C $ and $ \beta E_i=E_i, i=1,2,3 $. indeed, it follows that $ \beta E_i={f_{4,\sC^C}(C,1)}^{-1}\alpha E_i={f_{4,\sC^C}(C,1)}^{-1}(C^*E_iC)=C(C^*E_iC)C^*=E_i $. Here, set
\begin{align*}
\mathfrak{J}_i:&=\left\lbrace X \in \mathfrak{J}(3,\C^C) \relmiddle{|}2E_{i+1} \circ X=2E_{i+2} \circ X=X \right\rbrace
\\
&=\left\lbrace F_i(x) \relmiddle{|} x \in \C^C \right\rbrace,i=1,2,3,
\end{align*}
then we have $ \beta X \in \mathfrak{J}_i $ for $ X \in \mathfrak{J}_i $, so that $ \beta $ induces $ C $-linear transformations of $ \C^C $ such that
\begin{align*}
\beta F_i(x)=F_i(\beta_i(x)),i=1,2,3,
\end{align*}
in addition, apply $ \beta $ on the both of $ F_i(x) \circ F_j(y)=(x,y)(E_{i+1}E_{i+2}) $, then we have $ (\beta_i x, \beta_i y)=(x,y) $, that is, $ \beta \in O(2,\C^C) $. Indeed, it follows from
\begin{align*}
\beta (F_i(x) \circ F_j(y))&=\beta F_i(x) \circ \beta F_j(y)=F_i(\beta_i(x)) \circ F_j(\beta_i(y))=(\beta_i x, \beta_i y)(E_{i+1}+E_{i+2}),
\\[1mm]
\beta(x,y)(E_{i+1}+E_{i+2})&=(x,y)\beta(E_{i+1}+E_{i+2})=(x,y)(E_{i+1}+E_{i+2})
\end{align*}
that $ (\beta_i (x), \beta_i (y))=(x,y) $, so that $ \beta_i \in O(2,\C^C) $.
Moreover, $ \beta_1, \beta_2 $ and $ \beta_3 $ are combined with the following relation
\begin{align*}
\beta_1(x)\beta_2(y)=\ov{\beta_3(\ov{xy})}, \;\; x,y \in \C^C.
\end{align*}
Indeed, apply $ \beta $ on the both of $ 2F_1(x) \circ F_2(y)=F_3(\ov{xy}) $, then it follows from
\begin{align*}
\beta(2F_1(x) \circ F_2(y))&=2\beta F_1(x) \circ \beta F_2(y)=2F_1(\beta_1(x)) \circ F_2(\beta_2(y))=F_3(\ov{\beta_1(x)\beta_2(y)}),
\\[1mm]
\beta F_3(\ov{xy})&=F_3(\beta_3(\ov{xy}))
\end{align*}
that $ \ov{\beta_1(x)\beta_2(y)}=\beta_3(\ov{xy}) $, that is, $ \beta_1(x)\beta_2(y)=\ov{\beta_3(\ov{xy})} $.

We will investigate the relation of $ \beta_1, \beta_2 $ and $ \beta_3 $ in more detail. Set $ p:=\beta_1(1) \in \C^C $ and $ q:=\beta_2(1) \in \C^C $. Then we see $ |p|=|q|=1 $ from $ (\beta_i (x), \beta_i (y))=(x,y) $, moreover using $ \beta_1(x)\beta_2(y)=\ov{\beta_3(\ov{xy})} $, set $ y=1 $, then we have $ \beta_1(x)\beta_2(1)=\ov{\beta_3(\ov{x})} $, that is, $ \beta_1(x)q=\ov{\beta_3(\ov{x})} $, and set $ x=1 $, then we have $ \beta_1(1)\beta_2(y)=\ov{\beta_3(\ov{y})} $, that is, $ p\beta_2(x)=\ov{\beta_3(\ov{x})} $. Hence we see $ \beta_2(x)=\ov{p}\beta_1(x)q $. Again, using $ \beta_1(x)\beta_2(y)=\ov{\beta_3(\ov{xy})} $, set $ y=1 $, then we have $ \beta_3(\ov{x})=\ov{\beta_1(x)q} $, that is, $ \beta_3(x)=\ov{\beta_1(\ov{x})q} $. With above, we rewrite the results obtained as follows:
\begin{align*}
|p|=|q|=1, \quad \beta_2(x)=\ov{p}\beta_1(x)q, \quad \beta_3(x)=\ov{\beta_1(\ov{x})q}\; \cdots\;(*).
\end{align*}
Furthermore set $ \beta_1(x)=p\sigma(x) $, where $ \sigma $ is a $ C $-linear transformation of $ \C^C $ with $ \sigma(1)=1 $. Then $ \sigma $ satisfies $ \sigma(xy)=\sigma(x)\sigma(y), x,y \in \C^C $, that is, $ \sigma $ is an automorphism of $ \C^C $. Indeed, it follows from the formulas $ (*) $ that $ \beta_2(x)=\sigma(x)q $ and $ \beta_3(x)=\ov{p\sigma(x)q} $, so that $ \sigma(xy)=\sigma(x)\sigma(y) $.
For $ \sigma \in \Aut(\C^C) $, since $ \sigma(\i)\sigma(\i)=\sigma(\i\i)=\sigma(-1)=-\sigma(1)=-1 $, that is, $ (\sigma(\i))^2=-1 $. Hence we have the following
\begin{align*}
{\rm (i)}\; \sigma(\i)=\i,\quad {\rm (ii)}\; \sigma(\i)=-\i,\quad {\rm (iii)}\,\sigma(\i)=i,\quad {\rm (iv)}\,\sigma(\i)=-i,
\end{align*}
where $\i \in C, i \in C $. In the case (i), we see $ \sigma=1 $. In the case (ii), we see $ \sigma=\varepsilon $. In the case (iii), apply $ \sigma $ on the both of $ \sigma(\i)=i $, then we have $ \i=i $ because $ \sigma $ is the $ C $-linear isomorphism of $ \C^C $. This is contradiction, so that this case is impossible. In the case (iv), as in the case (iii), this case is also impossible. Thus we see
\begin{align*}
\sigma(x)=x,\quad \sigma(x)=\ov{x}, \;\; x \in \C^C.
\end{align*}
Therefore we have
\begin{align*}
\beta_1(x)=px,\; \beta_2(x)=xq, \; \beta_3(x)=\ov{q}x\ov{p}\;\;\;\text{or}\;\;\;\beta_1(x)=p\ov{x},\; \beta_2(x)=\ov{x}q, \; \beta_3(x)=\ov{qxp}.
\end{align*}
Here, we choose $ r \in \C^C $ such that $ r^3=\ov{p}q $ and construct a matrix $ D:=\begin{pmatrix}
q\ov{r} & 0 & 0 \\
0 & \ov{pr} & 0 \\
0 & 0 & \ov{r}
\end{pmatrix} $. 
Then we have $ D \in SU(3, \C^C) $.
In the former case, since we have
\begin{align*}
\beta X&=\begin{pmatrix}
\xi_1 & \beta_3(x_3) & \ov{\beta_2(x_2)} \\
\ov{\beta_2(x_3)} & \xi_2 & \beta_1(x_1) \\
\beta_2(x_2) & \ov{\beta_1(x_1)} & \xi_3
\end{pmatrix}=
\begin{pmatrix}
\xi_1 & \ov{q}x_3\ov{p} & \ov{q}\ov{x}_2 \\
p\ov{x}_3q & \xi_2 & px_1 \\
x_2q & \ov{x}_1\ov{p} & \xi_3
\end{pmatrix},\;\; X \in \mathfrak{J}(3,\C^C),
\end{align*}
we see $ \beta X=D^*XD $. In the latter case, as in the previous case, we have
\begin{align*}
\beta X&=\begin{pmatrix}
\xi_1 & \beta_3(x_3) & \ov{\beta_2(x_2)} \\
\ov{\beta_2(x_3)} & \xi_2 & \beta_1(x_1) \\
\beta_2(x_2) & \ov{\beta_1(x_1)} & \xi_3
\end{pmatrix}=
\begin{pmatrix}
\xi_1 & \ov{qx_3p} & \ov{q}x_2 \\
px_3q & \xi_2 & p\ov{x}_1 \\
\ov{x}_2q & x_1\ov{p} & \xi_3
\end{pmatrix},\;\; X \in \mathfrak{J}(3,\C^C),
\end{align*}
we see $ \beta X=D^*\ov{X}D $.

Hence, from $ \beta={f_{4,\sC^C}(C,1)}^{-1}\alpha $, we have the following
\begin{align*}
&\alpha X=f_{4,\sC^C}(C,1)\beta X=C^*(D^*XD)C=(DC)^*X(DC)=f_{4,\sC^C}(DC,1) X
\\
& \hspace{50mm} \text{or}
\\
&\alpha X=f_{4,\sC^C}(C,1)\beta X=C^*(D^*\ov{X}D)C=(DC)^*\ov{X}(DC)=f_{4,\sC^C}(DC,\varepsilon) X.
\end{align*}
Hence, since $ DC \in SU(3,\C^C) $, $ f_{4,\sC^C} $ is surjective.

Finally, we will determine $ \Ker\,f_{4,\sC^C} $. It follows from the definition of kernel that
\begin{align*}
\Ker\,f_{4,\sC^C}&=\left\lbrace (A,1) \in SU(3,\C^C) \rtimes \Z_2 \relmiddle{|} f_{4,\sC^C}(A,1)=1 \right\rbrace
\\
& \cup \left\lbrace (A,\varepsilon) \in SU(3,\C^C) \rtimes \Z_2 \relmiddle{|} f_{4,\sC^C}(A,\varepsilon)=1 \right\rbrace.
\end{align*}
In the former case, let $ (A,1) \in \Ker\,f_{4,\sC^C} $. Then we will find $ A \in SU(3,\C^C) $ satisfying $ A^*XA=X $ for any $ X \in \mathfrak{J}(3,\C^C) $, so let $ E_1, E_2, E_3, F_1(1), F_3(1) $ as $ X $. Then we see that $ A $ is of the form $ \diag(a,a,a) $ with $ a^3=1 $. Hence we have
\begin{align*}
\left\lbrace (A,1) \in SU(3,\C^C) \rtimes \Z_2 \relmiddle{|} f_{4,\sC^C}(A,1)=1 \right\rbrace  \subset \left\lbrace (E,1), ({\boldsymbol\omega} E,1), ({\boldsymbol\omega}^2E,1)\right\rbrace
\end{align*}
and vice versa, where $ {\boldsymbol\omega}^3=1,{\boldsymbol\omega} \in \C,{\boldsymbol\omega} \not=1 $. In the latter case, from $ f_{4,\sC^C}(A,\varepsilon)=1 $ we have $ f_{4,\sC^C}(A,1)f_{4,\sC^C}(E,\varepsilon)=1 $, that is, $ f_{4,\sC^C}(A,1)=f_{4,\sC^C}(E,\varepsilon) $. Hence there exists no $ (A,\varepsilon) \in \Ker\,f_{4,\sC^C} $ such that $ f_{4,\sC^C}(A,\varepsilon)=1 $, so that
\begin{align*}
\left\lbrace (A,\varepsilon) \in SU(3,\C^C) \rtimes \Z_2 \relmiddle{|} f_{4,\sC^C}(A,\varepsilon)=1 \right\rbrace=\emptyset.
\end{align*}
Thus we obtain $ \Ker\,f_{4,\sC^C}=\left\lbrace (E,1), ({\boldsymbol\omega} E,1), ({\boldsymbol\omega}^2E,1)\right\rbrace=(\Z_3,1) $.

Therefore, from $ SU(3,\C^C) \cong SL(3,C) $, we have the required isomorphism
\begin{align*}
(F_{4,\sC})^C \cong SL(3,C)/\Z_3 \rtimes \Z_2,
\end{align*}
where $ \Z_3=\{(E,1),(\omega E,1),(\omega^2E,1)\}, \omega \in C,\omega^3=1, \omega \not=1 $.
\end{proof}

\vspace{3mm}

\uwave{\hspace{130mm}}
\vspace{3mm}

$ \bullet $ Memoir \,-\, An existence of $ r \in \C^C $ satisfying the equation $ r^3=\ov{p}q \in \C^C $.
\vspace{2mm}

First, from $ |p|=|q|=1 $, we see $ |\ov{p}q|=1 $, so that since $ r $ satisfies $ r^3=\ov{p}q $, we have $ |r|^3=1 $, that is, $ |r|=1, |r|=\omega $ or $ |r|=\omega^2 $.

Set $ r:=x+ye_1, x,y \in C $ and $ \ov{p}q:=a+be_1, a,b \in C $. Then it follows from
\begin{align*}
r^3=(x+ye_1)^3=(x^3-3xy^2)+(3x^2y-y^3)e_1,\;\; x^3-3xy^2, 3x^2y-y^3 \in C
\end{align*}
that
\begin{align*}
\left\lbrace \begin{array}{l}
x^3-3xy^2=a \,\cdots \, (1)
\vspace{1mm}\\
3x^2y-y^3=b \,\cdots \, (2).
\end{array}\right.
\end{align*}

In the case where $ |r|=1 $, that is, $ x^2+y^2=1 $, (1), (2) can be deformed as follows:
\begin{align*}
\left\lbrace \begin{array}{l}
4x^3-3x=a \,\cdots \, (1)'
\vspace{1mm}\\
-4y^3+3y=b \,\cdots \, (2)'.
\end{array}\right.
\end{align*}
Hence both of (1)' and (2)' have roots over $ C $, so that there exists $ r \in \C^C $ satisfying the equation $ r^3=\ov{p}q \in \C^C $.

In the case where $ |r|=\omega $, that is, $ x^2+y^2=\omega^2 $, (1), (2) can be deformed as follows:
\begin{align*}
\left\lbrace \begin{array}{l}
4x^3-3\omega^2 x=a \,\cdots \, (1)'
\vspace{1mm}\\
-4y^3+3\omega^2 y=b \,\cdots \, (2)'.
\end{array}\right.
\end{align*}
Hence, as in the case above, there exists $ r \in \C^C $ satisfying the equation $ r^3=\ov{p}q \in \C^C $.

In the case where $ |r|=\omega^2 $, that is, $ x^2+y^2=\omega $, (1), (2) can be deformed as follows:
\begin{align*}
\left\lbrace \begin{array}{l}
4x^3-3\omega x=a \,\cdots \, (1)'
\vspace{1mm}\\
-4y^3+3\omega y=b \,\cdots \, (2)'.
\end{array}\right.
\end{align*}
Hence, as in the case above, there exists $ r \in \C^C $ satisfying the equation $ r^3=\ov{p}q \in \C^C $.

With above, we can prove that an existence of $ r \in \C^C $ satisfying the equation $ r^3=\ov{p}q \in \C^C $.
\vspace{2mm}

$ \bullet $ The isomorphism $ SU(3,\C^C) \cong SL(3,C) $ in relation to $ \Ker\,f_{4,\sC^C} $.

We define a mapping $ f:SL(3,C) \to SU(3,\C^C) $ by
\begin{align*}
f(A)=\iota A+\ov{\iota}\,{}^t\!A^{-1}.
\end{align*}
Then $ f $ induces the isomorphism $ SL(3,C) \cong SU(3,\C^C) $. In addition, the inverse mapping $ g:SU(3,\C^C) \to SL(3,C) $ of $ f $ is given by
\begin{align*}
g(P+i\i Q)=P+Q.
\end{align*}
Its details are omitted. Now, in relation $ \Ker\,f_{\sC^C} $, under the mapping $ g $, we can confirm the following correspondence:
\begin{align*}
1 \mapsto 1,\quad \boldsymbol{\omega} \mapsto \omega^2, \quad \boldsymbol{\omega}^2 \mapsto \omega.
\end{align*}

\vspace{3mm}

\uwave{\hspace{130mm}}
\vspace{3mm}

Let the group $ E_{6,\sC} $ be the compact group as the closed subgroup of the unitary group $ U(9,C)=U(\mathfrak{J}(3,\C^C)):=\{\alpha \in \Iso_C(\mathfrak{J}(3,\C^C))\,|\, \langle \alpha X, \alpha Y\rangle=\langle X, Y \rangle\} $.

Then we have the following proposition.

\begin{proposition}\label{proposition 4.1.4}
Any element $ X \in \mathfrak{J}(3,\C^C) $ can be transformed to a diagonal form by some element $ \alpha \in (E_{6,\sC})_0 ${\rm :}
\begin{align*}
\alpha X=\begin{pmatrix}
\xi_1 & 0 & 0 \\
0 & \xi_2 & 0 \\
0 & 0 & \xi_3
\end{pmatrix}, \xi_i \in C
\end{align*}

Moreover, we can choose $ \alpha \in (E_{6,\sC})_0 $ so that two of $ \xi_1, \xi_2, \xi_3 $ are non-negative real numbers.
\end{proposition}
\begin{proof}
As in \cite[Proposition 3.8.2]{iy0}, we can also prove this proposition.
\end{proof}

Now, we prove the lemma used in the proof of theorem below.

\begin{lemma}\label{lemma 4.1.5}
The Lie algebra $ (\mathfrak{e}_{6,\sC})^C $ of the group $ (E_{6,\sC})^C $ is given by
\begin{align*}
(\mathfrak{e}_{6,\sC})^C:&=\left\lbrace \phi \in \Hom_C(\mathfrak{J}(3,\C^C)) \relmiddle{|} (\phi X, X, X)=0 \right\rbrace
\\
&=\left\lbrace \phi=\delta+\tilde{T} \relmiddle{|} \delta \in (\mathfrak{f}_{4,\sC})^C, T \in \mathfrak{J}(3,\C^C), \tr(T)=0 \right\rbrace.
\end{align*}

In particular, we have $ \dim_C((\mathfrak{e}_{6,\sC})^C)=8+(2+3\times 2)=16 $.
\end{lemma}
\begin{proof}
As in \cite[Theorem 3.2.1]{iy0}, we can also prove this lemma.
\end{proof}

We defined a space $ (EIV_{\sC})^C  $ by
\begin{align*}
 (EIV_{\sC})^C:=\left\lbrace X \in \mathfrak{J}(3,\C^C) \relmiddle{|}\det\,X=1 \right\rbrace.
\end{align*}

Then we have the following theorem.

\begin{theorem}\label{theorem 4.1.6}
The homogeneous space $ (E_{6,\sC})^C/(F_{4,\sC})^C $ is homeomorphic to the space $(EIV_{\sC})^C ${\rm :} $ (E_{6,\sC})^C/(F_{4,\sC})^C \simeq (EIV_{\sC})^C $.

Moreover, the group $ (E_{6,\sC})^C $ has at most two connected components.
\end{theorem}
\begin{proof}
First, the group $ (E_{6,\sC})^C $  acts on  $ (EIV_{\sC})^C $, obviously. We will prove that the action of $ (E_{6,\sC})^C $ on $ (EIV_{\sC})^C $ is transitive.

For a given $ X \in (EIV_{\sC})^C $, $ X $ can be transformed to a diagonal form by some $ \alpha \in (E_{6,\sC})_0 \subset E_{6,\sC} \subset (E_{6,\sC})^C $ (Proposition \ref{proposition 4.1.4}):
  \begin{align*}
  \alpha X=\begin{pmatrix}
  \xi_1 & 0 & 0 \\
  0 & \xi_2 & 0 \\
  0 & 0 & \xi_3
  \end{pmatrix}=:X',\,\, \xi_i \in C.
  \end{align*}
Moreover, we can choose $ \alpha \in ((E_{6,\sR})^C)_0 $ so that two of $ \xi_1,\xi_2,\xi_3 $ are non-negative real numbers (Proposition \ref{proposition 4.1.4}).
Hence, from $ \det\,(\alpha X)=\det \,X=1 $, that is, $ \xi_1\xi_2\xi_3=1 $, we have $ \xi_i >0, i=1,2,3 $.

Let the elements $ (s/2)(E_1-E_2)^\sim, (t/2)(E_2-E_3)^\sim \in (\mathfrak{e}_{6,\sC})^C, s,t \in \R $ (Theorem \ref{lemma 4.1.5}). We denote $ \exp((s/2)(E_1-E_2)^\sim), \exp((t/2)(E_2-E_3)^\sim) \in ((E_{6,\sC})^C)_0$ by $ \alpha_{12}(s), \alpha_{23}(t) $, respectively. Moreover, the explicit form of the actions of $ \alpha_{12}(s), \alpha_{23}(t) $ to $ \mathfrak{J}(3,\C^C) $ are respectively given as follows:
  \begin{align*}
  \alpha_{12}(s)X=
  \begin{pmatrix}
  e^s\xi_1 & x_3 & e^{s/2}x_2 \\
  x_3 & e^{-s}\xi_2 & e^{-s/2}x_1 \\
  e^{s/2}x_2 & e^{-s/2}x_1 & \xi_3
  \end{pmatrix},\;
  \alpha_{23}(t)X=
  \begin{pmatrix}
  \xi_1 & e^{t/2}x_3 & e^{-t/2}x_2 \\
  e^{t/2}x_3 & e^t\xi_2 & x_1 \\
  e^{-t/2}x_2 & x_1 & e^{-t}\xi_3
  \end{pmatrix},\;\; X \in \mathfrak{J}(3,\C^C).
  \end{align*}
Then, apply $ \alpha_{12}(s) $ on $ X' $, then $ \alpha_{12}(s)(X') $ is of the form $ \diag(e^s\xi_1, e^{-s}\xi_2, \xi_3) $:
$ \alpha_{12}(s)(X')=\diag(e^s\xi_1,\allowbreak  e^{-s}\xi_2, \xi_3) $. Since we can choose $ s_0 \in \R $ such that $ e^{s_0}\xi_1=1 $, together with $ \det(\alpha_{12}(s)(X'))=1 $, we have $ \alpha_{12}(s_0)(X')=\diag(1, \xi, 1/\xi)=:X'' $.
In addition,  apply $ \alpha_{23}(t) $ on $ X'' $, then $  \alpha_{23}(t)(X'') $ is of the form $ \diag(1, e^t\xi, 1/(e^t\xi)) $: $ \alpha_{23}(t)(X'')=\diag(1, e^t\xi, 1/(e^t\xi)) $.
As in the case above, since we can choose $ t_0 \in \R $ such that $ e^{t_0}\xi=1 $, we have $ \alpha_{23}(t_0)(X'')=\diag(1, 1, 1)=E $. This shows the transitivity of action to $ (EIV_{\sC})^C $ by the group $ (E_{6,\sC})^C $. The isotropy subgroup of $ (E_{6,\sC})^C $ at $ E $ is the group $ (F_{4,\sC})^C $.

Thus we have the required homeomorphim
  \begin{align*}
  (E_{6,\sC})^C/(F_{4,\sC})^C \simeq (EIV_{\sC})^C.
  \end{align*}
Moreover, since $ (F_{4,\sC})^C $ has two connected components (Theorem \ref{theorem 4.1.3}) and $ (EIV_{\sC})^C=((E_{6,\sC})^C)_0E $ is connected,
the group $ (E_{6,\sC})^C $ has at most two connected components.
\end{proof}

Let the discrete group $ \Z_2=\{1,\varepsilon\} $. This group acts on the group $ SU(3,\C^C) \times SU(3,\C^C) $ by
\begin{align*}
1(A,B)=(A,B), \quad \varepsilon(A,B)=(\ov{B}, \ov{A}),
\end{align*}
and then let $ (SU(3,\C^C) \times SU(3,\C^C)) \rtimes \Z_2 $ be the semi-direct product of $ SU(3,\C^C) \times SU(3,\C^C) $ and $ \Z_2 $ with the multiplication
\begin{align*}
((A_1,B_1),1)((A_2,B_2),1)&=((A_1A_2,B_1B_2),1), \quad
((A_1,B_1),1)((A_2,B_2),\varepsilon)=((A_1A_2,B_1B_2),\varepsilon),
\\
((A_1,B_1),\varepsilon)((A_2,B_2),1)&=((A_1\ov{B}_2,B_1\ov{A}_2),\varepsilon), \quad
((A_1,B_1),\varepsilon)((A_2,B_2),\varepsilon)=((A_1\ov{B}_2,B_1\ov{A}_2),1).
\end{align*}

Now, we prove the main theorem below.

\begin{theorem}\label{theorem 4.1.7}
The group $ (E_{6,\sC})^C $ is isomorphic to the semi-direct product of the groups $ (SL(3,C) \times SL(3,C))/\Z_3 $ and $ \Z_2, \Z_3=\{(E,E),(\omega E, \omega E), (\omega^2E,\omega^2E)\}, \Z_2=\{1,\varepsilon\} ${\rm :} $ (E_{6,\sC})^C \cong (SL(3,C) \times SL(3,C))/\Z_3 \rtimes \Z_2 $.
\end{theorem}
\begin{proof}
Let the group $ SL(3,C) $ as the group $ SU(3,\C^C) $. We define a mapping $ f_{6,\sC^C}: (SU(3,\C^C) \times SU(3,\C^C)) \rtimes \{1,\varepsilon\} \to (E_{6,\sC})^C $ by
\begin{align*}
f_{6,\sC^C}((A,B),1)X&=h(A,B)Xh(A,B)^*,
\\
f_{6,\sC^C}((A,B),\varepsilon)X&=h(A,B)\ov{X}h(A,B)^*, \;\; X \in \mathfrak{J}(3,\C^C),
\end{align*}
where the mapping $ h: M(3,\C^C) \times M(3,\C^C) \to M(3,\C^C) $ is defined by $ h(A,B)=\ov{\iota}A+\iota B, \iota:=(1/2)(1+i\i) $ with the properties $ \tau h(A,B)=h(\tau B, \tau A), \ov{h(A,B)}=h(\ov{B}, \ov{A}), h(A,B)^*=h(B^*,A^*) $.

First, we will prove that $ f_{6,\sC^C} $ is well-defined. It follows from
\begin{align*}
\det\,h(A,B)=\det(\ov{\iota}A+\iota B)=\ov{\iota}\det\,A+\iota\det\,B=\ov{\iota}+\iota=1
\end{align*}
that
\begin{align*}
\det(f_{6,\sC^C}((A,B),1)X)&=\det(h(A,B)Xh(A,B)^*)=(\det\,h(A,B))(\det\,X)(\det\,h(A,B)^*)
\\
&=\det\,X.
\end{align*}
Hence we have $ f_{6,\sC^C}((A,B),1) \in (E_{6,\sC})^C $. Subsequently, note that it is shown that $ f_{6,\sC^C} $ is a homomorphism below,
it follows from $ f_{6,\sC^C}((E,E),\varepsilon) \in (F_{4,\sC})^C \subset (E_{6,\sC})^C $(Theorem \ref{theorem 4.1.3})
that
\begin{align*}
f_{6,\sC^C}((A,B),\varepsilon)=f_{6,\sC^C}((A,B),1)((E,E),\varepsilon))
=f_{6,\sC^C}((A,B),1)f_{6,\sC^C}((E,E),\varepsilon) \in (E_{6,\sC})^C.
\end{align*}
With above, $ f_{6,\sC^C} $ is well-defined. It is easy to verify that $  f_{6,\sC^C} $ is a homomorphism. Indeed, note that the mapping $ h $ is a homomorphism, we have the following
\begin{align*}
f_{6,\sC^C}((A_1,B_1),1)f_{6,\sC^C}((A_2,B_2),1)X&=f_{6,\sC^C}((A_1,B_1),1)(h(A_2,B_2)Xh(A_2,B_2)^*)
\\
&=h(A_1,B_1)(h(A_2,B_2)Xh(A_2,B_2)^*)h(A_1,B_1)^*
\\
&=h(A_1A_2,B_1B_2)Xh((B_1B_2)^*,(A_1A_2)^*)
\\
&=h(A_1A_2,B_1B_2)Xh(A_1A_2,B_1B_2)^*
\\
&=f_{6,\sC^C}((A_1A_2,B_1B_2),1)X
\\
&=f_{6,\sC^C}((A_1,B_1),1)((A_2,B_2),1))X,
\\[1mm]
f_{6,\sC^C}((A_1,B_1),1)f_{6,\sC^C}((A_2,B_2),\varepsilon)X&=f_{6,\sC^C}((A_1,B_1),1)(h(A_2,B_2)\ov{X}h(A_2,B_2)^*)
\\
&=h(A_1,B_1)(h(A_2,B_2)\ov{X}h(A_2,B_2)^*)h(A_1,B_1)^*
\\
&=f_{6,\sC^C}((A_1A_2,B_1B_2),\varepsilon)X
\\
&=f_{6,\sC^C}((A_1,B_1),1)((A_2,B_2),\varepsilon))X,
\\[1mm]
f_{6,\sC^C}((A_1,B_1),\varepsilon)f_{6,\sC^C}((A_2,B_2),1)X&=f_{6,\sC^C}((A_1,B_1),\varepsilon)(h(A_2,B_2)Xh(A_2,B_2)^*)
\\
&=h(A_1,B_1)(\ov{h(A_2,B_2)Xh(A_2,B_2)^*})h(A_1,B_1)^*
\\
&=h(A_1,B_1)h(\ov{B}_2,\ov{A}_2)\ov{X}h({}^t\!{A_2},{}^t\!{B_2})h({B_1}^*,{A_1}^*)
\\
&=h(A_1\ov{B}_2,B_1\ov{A}_2)\ov{X}h({}^t\!{A_2}{B_1}^*,{}^t\!{B_2}{A_1}^*)
\\
&=h(A_1\ov{B}_2,B_1\ov{A}_2)\ov{X}h(A_1\ov{B}_2,B_1\ov{A}_2)^*
\\
&=f_{6,\sC^C}((A_1\ov{B}_2,B_1\ov{A}_2),\varepsilon)X
\\
&=f_{6,\sC^C}((A_1,B_1),\varepsilon)((A_2,B_2),\varepsilon))X,
\\[1mm]
f_{6,\sC^C}((A_1,B_1),\varepsilon)f_{6,\sC^C}((A_2,B_2),\varepsilon)X&=f_{6,\sC^C}((A_1,B_1),\varepsilon)(h(A_2,B_2)\ov{X}h(A_2,B_2)^*)
\\
&=h(A_1,B_1)(\ov{h(A_2,B_2)\ov{X}h(A_2,B_2)^*})h(A_1,B_1)^*
\\
&=h(A_1,B_1)h(\ov{B}_2,\ov{A}_2)Xh({}^t\!{A_2},{}^t\!{B_2})h({B_1}^*,{A_1}^*)
\\
&=h(A_1\ov{B}_2,B_1\ov{A}_2)Xh({}^t\!{A_2}{B_1}^*,{}^t\!{B_2}{A_1}^*)
\\
&=h(A_1\ov{B}_2,B_1\ov{A}_2)Xh(A_1\ov{B}_2,B_1\ov{A}_2)^*
\\
&=f_{6,\sC^C}((A_1\ov{B}_2,B_1\ov{A}_2),1)X
\\
&=f_{6,\sC^C}((A_1,B_1),\varepsilon)((A_2,B_2),\varepsilon))X.
\end{align*}

Next, we will determine $ \Ker\,f_{6,\sC^C} $. It follows from the definition of kernel that
\begin{align*}
\Ker\,f_{6,\sC^C}&=\left\lbrace ((A,B),1) \in (SU(3,\C^C) \times SU(3,\C^C )) \rtimes \Z_2 \relmiddle{|} f_{6,\sC^C}((A,B),1)=1 \right\rbrace
\\
& \cup \left\lbrace ((A,B),\varepsilon) \in(SU(3,\C^C) \times SU(3,\C^C )) \rtimes \Z_2 \relmiddle{|} f_{6,\sC^C}((A,B),\varepsilon)=1 \right\rbrace.
\end{align*}
In the former case, let $ ((A,B),1) \in \Ker\,f_{6,\sC^C} $. Then we will find $ (A,B) \in SU(3,\C^C) \times SU(3,\C^C ) $ satisfying $ h(A,B)Xh(A,B)^*=X $ for any $ X \in \mathfrak{J}(3,\C^C) $, that is, $ h(A,B)X=Xh(B,A) $, so let $ E_1,E_2, E_3, F_1(1), F_3(1) $ as $ X $. Then we see that $ A=B $ and $ A, B $ is of the form $ \diag(a,a,a) $ with $ a^3=1 $. Hence we have
\begin{align*}
&\left\lbrace ((A,B),1) \in (SU(3,\C^C) \times SU(3,\C^C)) \rtimes \Z_2 \relmiddle{|} f_{6,\sC^C}((A,B),1)=1 \right\rbrace
\\
&\hspace{50mm}\cap
\\
&\hspace{20mm}\left\lbrace ((E,E),1), (({\boldsymbol\omega} E,{\boldsymbol\omega} E),1), (({\boldsymbol\omega}^2 E,{\boldsymbol\omega}^2 E),1)\right\rbrace
\end{align*}
and vice versa, where $ {\boldsymbol\omega}^3=1,{\boldsymbol\omega} \in \C, {\boldsymbol\omega} \not=1 $. In the latter case, from $ f_{6,\sC^C}((A,B),\varepsilon)=1 $ we have $ f_{6,\sC^C}((A,B),1)f_{6,\sC^C}(E,\varepsilon)=1 $, that is, $ f_{6,\sC^C}((A,B),1)=f_{6,\sC^C}((E,E),\varepsilon) $. Hence there exists no $ ((A,B),\varepsilon) \in \Ker\,f_{6,\sC^C} $ such that $ f_{6,\sC^C}((A,B),\varepsilon)=1 $, so that
\begin{align*}
\left\lbrace ((A,B),\varepsilon) \in (SU(3,\C^C) \times SU(3,\C^C)) \rtimes \Z_2 \relmiddle{|} f_{6,\sC^C}((A,B),\varepsilon)=1 \right\rbrace=\emptyset.
\end{align*}
Thus we obtain $ \Ker\,f_{6,\sC^C}=\left\lbrace ((E,E),1), ({\boldsymbol\omega} E,{\boldsymbol\omega} E),1), (({\boldsymbol\omega}^2 E,{\boldsymbol\omega}^2 E),1)\right\rbrace=(\Z_3,1) $.

Finally, we will prove that $ f_{6,\sC^C} $ is surjective. Since $ \Ker\,f_{6,\sC^C} $ is discrete, its connected component containing unit element is $ \{E\} $ only. Hence the Lie algebra of $ \Ker\,f_{6,\sC^C} $ is $ \{0\} $. Thus the differential mapping $ {f_{6,\sC^C}}_* : \mathfrak{su}(3,\C^C) \oplus \mathfrak{su}(3,\C^C) \to (\mathfrak{e}_{6,\sC})^C $ of the mapping $ f_{6,\sC^C} $ induces an injective homomorphism, and together with $ \dim_C(\mathfrak{su}(3,\C^C) \oplus \mathfrak{su}(3,\C^C))=8+8=16=\dim_C((\mathfrak{e}_{6,\sC})^C) $ (Lemma \ref{lemma 4.1.5}), $ {f_{6,\sC^C}}_* $ is surjective. Therefore, since the group $ ((E_{6,\sC^C})^C)_0 $ is connected, the mapping $ f_{6,\sC^C}: SU(3,\C^C) \times SU(3,\C^C) \to ((E_{6,\sC^C})^C)_0 $ induces a surjective homomorphism. However, $ \varepsilon=f_{6,\sC^C}((E,E),\varepsilon) $ does not be contained in $ ((E_{6,\sC^C})^C)_0 $. Indeed, if $ \varepsilon \in ((E_{6,\sC^C})^C)_0 $, since $ f_{6,\sC^C}: SU(3,\C^C) \times SU(3,\C^C) \to ((E_{6,\sC^C})^C)_0 $ is surjective, there exists $ (A,B) \in SU(3,\C^C) \times SU(3,\C^C) $ such that $ \varepsilon=f_{6,\sC^C}(A,B) $, that is, $ \ov{X}=h(A,B)Xh(A,B)^* $ for all $ X \in \mathfrak{J}(3,\C^C) $.
Let $ E_1,E_2,E_3,F_1(1), F_3(1) $ as $ X $, then as computed in $ \Ker\,f_{6,\sC^C} $, we have $ h(A,B)=E, h(A,B)=\omega E$ or $ h(A,B)=\omega^2 E $. Here, let $ F_1(\i) \in \mathfrak{J}(3,\C^C)$. Then, for those $ h(A,B) $, we have $ \ov{X}\not=h(A,B)Xh(A,B)^* $. Hence  there exists no $ (A,B) \in SU(3,\C^C) \times SU(3,\C^C) $ such that $ \varepsilon=f_{6,\sC^C}(A,B) $, that is, $ \varepsilon \not\in ((E_{6,\sC^C})^C)_0 $.
Thus $ (E_{6,\sC})^C $ has just two connected components (Theorem \ref{theorem 4.1.6}): $ (E_{6,\sC})^C=((E_{6,\sC})^C)_0 \cup ((E_{6,\sC})^C)_0 \cdot \varepsilon $. Let $ \beta=\alpha \varepsilon \in ((E_{6,\sC})^C)_0 \cdot \varepsilon, \alpha \in ((E_{6,\sC})^C)_0 $. Then there exists $ (A,B) \in SU(3,\C^C) \times SU(3,\C^C) $ such that $ \beta=\alpha\varepsilon=f_{6,\sC^C}((A,B),1)\varepsilon=f_{6,\sC^C}((A,B),\varepsilon) $. With above, the proof of surjective is completed.

Therefore we have the required isomorphism
\begin{align*}
(E_{6,\sC})^C \cong (SU(3,\C^C) \times SU(3,\C^C))/\Z_3 \rtimes \Z_2.
\end{align*}
\end{proof}

We prove the lemma used in the subsequent sections.

\begin{lemma}\label{lemma 4.1.8}
The mapping $ f_{6,\sC^C}: SU(3,\C^C) \times SU(3,\C^C)) \rtimes \{1,\varepsilon\} \to (E_{6,\sC})^C  $ of Theorem {\rm \ref{theorem 4.1.7}} satisfies the following
\begin{align*}
&(1)\,{}^t(f_{6,\sC^C}((A,B),1))^{-1}=f_{6,\sC^C}((B,A),1). \;\; (2)\,\tau f_{6,\sC^C}((A,B),1)\tau=f_{6,\sC^C}((\tau B,\tau A),1).
\\[1mm]
&(3)\gamma_{{}_{\scalebox{0.8}{\sC}}}f_{6,\sC^C}((A,B),1)\gamma_{{}_{\scalebox{0.8}{\sC}}}=f_{6,\sC^C}((\gamma_{{}_{\scalebox{0.8}{\sC}}} B,\gamma_{{}_{\scalebox{0.8}{\sC}}} A),1).\;\;
(4) \sigma f_{6,\sC^C}((A,B),1)\sigma=f_{6,\sC^C}((I_1AI_1,I_1BI_1),1).
\\[1mm]
&(5)\,{}^t(f_{6,\sC^C}((E,E),\varepsilon))^{-1}=f_{6,\sC^C}((E,E),\varepsilon). \;\;
(6)\,\tau f_{6,\sC^C}((E,E),\varepsilon)\tau=f_{6,\sC^C}((E,E),\varepsilon).
\\[1mm]
&(7)\gamma_{{}_{\scalebox{0.8}{\sC}}}f_{6,\sC^C}((E,E),\varepsilon)\gamma_{{}_{\scalebox{0.8}{\sC}}}=f_{6,\sC^C}((E,E),\varepsilon).\;\;
(8) \sigma f_{6,\sC^C}((E,E),\varepsilon)\sigma=f_{6,\sC^C}((E,E),\varepsilon),
\end{align*}
where $ I_1:=\diag(-1,1,1) $.
\end{lemma}
\begin{proof}
(1) It follows from
\begin{align*}
({}^t(f_{6,\sC^C}((A,B),1))X,Y)&=(X, f_{6,\sC^C}((A,B),1)Y)=(X,h(A,B)Yh(A,B)^*)
\\
&=(h(A,B)^*Xh(A,B),Y)=(h(B^*,A^*)Xh(B^*,A^*)^*,Y)
\\
&=(f_{6,\sC^C}((B^*,A^*),1)X,Y)
\\
&=(f_{6,\sC^C}((B^{-1},A^{-1}),1)X,Y)
\end{align*}
that $ {}^t(f_{6,\sC^C}((A,B),1))=f_{6,\sC^C}((B^{-1},A^{-1}),1) $. Hence we have
\begin{align*}
{}^t(f_{6,\sC^C}((A,B),1))^{-1}=f_{6,\sC^C}((B,A),1).
\end{align*}
\vspace{1mm}

(2) It follows from
\begin{align*}
\tau f_{6,\sC^C}((A,B),1)\tau X&=\tau(h(A,B)(\tau X)h(A,B)^*)=\tau h(A,B)X\tau h(B^*,A^*)
\\
&=h(\tau B,\tau A)Xh(\tau A^*,\tau B^*)=h(\tau B,\tau A)Xh(\tau B,\tau A)^*
\\
&=f_{6,\sC^C}((\tau B,\tau A),1)X
\end{align*}
that $ \tau f_{6,\sC^C}((A,B),1)\tau=f_{6,\sC^C}((B,A),1) $.
\vspace{2mm}

(3) It follows from
\begin{align*}
\gamma_{{}_{\scalebox{0.8}{\sC}}}f_{6,\sC^C}((A,B),1)\gamma_{{}_{\scalebox{0.8}{\sC}}}X&=\ov{(h(A,B)\ov{X}h(A,B)^*)}
\\
&=\ov{h(A,B)}\, X \,\ov{h(B^*,A^*)}
\\
&=h(\ov{B},\ov{A}) X h(\ov{B},\ov{A})^*
\\
&=f_{6,\sC^C}((\ov{B},\ov{A}),1)X
\\
&=f_{6,\sC^C}((\gamma_{{}_{\scalebox{0.8}{\sC}}} B,\gamma_{{}_{\scalebox{0.8}{\sC}}} A),1)X
\end{align*}
that $ \gamma_{{}_{\scalebox{0.8}{\sC}}}f_{6,\sC^C}((A,B),1)\gamma_{{}_{\scalebox{0.8}{\sC}}}=f_{6,\sC^C}((\gamma_{{}_{\scalebox{0.8}{\sC}}} B,\gamma_{{}_{\scalebox{0.8}{\sC}}} A),1) $.
\vspace{2mm}

(4)  It follows from
\begin{align*}
\sigma f_{6,\sC^C}((A,B),1)\sigma X&=\sigma(h(A,B)(\sigma X)h(A,B)^*)
\\
&=I_1(h(A,B)(I_1XI_1)h(A,B)^*)I_1
\\
&=I_1((\ov{\iota}A+\iota B)(I_1XI_1)(\ov{\iota}B^*+\iota A^*))I_1
\\
&=(\ov{\iota}(I_1AI_1)+\iota (I_1BI_1))X(\ov{\iota}(I_1B^*I_1)+\iota (I_1A^*I_1))
\\
&=(\ov{\iota}(I_1AI_1)+\iota (I_1BI_1))X(\ov{\iota}(I_1BI_1)^*+\iota (I_1AI_1)^*)
\\
&=h(I_1AI_1, I_1BI_1)Xhh(I_1AI_1, I_1BI_1)^*
\\
&=f_{6,\sC^C}((I_1AI_1,I_1BI_1),1)X
\end{align*}
that $ \sigma f_{6,\sC^C}((A,B),1)\sigma=f_{6,\sC^C}((I_1AI_1,I_1BI_1),1) $.
\vspace{2mm}

(5) Since $ f_{6,\sC^C}((E,E),\varepsilon)=\varepsilon \in (F_{4,\sC})^C$ (Theorem \ref{theorem 4.1.3}), note that the transpose $ {}^t\!\varepsilon $ of $ \varepsilon $ is defined by $ ({}^t\!\varepsilon X,Y)=(X,\varepsilon Y), X,Y \in \mathfrak{J}(3,\C^C) $, it follows from
\begin{align*}
({}^t\!\varepsilon X,Y)=(X,\varepsilon Y)=(X,\ov{Y})=(\ov{X},Y)=(\varepsilon X,Y)
\end{align*}
that $ {}^t\!\varepsilon=\varepsilon $. Moreover, from $ \varepsilon^2=1 $, we have $ {}^t\!\varepsilon^{-1}=\varepsilon $. Hence this implies the required result.

The formulas (6), (7), (8) are easy verified, so the proofs are omitted.
\end{proof}

\subsection{The group $ E_{6,\sC} $}

The structure of the group $ E_{6,\sC} $ has been determined by Ichiro Yokota and his school's members (\cite[Theorem 12]{iy9}). In the present article, although it is essentially the same as the proof of \cite[Theorem 12]{iy9}, we will determine the structure of the group $ E_{6,\sC} $ using the result of $ (E_{6,\sC})^C $.

As mentioned in the surface of this section, since $ (E_{6,\sC})^C $ has an involutive automorphism $ \tilde{\tau\lambda} $, we can define a subgroup $ ((E_{6,\sC})^C)^{\tau\lambda} $ of $ (E_{6,\sC})^C $:
\begin{align*}
((E_{6,\sC})^C)^{\tau\lambda}:=\left\lbrace \alpha \in (E_{6,\sC})^C \relmiddle{|} \tilde{\tau\lambda}(\alpha)=\alpha \right\rbrace.
\end{align*}

First, we prove the following theorem.

\begin{theorem}\label{theorem 4.2.1}
The group $ ((E_{6,\sC})^C)^{\tau\lambda} $ coincides with the group $ E_{6,\sC} ${\rm :} $ ((E_{6,\sC})^C)^{\tau\lambda}=E_{6,\sC} $.
\end{theorem}
\begin{proof}
Let $ \alpha \in ((E_{6,\sC})^C)^{\tau\lambda} $. Then it follows from $ \tau{}^t\!\alpha^{-1}\tau=\alpha $ that
\begin{align*}
\langle \alpha X, \alpha Y \rangle=(\tau \alpha X, \alpha Y)=({}^t\!\alpha^{-1}\tau X,\alpha Y)=(\tau X,\alpha^{-1}\alpha Y)=(\tau X,Y)=\langle X , Y \rangle.
\end{align*}
Hence we see $ \alpha \in E_{6,\sC} $. Conversely, let $ \beta \in E_{6,\sC} $. Then it follows from $ \langle \beta X, \beta Y \rangle=\langle X , Y \rangle $ that
\begin{align*}
(\tau X,Y)=\langle X , Y \rangle=\langle \beta X, \beta Y \rangle=(\tau\beta X, \beta Y)=({}^t\!\beta\beta\tau X, Y),
\end{align*}
that is, $ \tau={}^t\!\beta\beta\tau $. Hence we have $ \tau{}^t\!\beta^{-1}\tau=\beta $, so that $ \beta \in ((E_{6,\sC})^C)^{\tau\lambda} $.

With above, the proof of this theorem is completed.
\end{proof}

Let the discrete group $ \Z_2=\{1,\varepsilon\} $. This group acts on $ SU(3) \times SU(3) $ by
\begin{align*}
1(A,B)=(A,B),\quad \varepsilon(A,B)=(\ov{B},\ov{A}),
\end{align*}
and then $ (SU(3) \times SU(3)) \rtimes \Z_2 $ be the semi-direct product of $ SU(3) \times SU(3) $ and $ \Z_2 $ with the same multiplication of $ (SU(3,\C^C) \times SU(3,\C^C)) \rtimes \Z_2 $.

Now, we determine the structure of the group $ E_{6,\sC} $.

\begin{theorem}\label{theorem 4.2.2}
The group $ E_{6,\sC} $ is isomorphic to the semi-direct product of the groups $ (SU(3) \times SU(3))/\Z_3 $ and $ \Z_2, \Z_3=\{(E,E),(\boldsymbol{\omega} E, \boldsymbol{\omega} E), (\boldsymbol{\omega}^2E,\boldsymbol{\omega}^2E)\}, \Z_2=\{1,\varepsilon\} ${\rm :} $ E_{6,\sC} \cong (SU(3) \times SU(3))/\Z_3 \allowbreak \rtimes \Z_2 $.
\end{theorem}
\begin{proof}
Let $ E_{6,\sC} $ as the group $ ((E_{6,\sC})^C)^{\tau\lambda} $ (Theorem \ref{theorem 4.2.1}). We define a mapping $ f_{6,\tau\lambda}: (SU(3) \times SU(3)) \rtimes \{1,\varepsilon\} \to ((E_{6,\sC})^C)^{\tau\lambda} $ by
\begin{align*}
f_{6,\tau\lambda}((A,B),1)X&=h(A,B)Xh(A,B)^*,
\\
f_{6,\tau\lambda}((A,B),\varepsilon)X&=h(A,B)\ov{X}h(A,B)^*, \;\; X \in \mathfrak{J}(3,\C^C).
\end{align*}
Note that this mapping is the restriction of the mapping $ f_{6,\sC^C} $.
First, we will prove that $ f_{6,\tau\lambda} $ is well-defined and a homomorphism. Since $ f_{6,\tau\lambda} $ is the restriction of the mapping $ f_{6,\sC^C} $, it is clear that $ f_{6,\tau\lambda}((A,B),1),f_{6,\tau\lambda}((A,B),\varepsilon) \allowbreak \in (E_{6,\sC})^C $ and $ f_{6,\tau\lambda} $ is a homomorphism. Here, we show the following claim needed in order to prove $ f_{6,\tau\lambda}((A,B),1), f_{6,\tau\lambda}((A,B),\varepsilon) \in ((E_{6,\sC})^C)^{\tau\lambda} $.
By using Lemma \ref{lemma 4.1.8} (1), (2), we have $ \tau {}^t(f_{6,\tau\lambda}((A,B),1))^{-1}\tau=f_{6,\tau\lambda}((A,B),1) $, so that $ f_{6,\tau\lambda}((A, \allowbreak B),1) \in ((E_{6,\sC})^C)^{\tau\lambda} $. As in that above, by using  Lemma \ref{lemma 4.1.8} (5), (6), we have $ f_{6,\tau\lambda}((E,E),\varepsilon) \in ((E_{6,\sC})^C)^{\tau\lambda} $, so that since $ f_{6,\tau\lambda}((A,B),\varepsilon)=f_{6,\tau\lambda}((A,B),1)f_{6,\tau\lambda}((E,E),\varepsilon) $, we see $ f_{6,\tau\lambda}((A,B),\varepsilon) \in ((E_{6,\sC})^C)^{\tau\lambda} $. With above, the proof of well-defined is proved.

Next, we will prove that $ f_{6,\tau\lambda} $ is surjective. Let $ \alpha \in E_{6,\sC}=((E_{6,\sC})^C)^{\tau\lambda} \subset (E_{6,\sC})^C $. Then there exists $ ((P,Q),1) \in (SU(3,\C^C) \times SU(3,\C^C)) \rtimes \Z_2 $ such that $ \alpha=f_{6,\sC^C}((P,Q),1) $ or there exists $ ((P,Q),\varepsilon) \in (SU(3,\C^C) \times SU(3,\C^C)) \rtimes \Z_2 $ such that $ \alpha=f_{6,\sC^C}((P,Q),\varepsilon) $ (Theorem \ref{theorem 4.1.7}). Moreover $ \alpha $ satisfies the condition $ \tau {}^t\!\alpha^{-1} \tau=\alpha $, that is, $ \tau{}^t(f_{6,\sC^C}((P,Q),1))^{-1}\tau=f_{6,\sC^C}((P,Q),1) $ or $ \tau{}^t(f_{6,\sC^C}((P,Q),\varepsilon))^{-1}\tau=f_{6,\sC^C}((P,Q),\varepsilon) $. Since it follows from  Lemma \ref{lemma 4.1.8} (1), (2), (5), (6) that
\begin{align*}
\tau{}^t(f_{6,\sC^C}((P,Q),1))^{-1}\tau&=f_{6,\sC^C}((\tau P,\tau Q),1)
\\
&\text{or}
\\
\tau{}^t(f_{6,\sC^C}((P,Q),\varepsilon))^{-1}\tau&=f_{6,\sC^C}((\tau P,\tau Q),\varepsilon),
\end{align*}
we have the following
\begin{align*}
\left\lbrace
\begin{array}{l}
\tau P=P \\
\tau Q=Q,
\end{array}
\right. \quad
\left\lbrace
\begin{array}{l}
\tau P=\omega P \\
\tau Q=\omega Q
\end{array}
\right. \quad \text{or} \quad
\left\lbrace
\begin{array}{l}
\tau P=\omega ^2 P \\
\tau Q=\omega ^2 Q.
\end{array}
\right.
\end{align*}
In the first case, we have $ P,Q \in SU(3) $. In the others cases, since we have $ P=Q=0 $, those are impossible. The proof of surjective is completed.

Finally, we will determine $ \Ker\,f_{6,\tau\lambda} $. Since $ f_{6,\tau\lambda} $ is the restriction of the mapping $ f_{6,\sC^C} $, it is easy to obtain $ \Ker\,f_{6,\tau\lambda}=\Ker\,f_{6,\sC^C}=\left\lbrace ((E,E),1), (({\boldsymbol\omega} E,{\boldsymbol\omega} E),1), (({\boldsymbol\omega}^2 E,{\boldsymbol\omega}^2 E),1)\right\rbrace \cong (\Z_3,1) $.

Therefore we have the required isomorphism
\begin{align*}
E_{6,\sC} \cong (SU(3) \times SU(3))/\Z_3 \allowbreak \rtimes \Z_2.
\end{align*}
\end{proof}

\subsection{The group $ E_{6(6), \sC'} $}

We consider a subspace $ (\C^C)_{\tau\gamma_{{}_{\scalebox{0.8}{\sC}}}} $ of $ \C^C $:
\begin{align*}
(\C^C)_{\tau\gamma_{{}_{\scalebox{0.8}{\sC}}}}:&=\left\lbrace z \in \C^C \relmiddle{|} {\tau\gamma_{{}_{\scalebox{0.8}{\sC}}}} z=z \right\rbrace
\\
&=\left\lbrace z=x+ y i\i \relmiddle{|} x,y \in \R \right\rbrace.
\end{align*}

Let $ \C' $ be the algebra of split complex numbers: $ \C':=\R \oplus \R\i',{\i'}^2=1 $. Then the correspondence
\begin{align*}
(\C^C)_{\tau\gamma_{{}_{\scalebox{0.8}{\sC}}}} \ni x+y i\i
\underset{l}{\mapsto} x+y\i' \in \C'
\end{align*}
gives an isomorphism $ (\C^C)_{\tau\gamma_{{}_{\scalebox{0.8}{\sC}}}} \simeq \C' $ as algebras.

We define a subalgebra $ (\mathfrak{J}(3,\C^C))_{\tau\gamma_{{}_{\scalebox{0.8}{\sC}}}} $ of $ \mathfrak{J}(3,\C^C) $ by
\begin{align*}
(\mathfrak{J}(3,\C^C))_{\tau\gamma_{{}_{\scalebox{0.8}{\sC}}}}:&=\left\lbrace X \in \mathfrak{J}(3,\C^C) \relmiddle{|} {\tau\gamma_{{}_{\scalebox{0.8}{\sC}}}} X=X \right\rbrace
\\
&=\left\lbrace
X=\begin{pmatrix}
\xi & x_3 & \ov{x}_2 \\
\ov{x}_3 & \xi_2 & x_1 \\
x_2 & \ov{x}_1 & \xi_3
\end{pmatrix} \relmiddle{|} \xi_i \in \R, \x_i \in (\C^C)_{\tau\gamma_{{}_{\scalebox{0.8}{\sC}}}} \right\rbrace.
\end{align*}
Let $ \mathfrak{J}(3,\C') $ be the split Jordan algebra. Then the correspondence
\begin{align*}
(\mathfrak{J}(3,\C^C))_{\tau\gamma_{{}_{\scalebox{0.8}{\sC}}}}\ni
\begin{pmatrix}
\xi & x_3 & \ov{x}_2 \\
\ov{x}_3 & \xi_2 & x_1 \\
x_2 & \ov{x}_1 & \xi_3
\end{pmatrix} \underset{f}{\mapsto}
\begin{pmatrix}
\xi & l(x_3) & \ov{l(x_2)} \\
\ov{l(x_3)} & \xi_2 & l(x_1) \\
l(x_2) & \ov{l(x_1)} & \xi_3
\end{pmatrix} \in \mathfrak{J}(3,\C')
\end{align*}
gives an isomorphism $ (\mathfrak{J}(3,\C^C))_{\tau\gamma_{{}_{\scalebox{0.8}{\sC}}}} \simeq  \mathfrak{J}(3,\C') $ as algebras with the properties $ \det(fX)=\det\,X $.

We define a group $ E_{6(6),(\sC^C)_{\tau\gamma_{{}_{\scalebox{0.8}{\sC}}}}} $ by
\begin{align*}
E_{6(6),(\sC^C)_{\tau\gamma_{{}_{\scalebox{0.8}{\sC}}}}}:=\left\lbrace \alpha \in \Iso_{\sR}((\mathfrak{J}(3,\C^C))_{\tau\gamma_{{}_{\scalebox{0.8}{\sC}}}}) \relmiddle{|}\det(\alpha X)=\det\, X\right\rbrace.
\end{align*}

Then we have the following proposition.

\begin{proposition}\label{proposition 4.3.1}
The group $ E_{6(6),(\sC^C)_{\tau\gamma_{{}_{\scalebox{0.8}{\sC}}}} } $ is isomorphic to the group $ E_{6,\sC'} ${\rm :} $ E_{6(6),(\sC^C)_{\tau\gamma_{{}_{\scalebox{0.8}{\sC}}}} } \cong E_{6(6),\sC'} $.
\end{proposition}
\begin{proof}
We consider the following diagram:
    \[
    \begin{CD}
        (\mathfrak{J}(3,\C^C))_{\tau\gamma_{{}_{\scalebox{0.8}{\sC}}}} @>{\alpha}>>  (\mathfrak{J}(3,\C^C))_{\tau\gamma_{{}_{\scalebox{0.8}{\sC}}}} \\
        @V{f}VV    @V{f}VV   \\
        \mathfrak{J}(3,\C')   @>{\beta}>>   \mathfrak{J}(3,\C')
    \end{CD}
    \]

Then we define a mapping $ \varphi_{6(6)}: E_{6(6),(\sC^C)_{\tau\gamma_{{}_{\scalebox{0.8}{\sC}}}} } \to E_{6(6),\sC'} $ by
\begin{align*}
\varphi_{6(6)}(\alpha)=f\alpha f^{-1}.
\end{align*}

First, we will prove that $ \varphi_{6(6)} $ is well-defined and a homomorphism. It follows from $ \det(fX)=\det\,X, X \in (\mathfrak{J}(3,\C^C))_{\tau\gamma_{{}_{\scalebox{0.8}{\sC}}}} $ that
\begin{align*}
\det(\varphi_{6(6)}(\alpha)X')&=\det(f\alpha f^{-1}X')=\det(\alpha f^{-1}X')=\det(f^{-1}X')
\\
&=\det\,X',\;\; X' \in \mathfrak{J}(3,\C').
\end{align*}
Hence we have $ \varphi_{6(6)}(\alpha) \in E_{6(6),\sC'} $, so that $ \varphi_{6(6)} $ is well-defined. In addition, it is clear that $ \varphi_{6(6)} $ is a homomorphism.

Next, we will prove that $ \varphi_{6(6)} $ is surjective. Let $ \beta \in E_{6(6),\sC'} $. Then there exists $ \alpha \in E_{6(6),(\sC^C)_{\tau\gamma_{{}_{\scalebox{0.8}{\sC}}}} } $ such that $ \beta=f \alpha f^{-1} $. Indeed, as shown in the proof of well-defined above, it is easy to verify that $ \alpha=f^{-1}\beta f \in E_{6(6),(\sC^C)_{\tau\gamma_{{}_{\scalebox{0.8}{\sC}}}} } $.

Finally, we will prove that $ \varphi_{6(6)} $ is injective, however it is clear.

Therefore we have the required isomorphism
\begin{align*}
E_{6(6),(\sC^C)_{\tau\gamma_{{}_{\scalebox{0.8}{\sC}}}} } \cong E_{6(6),\sC'}.
\end{align*}
\end{proof}

Since the group $ (E_{6,\sC})^C $ has an involutive automorphism $ \tilde{\tau\gamma_{{}_{\scalebox{0.8}{\sC}}}} $, we can define a subgroup $ ((E_{6,\sC})^C)^{\tau\gamma_{{}_{\scalebox{0.8}{\sC}}}} $ of $ (E_{6,\sC})^C $ by
\begin{align*}
((E_{6,\sC})^C)^{\tau\gamma_{{}_{\scalebox{0.8}{\sC}}}}:=\left\lbrace\alpha \in (E_{6,\sC})^C \relmiddle{|}\tilde{\tau\gamma_{{}_{\scalebox{0.8}{\sC}}}}(\alpha)=\alpha \right\rbrace.
\end{align*}

Then we prove the following theorem.

\begin{theorem}\label{theorem 4.3.2}
The group $ ((E_{6,\sC})^C)^{\tau\gamma_{{}_{\scalebox{0.8}{\sC}}}} $ coincides with the group $ E_{6(6),(\sC^C)_{\tau\gamma_{\scalebox{0.8}{\sC}}} } ${\rm :} $ ((E_{6,\sC})^C)^{\tau\gamma_{{}_{\scalebox{0.8}{\sC}}}}=E_{6(6),(\sC^C)_{\tau\gamma_{\scalebox{0.8}{\sC}}} } $.

In particular, we have the isomorphism $ ((E_{6,\sC})^C)^{\tau\gamma_{{}_{\scalebox{0.8}{\sC}}}} \cong E_{6(6),\sC'} $.
\end{theorem}
\begin{proof}
Let $ \alpha \in ((E_{6,\sC})^C)^{\tau\gamma_{{}_{\scalebox{0.8}{\sC}}}} $. Since $ (\tau\gamma_{{}_{\scalebox{0.8}{\sC}}})\alpha=\alpha(\tau\gamma_{{}_{\scalebox{0.8}{\sC}}}) $, $ \alpha $ induces an $ \R $-linear isomorphism of $ (\mathfrak{J}(3,\C^C))_{\tau\gamma_{{}_{\scalebox{0.8}{\sC}}}} $. Moreover, since it is clear that $ \det(\alpha X)=\det\,X, X \in (\mathfrak{J}(3,\C^C))_{\tau\gamma_{{}_{\scalebox{0.8}{\sC}}}} $, we have $ \alpha \in E_{6(6),(\sC^C)_{\tau\gamma_{\scalebox{0.8}{\sC}}} } $. Conversely, let $ \beta \in  E_{6(6),(\sC^C)_{\tau\gamma_{\scalebox{0.8}{\sC}}} } $.
Since $ \mathfrak{J}(3,\C^C) $ is decomposed as $ (\mathfrak{J}(3,\C^C))_{\tau\gamma_{{}_{\scalebox{0.8}{\sC}}}} \oplus i(\mathfrak{J}(3,\C^C))_{\tau\gamma_{{}_{\scalebox{0.8}{\sC}}}} $: $ \mathfrak{J}(3,\C^C)=(\mathfrak{J}(3,\C^C))_{\tau\gamma_{{}_{\scalebox{0.8}{\sC}}}} \oplus i(\mathfrak{J}(3,\C^C))_{\tau\gamma_{{}_{\scalebox{0.8}{\sC}}}} $, that is, $ \mathfrak{J}(3,\C^C) $ is the complexification of $ (\mathfrak{J}(3,\C^C))_{\tau\gamma_{{}_{\scalebox{0.8}{\sC}}}} $, we can define an action to $ \mathfrak{J}(3,\C^C) $ of the group $  E_{6(6),(\sC^C)_{\tau\gamma_{\scalebox{0.8}{\sC}}} } $ by
\begin{align*}
\beta X=\beta (X_1+iX_2)=\beta X_1+i\beta X_2,\;\; X:=X_1+iX_2 \in \mathfrak{J}(3,\C^C), X_i \in (\mathfrak{J}(3,\C^C))_{\tau\gamma_{{}_{\scalebox{0.8}{\sC}}}}.
\end{align*}
Then it follows from
\begin{align*}
\beta X \times \beta Y&=\beta (X_1+iX_2) \times \beta (Y_1+iY_2)=(\beta X_1+i\beta X_2) \times (\beta Y_1+i\beta Y_2)
\\
&=(\beta X_1 \times \beta Y_1 - \beta X_2 \times \beta Y_2)+i(\beta X_1 \times \beta Y_2+\beta X_2 \times \beta Y_1)
\\
&=({}^t\!\beta^{-1}(X_1 \times Y_1)-{}^t\!\beta^{-1}(X_2 \times Y_2))+i({}^t\!\beta^{-1}(X_1 \times Y_2)+{}^t\!\beta^{-1}(X_2 \times Y_1))
\\
&={}^t\!\beta^{-1}((X_1 \times Y_1)-(X_2 \times Y_2))+i((X_1 \times Y_2)+(X_2 \times Y_1))
\\
&={}^t\!\beta^{-1}((X_1+iX_2)\times (Y_1+iY_2))
\\
&={}^t\!\beta^{-1}(X \times Y)
\end{align*}
that $ \beta \in (E_{6,\sC})^C $. Moreover, it is easy to verify that $ (\tau\gamma_{{}_{\scalebox{0.8}{\sC}}})\beta=\beta(\gamma_{{}_{\scalebox{0.8}{\sC}}}\tau) $. Indeed, for $ X:=X_1+iX_2 \in (\mathfrak{J}(3,\C^C))_{\tau\gamma_{{}_{\scalebox{0.8}{\sC}}}} \oplus i(\mathfrak{J}(3,\C^C))_{\tau\gamma_{{}_{\scalebox{0.8}{\sC}}}}=\mathfrak{J}(3,\C^C) $, it follows that
\begin{align*}
(\tau\gamma_{{}_{\scalebox{0.8}{\sC}}})\beta(\gamma_{{}_{\scalebox{0.8}{\sC}}}\tau)X&=(\tau\gamma_{{}_{\scalebox{0.8}{\sC}}})\beta(\gamma_{{}_{\scalebox{0.8}{\sC}}}\tau)(X_1+iX_2)=(\tau\gamma_{{}_{\scalebox{0.8}{\sC}}})\beta(X_1-iX_2)
\\
&=(\tau\gamma_{{}_{\scalebox{0.8}{\sC}}})(\beta X_1-i\beta X_2)=\beta X_1+i\beta X_2=\beta(X_1+iX_2)
\\
&=\beta X,
\end{align*}
that is, $ (\tau\gamma_{{}_{\scalebox{0.8}{\sC}}})\beta(\gamma_{{}_{\scalebox{0.8}{\sC}}}\tau)=\beta $.
With above, we have $ \beta \in ((E_{6,\sC})^C)^{\tau\gamma_{{}_{\scalebox{0.8}{\sC}}}} $.

Thus we have the required result $ ((E_{6,\sC})^C)^{\tau\gamma_{{}_{\scalebox{0.8}{\sC}}}}=E_{6(6),(\sC^C)_{\tau\gamma_{\scalebox{0.8}{\sC}}} } $.

Therefore, together with Proposition \ref{proposition 4.3.1}, we have the isomorphism
\begin{align*}
((E_{6,\sC})^C)^{\tau\gamma_{{}_{\scalebox{0.8}{\sC}}}} \cong E_{6(6),\sC'}.
\end{align*}
\end{proof}

We prepare a more little. Let the group $ SU(3,\C'):=\{A \in M(3,\C')\,|\,AA^*=E, \det\,A=1\} $. Moreover, we define a subgroup $ (SU(3,\C^C))^{\tau\gamma_{{}_{\scalebox{0.8}{\sC}}}} $ of $ SU(3,\C^C) $ by
\begin{align*}
(SU(3,\C^C))^{\tau\gamma_{{}_{\scalebox{0.8}{\sC}}}}:&=\left\lbrace A \in SU(3,\C^C) \relmiddle{|} \tau\gamma_{{}_{\scalebox{0.8}{\sC}}} A=A \right\rbrace
\\
&=\left\lbrace A \in M(3,(\C^C)_{\tau\gamma_{{}_{\scalebox{0.8}{\sC}}}})\relmiddle{|}A^*A=E, \det\,A=1 \right\rbrace.
\end{align*}

Then we prove the following proposition.

\begin{proposition}\label{proposition 4.3.3}
The group $ (SU(3,\C^C))^{\tau\gamma_{{}_{\scalebox{0.8}{\sC}}}} $ is isomorphic to the group $ SU(3,\C') ${\rm :}$ (SU(3,\C^C))^{\tau\gamma_{{}_{\scalebox{0.8}{\sC}}}} \cong SU(3,\C') $.

\end{proposition}
\begin{proof}
Then the correspondence
\begin{align*}
(\C^C)_{\tau\gamma_{\scalebox{0.8}{\sC}}} \ni x+y i\i
\underset{l}{\mapsto} x+y\i' \in \C'
\end{align*}
gives an isomorphism $  (SU(3,\C^C))^{\tau\gamma_{{}_{\scalebox{0.8}{\sC}}}} \cong SU(3,\C') $.

\end{proof}

Let the discrete group $ \Z_2=\{1,\varepsilon\} $. This group acts on $ (SU(3,\C^C))^{\tau\gamma_{{}_{\scalebox{0.8}{\sC}}}} \times (SU(3,\C^C))^{\tau\gamma_{{}_{\scalebox{0.8}{\sC}}}} $ by
\begin{align*}
1(A,B)=(A,B),\quad \varepsilon(A,B)=(\ov{B},\ov{A}),
\end{align*}
and then $ ((SU(3,\C^C))^{\tau\gamma_{{}_{\scalebox{0.8}{\sC}}}} \times (SU(3,\C^C))^{\tau\gamma_{{}_{\scalebox{0.8}{\sC}}}}) \rtimes \Z_2 $ be the semi-direct product of $ (SU(3,\C^C))^{\tau\gamma_{{}_{\scalebox{0.8}{\sC}}}} \times (SU(3,\C^C))^{\tau\gamma_{{}_{\scalebox{0.8}{\sC}}}} $ and $ \Z_2 $ with the same multiplication of $ (SU(3,\C^C) \times SU(3,\C^C)) \rtimes \Z_2 $.

Now, we determine the structure of the group $ E_{6(6),\sC'} $.

\begin{theorem}\label{theorem 4.3.4}
The group $ E_{6(6),\sC'} $ is isomorphic to the semi-direct product of the groups $ SU(3,\C') \times SU(3,\C') $ and $ \Z_2, \Z_2=\{1,\varepsilon \} ${\rm :} $ E_{6(6),\sC'} \cong (SU(3,\C') \times SU(3,\C')) \rtimes \Z_2 $.
\end{theorem}
\begin{proof}
Let the group $ E_{6(6),\sC'} $ as the group $ ((E_{6,\sC})^C)^{\tau\gamma_{{}_{\scalebox{0.8}{\sC}}}} $ (Theorem \ref{theorem 4.3.2}) and the group $ SU(3,\C') $ as $ (SU(3,\C^C))^{\tau\gamma_{{}_{\scalebox{0.8}{\sC}}}} $ (Proposition 4.3.3).  We define a mapping $ f_{6,\tau\gamma_{{}_{\scalebox{0.8}{\sC}}}}\!: ((SU(3,\C^C))^{\tau\gamma_{{}_{\scalebox{0.8}{\sC}}}} \times (SU(3,\C^C))^{\tau\gamma_{{}_{\scalebox{0.8}{\sC}}}}) \allowbreak \rtimes \{1,\varepsilon\} \to ((E_{6,\sC})^C)^{\tau\gamma_{{}_{\scalebox{0.8}{\sC}}}} $ by
\begin{align*}
f_{6,\tau\gamma_{{}_{\scalebox{0.8}{\sC}}}}((A,B),1)X&=h(A,B)Xh(A,B)^*,
\\
f_{6,\tau\gamma_{{}_{\scalebox{0.8}{\sC}}}}((A,B),\varepsilon)X&=h(A,B)\ov{X}h(A,B)^*, \;\; X \in \mathfrak{J}(3,\C^C).
\end{align*}
Note that this mapping is the restriction of the mapping $ f_{6,\sC^C} $.
First, we will prove that $ f_{6,\tau\lambda} $ is well-defined and a homomorphism. Note that this mapping is the restriction of the mapping $ f_{6,\sC^C} $.
First, we will prove that $ f_{6,\tau\gamma_{{}_{\scalebox{0.8}{\sC}}}} $ is well-defined and a homomorphism. Since $ f_{6,\tau\gamma_{{}_{\scalebox{0.8}{\sC}}}} $ is the restriction of the mapping $ f_{6,\sC^C} $, it is clear that $ f_{6,\tau\gamma_{{}_{\scalebox{0.8}{\sC}}}}((A,B),1),f_{6,\tau\gamma_{{}_{\scalebox{0.8}{\sC}}}}((A,B),\varepsilon) \allowbreak \in (E_{6,\sC})^C $ and $ f_{6,\tau\gamma_{{}_{\scalebox{0.8}{\sC}}}} $ is a homomorphism.
By using Lemma \ref{lemma 4.1.8} (2), (3), we have
\begin{align*}
(\tau\gamma_{{}_{\scalebox{0.8}{\sC}}})f_{6,\tau\gamma_{{}_{\scalebox{0.8}{\sC}}}}((A,B),1)(\gamma_{{}_{\scalebox{0.8}{\sC}}}\tau)&=f_{6,\tau\gamma_{{}_{\scalebox{0.8}{\sC}}}}((\tau\gamma_{{}_{\scalebox{0.8}{\sC}}} A,\tau\gamma_{{}_{\scalebox{0.8}{\sC}}} B),1)=f_{6,\tau\gamma_{{}_{\scalebox{0.8}{\sC}}}}((A,B),1),
\end{align*}
so that $ f_{6,\tau\gamma_{{}_{\scalebox{0.8}{\sC}}}}((A,B),1) \in ((E_{6,\sC})^C)^{\tau\gamma_{{}_{\scalebox{0.8}{\sC}}}} $. As in that above, by using Lemma \ref{lemma 4.1.8} (6), (7), it is clear $ (\tau\gamma_{{}_{\scalebox{0.8}{\sC}}}) f_{6,\tau\gamma_{{}_{\scalebox{0.8}{\sC}}}}((E,E),\varepsilon)(\gamma_{{}_{\scalebox{0.8}{\sC}}}\tau)=f_{6,\tau\gamma_{{}_{\scalebox{0.8}{\sC}}}}((E,E),\varepsilon) $, so that $ f_{6,\tau\gamma_{{}_{\scalebox{0.8}{\sC}}}}((E,E),\varepsilon) \in ((E_{6,\sC})^C)^{\tau\gamma_{{}_{\scalebox{0.8}{\sC}}}} $. Hence, since $ f_{6,\tau\gamma_{{}_{\scalebox{0.8}{\sC}}}}((A,B),\varepsilon)=f_{6,\tau\gamma_{{}_{\scalebox{0.8}{\sC}}}}((A,B),1)f_{6,\tau\gamma_{{}_{\scalebox{0.8}{\sC}}}}((E,E),\varepsilon) $, we see $ f_{6,\tau\gamma_{{}_{\scalebox{0.8}{\sC}}}}((A,B),\varepsilon) \in ((E_{6,\sC})^C)^{\tau\gamma_{{}_{\scalebox{0.8}{\sC}}}} $. With above, the proof of well-defined is proved.

Next, we will prove that $ f_{6.\tau\gamma_{{}_{\scalebox{0.8}{\sC}}}} $ is surjective. Let $ \alpha \in ((E_{6,\sC})^C)^{\tau\gamma_{{}_{\scalebox{0.8}{\sC}}}} \subset (E_{6,\sC})^C $. Then there exists $ ((P,Q),1) \in (SU(3,\C^C) \times SU(3,\C^C)) \rtimes \Z_2 $ such that $ \alpha=f_{6,\sC^C}((P,Q),1) $ or there exists $ ((P,Q),\varepsilon) \in (SU(3,\C^C) \times SU(3,\C^C)) \rtimes \Z_2 $ such that $ \alpha=f_{6,\sC^C}((P,Q),\varepsilon) $ (Theorem \ref{theorem 4.1.7}). Moreover, $ \alpha $ satisfies the condition $ (\tau\gamma_{{}_{\scalebox{0.8}{\sC}}})\alpha(\gamma_{{}_{\scalebox{0.8}{\sC}}}\tau)=\alpha $, that is, $ (\tau\gamma_{{}_{\scalebox{0.8}{\sC}}})f_{6,\tau\gamma_{{}_{\scalebox{0.8}{\sC}}}}((A,B),1)(\gamma_{{}_{\scalebox{0.8}{\sC}}}\tau)=f_{6,\tau\gamma_{{}_{\scalebox{0.8}{\sC}}}}((A,B),1) $ or $ (\tau\gamma_{{}_{\scalebox{0.8}{\sC}}})f_{6,\tau\gamma_{{}_{\scalebox{0.8}{\sC}}}}((A,B),\varepsilon)(\gamma_{{}_{\scalebox{0.8}{\sC}}}\tau)=f_{6,\tau\gamma_{{}_{\scalebox{0.8}{\sC}}}}((A,B),\varepsilon) $.
Since it follows from Lemma \ref{lemma 4.1.8} (2), (3), (6), (7) that
\begin{align*}
(\tau\gamma_{{}_{\scalebox{0.8}{\sC}}})f_{6,\tau\gamma_{{}_{\scalebox{0.8}{\sC}}}}((A,B),1)(\gamma_{{}_{\scalebox{0.8}{\sC}}}\tau)&=f_{6,\tau\gamma_{{}_{\scalebox{0.8}{\sC}}}}((\tau\gamma_{{}_{\scalebox{0.8}{\sC}}} A,\tau\gamma_{{}_{\scalebox{0.8}{\sC}}} B),1)
\\
&\;\text{or}
\\
(\tau\gamma_{{}_{\scalebox{0.8}{\sC}}})f_{6,\tau\gamma_{{}_{\scalebox{0.8}{\sC}}}}((A,B),\varepsilon)(\gamma_{{}_{\scalebox{0.8}{\sC}}}\tau)&=f_{6,\tau\gamma_{{}_{\scalebox{0.8}{\sC}}}}((\tau\gamma_{{}_{\scalebox{0.8}{\sC}}} A,\tau\gamma_{{}_{\scalebox{0.8}{\sC}}} B),\varepsilon),
\end{align*}
we have the following
\begin{align*}
\left\lbrace
\begin{array}{l}
\tau\gamma_{{}_{\scalebox{0.8}{\sC}}} P=P \\
\tau\gamma_{{}_{\scalebox{0.8}{\sC}}} Q=Q,
\end{array}
\right. \quad
\left\lbrace
\begin{array}{l}
\tau\gamma_{{}_{\scalebox{0.8}{\sC}}} P=\omega P \\
\tau\gamma_{{}_{\scalebox{0.8}{\sC}}} Q=\omega Q
\end{array}
\right. \quad \text{or} \quad
\left\lbrace
\begin{array}{l}
\tau\gamma_{{}_{\scalebox{0.8}{\sC}}} P=\omega ^2 P \\
\tau\gamma_{{}_{\scalebox{0.8}{\sC}}} Q=\omega ^2 Q.
\end{array}
\right.
\end{align*}
In the first case, we have $ P,Q \in (SU(3,\C^C))^{\tau\gamma_{{}_{\scalebox{0.8}{\sC}}}} $. In the others cases, since we have $ P=Q=0 $, those are impossible. The proof of surjective is completed.

Finally, we will determine $ \Ker\,f_{6,\tau\gamma_{{}_{\scalebox{0.8}{\sC}}}} $. Since $ f_{6,\tau\gamma_{{}_{\scalebox{0.8}{\sC}}}} $ is the restriction of the mapping $ f_{6,\sC^C} $, we have
\begin{align*}
\Ker\,f_{6,\tau\gamma_{{}_{\scalebox{0.8}{\sC}}}} \subset \Ker\,f_{6,\sC^C}=\left\lbrace ((E,E),1), (({\boldsymbol\omega} E,{\boldsymbol\omega} E),1), (({\boldsymbol\omega}^2 E,{\boldsymbol\omega}^2 E),1)\right\rbrace.
\end{align*}
However, since $ \boldsymbol{\omega}E, \boldsymbol{\omega^2}E \notin (SU(3,\C^C))^{\tau\gamma_{{}_{\scalebox{0.8}{\sC}}}} $, we obtain
$ \Ker\,f_{6,\tau\gamma_{{}_{\scalebox{0.8}{\sC}}}}=\{((E,E),1)\} $.

Thus we have the isomorphism
\begin{align*}
((E_{6,\sC^C})^C)^{\tau\gamma_{{}_{\scalebox{0.8}{\sC}}}} \cong ((SU(3,\C^C))^{\tau\gamma_{{}_{\scalebox{0.8}{\sC}}}} \times (SU(3,\C^C))^{\tau\gamma_{{}_{\scalebox{0.8}{\sC}}}} ) \rtimes \Z_2.
\end{align*}

Therefore, from Theorem \ref{theorem 4.3.2} and Proposition \ref{proposition 4.3.3}, we have the required isomorphism
\begin{align*}
E_{6(6),\sC'} \cong (SU(3,\C') \times SU(3,\C')) \rtimes \Z_2.
\end{align*}
\end{proof}

\subsection{The group $ E_{6(-14),\sC} $}

As mentioned in the surface of this section, since $ (E_{6,\sC})^C $ has an involutive automorphism $ \tilde{\tau\lambda\sigma} $, we can define a subgroup $ ((E_{6,\sC})^C)^{\tau\lambda\sigma} $ of $ (E_{6,\sC})^C $:
\begin{align*}
((E_{6,\sC})^C)^{\tau\lambda\sigma}:=\left\lbrace \alpha \in (E_{6,\sC})^C \relmiddle{|} \tilde{\tau\lambda\sigma}(\alpha)=\alpha \right\rbrace.
\end{align*}

First, we prove the following theorem.

\begin{theorem}\label{theorem 4.4.1}
The group $ ((E_{6,\sC})^C)^{\tau\lambda\sigma} $ coincides with the group $ E_{6(-14),\sC} ${\rm :} $ ((E_{6,\sC})^C)^{\tau\lambda\sigma}=E_{6(-14),\sC} $.
\end{theorem}
\begin{proof}
Let $ \alpha \in ((E_{6,\sC})^C)^{\tau\lambda\sigma} $. It follows from $ (\tau\sigma){}^t\!\alpha^{-1}(\sigma\tau)=\alpha $ that
\begin{align*}
\langle \alpha X, \alpha Y \rangle_\sigma&=(\tau\sigma\alpha X,\alpha Y)=({}^t\!\alpha^{-1}(\sigma\tau)X,\alpha Y)=((\sigma\tau)X,\alpha^{-1
}\alpha Y)=(\tau\sigma X,Y)
\\
&=\langle X ,Y \rangle_\sigma.
\end{align*}
Hence we see $ \alpha \in E_{6(-14),\sC} $. Conversely, let $ \beta \in E_{6(-14),\sC} $. Then it follows from $ \langle \beta X, \beta Y \rangle_\sigma=\langle X ,Y \rangle_\sigma $ that
\begin{align*}
(\tau\sigma X,Y)=\langle X ,Y \rangle_\sigma=\langle \beta X, \beta Y \rangle_\sigma=(\tau\sigma\beta X,\beta Y)=({}^t\!\beta\tau\sigma \alpha X,Y)
\end{align*}
that $ \tau\sigma={}^t\!\beta\tau\sigma \beta $, that is, $ (\tau\sigma){}^t\!\beta^{-1}(\sigma\tau)=\beta $, so that $ \beta \in ((E_{6,\sC})^C)^{\tau\lambda\sigma} $.

With above, the proof of this theorem is completed.
\end{proof}

Let the discrete group $ \Z_2=\{1,\varepsilon\} $. This group acts on $ SU(1,2) \times SU(1,2) $ by
\begin{align*}
1(A,B)=(A,B),\quad \varepsilon(A,B)=(\ov{B},\ov{A}),
\end{align*}
and then $ (SU(1,2) \times SU(1,2)) \rtimes \Z_2 $ be the semi-direct product of $ SU(1,2) \times SU(1,2) $ and $ \Z_2 $ with the same multiplication of $ (SU(3,\C^C) \times SU(3,\C^C)) \rtimes \Z_2 $.

Now, we determine the structure of the group $ E_{6(-14),\sC} $.

\begin{theorem}\label{theorem 4.4.2}
The group $ E_{6(-14),\sC} $ is isomorphic to the semi-direct product of the groups $ (SU(1,2) \times SU(1,2))/\Z_3 $ and $ \Z_2, \Z_3=\{(E,E),(\boldsymbol{\omega} E, \boldsymbol{\omega} E), (\boldsymbol{\omega}^2E,\boldsymbol{\omega}^2E)\}, \Z_2=\{1,\varepsilon\} ${\rm :} $ E_{6(-14),\sC} \cong (SU(1,2) \times SU(1,2))/\Z_3 \allowbreak \rtimes \Z_2 $.
\end{theorem}
\begin{proof}
Let $ E_{6(-14),\sC} $ as the group $ ((E_{6,\sC})^C)^{\tau\lambda\sigma} $ (Theorem \ref{theorem 4.4.1}). We define a mapping $ f_{6,\tau\lambda\sigma}: (SU(1,2) \times SU(1,2)) \rtimes \{1,\varepsilon\} \to ((E_{6,\sC})^C)^{\tau\lambda\sigma} $ by
\begin{align*}
f_{6,\tau\lambda\sigma}((A,B),1)X&=h(A,B)Xh(A,B)^*,
\\
f_{6,\tau\lambda\sigma}((A,B),\varepsilon)X&=h(A,B)\ov{X}h(A,B)^*, \;\; X \in \mathfrak{J}(3,\C^C).
\end{align*}
Note that this mapping is the restriction of the mapping $ f_{6,\sC^C} $.
First, we will prove that $ f_{6,\tau\lambda\sigma} $ is well-defined and a homomorphism. Since $ f_{6,\tau\lambda\sigma} $ is the restriction of the mapping $ f_{6,\sC^C} $, it is clear that $ f_{6,\tau\lambda\sigma}((A,B),1),f_{6,\tau\lambda\sigma}((A,B),\varepsilon) \allowbreak \in (E_{6,\sC})^C $ and $ f_{6,\tau\lambda\sigma} $ is a homomorphism.
By using Lemma \ref{lemma 4.1.8} (1), (2), (4), we have $ \tau {}^t(f_{6,\tau\lambda}((A,B),1))^{-1}\tau=f_{6,\tau\lambda}((A,B),1) $, so that $ f_{6,\tau\lambda}((A,B),1) \in ((E_{6,\sC})^C)^{\tau\lambda} $. As in that above, by using Lemma \ref{lemma 4.1.8} (5), (6), (8), we have $ f_{6,\tau\lambda}((E,E),\varepsilon) \in ((E_{6,\sC})^C)^{\tau\lambda} $, so that since $ f_{6,\tau\lambda}((A,B),\varepsilon)=f_{6,\tau\lambda}((A,B),1)f_{6,\tau\lambda}((E,E),\varepsilon) $, we see $ f_{6,\tau\lambda}((A,B),\varepsilon) \in ((E_{6,\sC})^C)^{\tau\lambda} $. With above, the proof of well-defined is proved.

Next, we will prove that $ f_{6,\tau\lambda\sigma} $ is surjective. Let $ \alpha \in E_{6,\tau\lambda\sigma}=((E_{6,\sC})^C)^{\tau\lambda\sigma} \subset (E_{6,\sC})^C $. Then there exists $ ((P,Q),1) \in (SU(3,\C^C) \times SU(3,\C^C)) \rtimes \Z_2 $ such that $ \alpha=f_{6,\sC^C}((P,Q),1) $ or there exists $ ((P,Q),\varepsilon) \in (SU(3,\C^C) \times SU(3,\C^C)) \rtimes \Z_2 $ such that $ \alpha=f_{6,\sC^C}((P,Q),\varepsilon) $ (Theorem \ref{theorem 4.1.7}). Moreover $ \alpha $ satisfies the condition $ \tau {}^t\!\alpha^{-1} \tau=\alpha $, that is, $ \tau{}^t(f_{6,\sC^C}((P,Q),1))^{-1}\tau=f_{6,\sC^C}((P,Q),1) $ or $ \tau{}^t(f_{6,\sC^C}((P,Q),\varepsilon))^{-1}\tau=f_{6,\sC^C}((P,Q),\varepsilon) $. Since it follows from Lemma \ref{lemma 4.1.8} (1), (2), (4), (5), (6), (8) that
\begin{align*}
\tau{}^t(f_{6,\sC^C}((P,Q),1))^{-1}\tau&=f_{6,\sC^C}((\tau P,\tau Q),1)
\\
&\text{or}
\\
\tau{}^t(f_{6,\sC^C}((P,Q),\varepsilon))^{-1}\tau&=f_{6,\sC^C}((\tau P,\tau Q),\varepsilon),
\end{align*}
we have the following
\begin{align*}
\left\lbrace
\begin{array}{l}
\tau P=P \\
\tau Q=Q,
\end{array}
\right. \quad
\left\lbrace
\begin{array}{l}
\tau P=\omega P \\
\tau Q=\omega Q
\end{array}
\right. \quad \text{or} \quad
\left\lbrace
\begin{array}{l}
\tau P=\omega ^2 P \\
\tau Q=\omega ^2 Q.
\end{array}
\right.
\end{align*}
In the first case, we have $ P,Q \in SU(1,2) $. In the others cases, since we have $ P=Q=0 $, those are impossible. The proof of surjective is completed.

Finally, we will determine $ \Ker\,f_{6,\tau\lambda\sigma} $. Since $ f_{6,\tau\lambda\sigma} $ is the restriction of the mapping $ f_{6,\sC^C} $, it is easy to obtain $ \Ker\,f_{6,\tau\lambda\sigma}=\Ker\,f_{6,\sC^C}=\left\lbrace ((E,E),1), (({\boldsymbol\omega} E,{\boldsymbol\omega} E),1), (({\boldsymbol\omega}^2 E,{\boldsymbol\omega}^2 E),1)\right\rbrace \cong (\Z_3,1) $.

Therefore we have the required isomorphism
\begin{align*}
E_{6,\tau\lambda\sigma} \cong (SU(1,2) \times SU(1,2))/\Z_3 \rtimes \Z_2.
\end{align*}
\end{proof}

\subsection{The group $ E_{6(-26),\sC} $}

As mentioned in the surface of this section, since $ (E_{6,\sC})^C $ has an involutive automorphism $ \tilde{\tau} $, we can define a subgroup $ ((E_{6,\sC})^C)^{\tau} $ of $ (E_{6,\sC})^C $:
\begin{align*}
((E_{6,\sC})^C)^{\tau}:=\left\lbrace \alpha \in (E_{6,\sC})^C \relmiddle{|} \tilde{\tau}(\alpha)=\alpha \right\rbrace.
\end{align*}

First, we prove the following theorem.

\begin{theorem}\label{theorem 4.5.1}
The group $ ((E_{6,\sC})^C)^{\tau} $ coincides with the group $ E_{6(-26),\sC} ${\rm :} $ ((E_{6,\sC})^C)^{\tau}=E_{6(-26),\sC} $.
\end{theorem}
\begin{proof}
Let $ \alpha \in ((E_{6,\sC})^C)^{\tau} $. Then it follows from $ \tau\alpha=\alpha\tau $ that $ \alpha X=\alpha (\tau X)=\tau (\alpha X), X \in \mathfrak{J}(3,\C) $, that is, $ \alpha X \in \mathfrak{J}(3,\C) $, so that $ \alpha $ indices an $ \R $-linear isomorphism of $ \mathfrak{J}(3,\C) $. Hence we see $ \alpha \in E_{6(-26),\sC} $. Conversely, let $ \beta \in E_{6(-26),\sC} $. Then we define an action of $ \beta $ to $ \mathfrak{J}(3,\C^C) $ by
\begin{align*}
\beta X=\beta(X_1+iX_2)=\beta X_1+ i\beta X_2,\;\; X \in \mathfrak{J}(3,\C^C), X_i \in \mathfrak{J}(3,\C).
\end{align*}
Hence $ \beta $ induces a $ C $-linear isomorphism of $ \mathfrak{J}(3,\C^C) $. Moreover, it follows that
\begin{align*}
\beta X \times \beta Y&=\beta(X_1+iX_2) \times \beta(Y_1+iY_2)=(\beta X_1+i\beta X_2) \times (\beta Y_1+i\beta Y_2)
\\
&=(\beta X_1 \times \beta Y_1-\beta X_2 \times \beta Y_2)+i(\beta X_1 \times \beta Y_2+\beta X_2 \times \beta Y_1)
\\
&=({}^t\!\beta^{-1}(X_1 \times Y_1)-{}^t\!\beta^{-1}(X_2 \times Y_2))+i({}^t\!\beta^{-1}(X_1 \times Y_2)+{}^t\!\beta^{-1}(X_2 \times Y_1))
\\
&={}^t\!\beta^{-1}((X_1 \times Y_1-X_2 \times Y_2)+i(X_1 \times Y_2+X_2 \times Y_1))
\\
&={}^t\!\beta^{-1}((X_1+iX_2)\times (Y_1+iY_2))
\\
&={}^t\!\beta^{-1}(X \times Y),
\end{align*}
so that $ \beta \in (E_{6,\sC})^C $. In addition, we have
\begin{align*}
\tau\beta X&=\tau\beta(X_1+iX_2)=\tau(\beta X_1+i\beta X_2)=\beta X_1-i\beta X_2=\beta(X_1-iX_2)=\beta\tau(X_1+iX_2)
\\
&=\beta\tau X,
\end{align*}
that is, $ \tau\beta=\beta\tau $. Hence we see $ \beta \in ((E_{6,\sC})^C)^{\tau} $.

With above, the proof of this theorem is completed.
\end{proof}

Let the discrete group $ \Z_2=\{1,\varepsilon\} $. This group acts on $ SU(3,\C^C) $ by
\begin{align*}
1 A=A,\quad \varepsilon A=\ov{A},
\end{align*}
and then $ SU(3,\C^C) \rtimes \Z_2 $ be the semi-direct product of $ SU(3,\C^C)$ and $ \Z_2 $ with the multiplication
\begin{align*}
(A_1,1)(A_2,1)&=(A_1A_2,1), \quad
(A_1,1)(A_2,\varepsilon)=(A_1A_2,\varepsilon),
\\
(A_1,\varepsilon)(A_2,1)&=(A_1\ov{A}_2,\varepsilon), \quad
(A_1,\varepsilon)(A_2,\varepsilon)=(A_1\ov{A}_2,1).
\end{align*}

Now, we determine the structure of the group $ E_{6(-26),\sC} $.

\begin{theorem}\label{theorem 4.5.2}
The group $ E_{6(-26),\sC} $ is isomorphic to the semi-direct product of the groups $ SU(3,\C^C)/\Z_3 $ and $ \Z_2, \Z_3=\{ E,\boldsymbol{\omega} E, \boldsymbol{\omega}^2E \}, \Z_2=\{1,\varepsilon\} ${\rm :} $ E_{6(-26),\sC} \cong SU(3,\C^C)/\Z_3 \rtimes \Z_2 $.
\end{theorem}
\begin{proof}
Let the group $ E_{6(-26),\sC} $ as the group $ ((E_{6,\sC})^C)^{\tau} $ (Theorem \ref{theorem 4.5.1}). Then we define a mapping $ f_{6,\tau}: SU(3,\C^C) \rtimes \{1,\varepsilon \} \to ((E_{6,\sC})^C)^{\tau} $ by
\begin{align*}
f_{6,\tau}(A,1)X&=h(A, \tau A)Xh(A,\tau A)^*,
\\
f_{6,\tau}(A,\varepsilon)X&=h(A, \tau A)\ov{X}h(A,\tau A)^*, \;\; X \in \mathfrak{J}(3,\C^C).
\end{align*}
Note that if $ A \in SU(3,\C^C) $, then it is easy to verify $ \tau A \in SU(3,\C^C)  $, so that the mapping $ f_{6,\tau} $ is the restriction of the mapping $ f_{6,\sC^C} $.

First, we will prove that $ E_{6,\tau} $ is well-defined and a homomorphism. As mentioned above, since $ f_{6,\tau} $ is the restriction of the mapping $ f_{6,\sC^C} $, immediately we see that $ f_{6,\tau}(A,1), f_{6,\tau}(A,\varepsilon) \in (E_{6,\sC^C})^C $ and $ f_{6,\tau} $ is a homomorphism. Subsequently, it follows from Lemma \ref{lemma 4.1.8} (2) that
\begin{align*}
\tau f_{6,\tau}(A,1)\tau X&=\tau(h(A, \tau A)(\tau X)h(A,\tau A)^*)=\tau h(A, \tau A)X \tau h(A,\tau A)^*
\\
&=h(A, \tau A)Xh(A,\tau A)^*
\\
&=f_{6,\tau}(A,1)X.
\end{align*}
Hence we see $ f_{6,\tau}(A,1) \in ((E_{6,\sC^C})^C)^\tau $. Moreover, from Lemma \ref{lemma 4.1.8} (6), it is clear $ \tau f_{6,\tau}(E,\varepsilon) \tau=f_{6,\tau}(E,\varepsilon) $, so that since $  f_{6,\tau}(A,\varepsilon)= f_{6,\tau}(A,1) f_{6,\tau}(E,\varepsilon) $, we see $  f_{6,\tau}(A,\varepsilon) \in ((E_{6,\sC^C})^C)^\tau $.
with above, the proof of well-defined is proved.

Next, we will prove that $ f_{6.\tau} $ is surjective. Let $ \alpha \in E_{6(-26), \sC}=((E_{6,\sC})^C)^\tau \subset (E_{6,\sC})^C $. Then there exists $ ((P,Q),1) \in (SU(3,\C^C) \times SU(3,\C^C)) \rtimes \Z_2 $ such that $ \alpha=f_{6,\sC^C}((P,Q),1) $ or there exists $ ((P,Q),\varepsilon) \in (SU(3,\C^C) \times SU(3,\C^C)) \rtimes \Z_2 $ such that $ \alpha=f_{6,\sC^C}((P,Q),\varepsilon) $ (Theorem \ref{theorem 4.1.7}). Moreover $ \alpha $ satisfies the condition $ \tau\alpha\tau=\alpha $, that is, $ \tau f_{6,\sC^C}((P,Q),1)\tau=f_{6,\sC^C}((P,Q),1) $ or $ \tau f_{6,\sC^C}((P,Q),\varepsilon)\tau \allowbreak =f_{6,\sC^C}((P,Q),\varepsilon) $. Since it follows from Lemma \ref{lemma 4.1.8} (2), (6) that
\begin{align*}
\tau f_{6,\sC^C}((P,Q),1)\tau&=f_{6,\sC^C}((\tau Q,\tau P),1)
\\
&\text{or}
\\
\tau f_{6,\sC^C}((P,Q),\varepsilon) \tau&=f_{6,\sC^C}((\tau Q,\tau P),\varepsilon),
\end{align*}
we have the following
\begin{align*}
\left\lbrace
\begin{array}{l}
\tau P=Q \\
\tau Q=P,
\end{array}
\right. \quad
\left\lbrace
\begin{array}{l}
\tau P=\omega Q \\
\tau Q=\omega P
\end{array}
\right. \quad \text{or} \quad
\left\lbrace
\begin{array}{l}
\tau P=\omega ^2 Q \\
\tau Q=\omega ^2 P.
\end{array}
\right.
\end{align*}
In the first case, we have $ Q=\tau P $, so that there exists $ A \in SU(3,\C^C) $ such that $ \alpha=f_{6,\sC^C}((A,\tau A), 1)=f_{6,\tau}(A,1) $ or  $ \alpha=f_{6,\sC^C}((A,\tau A), \varepsilon)=f_{6,\tau}(A,\varepsilon) $.
In the others cases, since we have $ P=Q=0 $, those are impossible.
The proof of surjective is completed.

Finally, we will determine $ \Ker\,f_{6,\tau} $. Since $ f_{6,\tau} $ is the restriction of the mapping $ f_{6,\sC^C} $, it is easy to obtain $ \Ker\,f_{6,\tau}=\left\lbrace (E,1), ({\boldsymbol\omega} E,1), ({\boldsymbol\omega}^2 E,1)\right\rbrace \cong (\Z_3,1) $.

Therefore, from Theorem \ref{theorem 4.5.1}, we have the required isomorphism
\begin{align*}
E_{6(-26),\sC} \cong SU(3,\C^C)/\Z_3 \rtimes \Z_2.
\end{align*}
Note that we regard the group $ SU(3,\C^C) $ as real Lie groups.
\end{proof}

\section{The complex Lie group $ (E_{6,\sH})^C $ and its real forms}

We define the group $ (E_{6,\sH})^C $ by
\begin{align*}
(E_{6,\sH})^C:&=\left\lbrace \alpha \in \Iso_{C}(\mathfrak{J}(3,\H^C))\relmiddle{|} \det(\alpha X)=\det\,X \right\rbrace
\\
&=\left\lbrace \alpha \in \Iso_{C}(\mathfrak{J}(3,\H^C))\relmiddle{|} (\alpha X,\alpha Y, \alpha Z)=(X,Y,Z) \right\rbrace
\\
&=\left\lbrace \alpha \in \Iso_{C}(\mathfrak{J}(3,\H^C))\relmiddle{|} \alpha X \times \alpha Y={}^t\!\alpha^{-1}(X \times Y) \right\rbrace.
\end{align*}
and its real forms are defined as follows:
\begin{align*}
E_{6,\sH}&:=\left\lbrace \alpha \in \Iso_{C}(\mathfrak{J}(3,\H^C))\relmiddle{|} \det(\alpha X)=\det\,X, \langle \alpha X, \alpha Y \rangle=\langle X,Y \rangle \right\rbrace,
\\
E_{6(6),\sH'}&:=\left\lbrace \alpha \in \Iso_{\sR}(\mathfrak{J}(3,\H'))\relmiddle{|} \det(\alpha X)=\det\,X \right\rbrace,
\\
E_{6(-14),\sH}&:=\left\lbrace \alpha \in \Iso_{C}(\mathfrak{J}(3,\H^C))\relmiddle{|} \det(\alpha X)=\det\,X, \langle \alpha X, \alpha Y \rangle_\sigma=\langle X,Y \rangle_\sigma \right\rbrace,
\\
E_{6(-26),\sH}&:=\left\lbrace \alpha \in \Iso_{\sR}(\mathfrak{J}(3,\H))\relmiddle{|} \det(\alpha X)=\det\,X \right\rbrace,
\end{align*}
where $ \langle X,Y \rangle=(\tau X,Y), \langle X,Y \rangle_\sigma=(\tau\sigma X,Y), \H'=\C' \oplus \C'\j, \j^2 =-1$,
and since $ E_{6(2),\sH}=E_{6,\sH}$, the definition of $ E_{6(2),\sH} $ is omitted.
\vspace{1mm}

As in $ (E_{6,\sC})^C $, the group $ (E_{6,\sH})^C $ has involutive automorphism $ \tilde{\tau\lambda}, \tilde{\tau\gamma_{\scalebox{0.8}{\sC}}}, \tilde{\tau\lambda\sigma} $ and $ \tilde{\tau} $.

\subsection{The group $ (E_{6,\sH})^C $}

The structure of the group $ (E_{6,\sH})^C $ has been determined by Ichiro Yokota and his school's members (\cite[Proposition 3.5.3]{iy7}). We state its result as theorem below with short proof.

\begin{theorem}{\rm \cite[Proposition 3.5.3]{iy7}}\label{theprem 5.1.1}
The group $ (E_{6,\sH})^C $ is isomorphic to the group $ SU^*(6,\C^C)/\Z_2, \Z_2 \allowbreak =\{E,-E \} ${\rm :} $ (E_{6,\sH})^C \cong SU^*(6,\C^C)/\Z_2 $.
\end{theorem}
\begin{proof}
We define a mapping $ f_{6,\sH^C}: SU^*(6,\C^C) \to (E_{6,\sH})^C $ by
\begin{align*}
f_{6,\sH^C}(A)X=k^{-1}(A(k X)A^*),\;\; X \in \mathfrak{J}(3,\H^C),
\end{align*}
where as for the mapping $ k $, see \cite[in the beginning of Subsection 3.5 (p.211)]{iy7} in detail.

This mapping induces the required isomorphism.
\end{proof}

We prove the lemma used in the subsequent subsections.

\begin{lemma}\label{lemma 5.1.2}
The mapping $ f_{6,\sH^C}: SU^*(6,\C^C) \to (E_{6,\sH})^C $ of Theorem {\rm \ref{theprem 5.1.1}} satisfies the following
\begin{align*}
&(1)\,{}^t\!(f_{6,\sH^C}(A))^{-1}=f_{6,\sH^C}((A^*)^{-1}).
\quad (2)\,\tau f_{6,\sH^C}(A)\tau=f_{6,\sH^C}(\tau A).
\\
&(3)\,\gamma_{{}_{\scalebox{0.8}{\sC}}} f_{6,\sH^C}(A)\gamma_{{}_{\scalebox{0.8}{\sC}}}=f_{6,\sH^C}(JAJ).
\quad \;\;\,(4)\,\sigma f_{6,\sH^C}(A)\sigma=f_{6,\sH^C}(I_2AI_2),
\end{align*}
where $ I_2:=\diag(-1,-1,1,1,1,1) \in M(6,\R) $.
\end{lemma}
\begin{proof}
(1) Using the formula $ k X^*=(k X)^*, X \in M(3,\H^C) $, it follows from
\begin{align*}
({}^t\!(f_{6,\sH^C}(A))X, Y)&=(X,f_{6,\sH^C}(A)Y )=(X,k^{-1}(A(k Y)A^*))
\\
&=(X,(k^{-1}A)X(k^{-1}A^*))=((k^{-1}A)^*X(k^{-1}A^*)^*,Y)
\\
&=((k^{-1}A^*)X(k^{-1}A),Y)=(f_{6,\sH^C}(A^*)X,Y)
\end{align*}
that $ {}^t\!(f_{6,\sH^C}(A))=f_{6,\sH^C}(A^*) $. Hence we see $ {}^t\!(f_{6,\sH^C}(A))^{-1}=f_{6,\sH^C}((A^*)^{-1}) $.
\vspace{1mm}

(2) It follows from
\begin{align*}
\tau f_{6,\sH^C}(A)\tau X&=\tau k^{-1}(A(k(\tau X)A^*))=k^{-1}(\tau A(\tau k(\tau X))\tau A^*)
\\
&=k^{-1}(\tau A(kX)\tau A^*)=f_{6,\sH^C}(\tau A)X
\end{align*}
that $ \tau f_{6,\sH^C}(A)\tau=f_{6,\sH^C}(\tau A) $.
\vspace{1mm}

(3) Using the formula $ \gamma_{{}_{\scalebox{0.8}{\sC}}}=f_{6,\sH^C}(J), J=\diag(J_1,J_1,J_1) \in SU^*(6,\C^C),J_1:=
\begin{pmatrix}
0 & 1 \\
-1 & 0
\end{pmatrix} $, we have the following
\begin{align*}
\gamma_{{}_{\scalebox{0.8}{\sC}}}f_{6,\sH^C}(A)\gamma_{{}_{\scalebox{0.8}{\sC}}} &=f_{6,\sH^C}(J)f_{6,\sH^C}(A)f_{6,\sH^C}(J)=f_{6,\sH^C}(JAJ).
\end{align*}
\vspace{1mm}

(4) Using the formula $ \sigma=f_{6,\sH^C}(I_2), I_2=\diag(-1,-1,1,1,1,1) \in SU^*(6,\C^C) $, we have the following
\begin{align*}
\sigma f_{6,\sH^C}(A)\sigma=f_{6,\sH^C}(I_2)f_{6,\sH^C}(A)f_{6,\sH^C}(I_2)=f_{6,\sH^C}(I_2AI_2).
\end{align*}
\end{proof}

\subsection{The group $ E_{6,\sH} $}

The structure of the group $ E_{6,\sH} $ has been also determined by Ichiro Yokota (\cite[Proposition 3.11.3]{iy0}). In the present article, we will determine the structure of the group $ E_{6,\sH} $ by a different approach using the result of $ (E_{6,\sH})^C $.

As mentioned in the surface of this section, since the group $ (E_{6,\sH})^C $ has an involutive automorphism $ \tilde{\tau\lambda} $, we can consider the subgroup $ ((E_{6,\sH})^C)^{\tau\lambda} $ of $ (E_{6,\sH})^C $:
\begin{align*}
((E_{6,\sH})^C)^{\tau\lambda}:=\left\lbrace \alpha \in (E_{6,\sH})^C \relmiddle{|} \tilde{\tau\lambda}(\alpha)=\alpha \right\rbrace.
\end{align*}

Then we have the following theorem.

\begin{theorem}\label{theorem 5.2.1}
The group $ ((E_{6,\sH^C})^C)^{\tau\lambda} $ coincides with to the group $ E_{6,\sH} ${\rm :} $ ((E_{6,\sH^C})^C)^{\tau\lambda}=E_{6,\sH} $.
\end{theorem}
\begin{proof}
Let $ \alpha \in ((E_{6,\sH^C})^C)^{\tau\lambda} $. It follows from
$ \tau{}^t\!\alpha^{-1}\tau=\alpha $ that
\begin{align*}
\langle \alpha X, \alpha Y \rangle =(\tau\alpha X, \alpha Y )=({}^t\!\alpha^{-1}\tau X,Y)=(\tau X,\alpha^{-1}\alpha Y)=(\tau X,Y )=\langle X, Y \rangle.
\end{align*}
Hence we see $ \alpha \in E_{6,\sH} $. Conversely, let $ \beta \in E_{6,\sH} $. Then it follows from $ \langle \beta X, \beta Y \rangle=\langle X , Y \rangle $ that
\begin{align*}
(\tau X,Y)=\langle X , Y \rangle=\langle \beta X, \beta Y \rangle=(\tau\beta X, \beta Y)=({}^t\!\beta\beta\tau X, Y),
\end{align*}
that is, $ \tau={}^t\!\beta\beta\tau $. Hence we have $ \tau{}^t\!\beta^{-1}\tau=\beta $, so that $ \beta \in ((E_{6,\sH})^C)^{\tau\lambda} $.

With above, the proof of this theorem is completed.
\end{proof}

Let the mapping $ \phi:SU(6,\C^C) \to SU^*(6,\C^C) $ defined by $ \phi(B)=\iota B-\ov{\iota}J\ov{B}J,\iota:=(1/2)(1+i\i) $, then this mapping gives an isomorphism $ SU(6,\C^C) \cong SU^*(6,\C^C) $ (\cite[Lemma 3.5.10]{iy7}). Hence the composition mapping $ f_{6,\sH^C}\phi $ of $ \phi $ and $ f_{6,\sH^C} $ induces the isomorphism $ (E_{6,\sH})^C \cong SU(6,\C^C)/\Z_2,\Z_2=\{E,-E\} $:
\begin{align*}
SU(6,\C^C) \overset{\phi}{\longrightarrow} SU^*(6,\C^C) \overset{f_{6,{\scalebox{0.8}{\sH}^C}}}{\longrightarrow} (E_{6,\sH})^C.
\end{align*}
We denote the composition mapping $ f_{6,\sH^C}\phi $ by $ g_{6,\sH^C} $: $ g_{6,\sH^C}:=f_{6,\sH^C}\phi $.

Here, we prove the lemma needed in the proof of theorem below.

\begin{lemma}\label{lemma 5.2.2}
For $ B \in SU(6,\C^C) $, the mapping $ g_{6,\sH^C} $ satisfies $ \tau {}^t\!(g_{6,\sH^C}(B))^{-1}\tau=g_{6,\sH^C}(\tau B) $.
\end{lemma}
\begin{proof}
First we have $ \phi(B)^*=\phi(-J{}^t\!BJ) $. Indeed, it follows from
\begin{align*}
\phi(B)^*&=(\iota B-\ov{\iota}J\ov{B}J)^*=\ov{\iota}B^*-\iota (-J){}^t\!B(-J)=\ov{\iota}B^*-\iota J{}^t\!B J
\\
&=\iota(-J{}^t\!BJ)-\ov{\iota}J\ov{(-J{}^t\!BJ)}J=\phi(-J{}^t\!BJ)
\end{align*}
that $ \phi(B)^*=\phi(-J{}^t\!BJ) $. Hence we have $ (\phi(B)^*)^{-1}=\phi(-J{}^t\!B^{-1}J) $, that is, $ (\phi(B)^*)^{-1}=\phi(-J\ov{B}J) $.
Moreover, we have $ \tau\phi(B)=\phi(-J\ov{\tau B}J) $. Indeed, it follows that
\begin{align*}
\tau\phi(B)&=\tau(\iota B-\ov{\iota}J\ov{B}J)=\ov{\iota}\tau B-\iota J\ov{\tau B}J=\iota(-J\ov{\tau B}J)-\ov{\iota}J\ov{(-J\ov{\tau B}J)}J
\\
&=\phi(-J\ov{\tau B}J).
\end{align*}
Note that if $ B \in SU(6,\C^C) $, then we have $ -J\ov{B}J, -J\ov{\tau B}J \in SU(6,\C^C) $.
Hence it follows from Lemma \ref{lemma 5.1.2} (1), (2) that
\begin{align*}
\tau{}^t\!(g_{6,\sH^C}(B))^{-1}\tau&=\tau{}^t\!(f_{6,\sH^C}(\phi(B)))^{-1}\tau=f_{6,\sH^C}(\tau(\phi(B)^*)^{-1})=f_{6,\sH^C}\phi(\tau B)
\\
&=g_{6,\sH^C}(\tau B).
\end{align*}

With above, this lemma is proved.
\end{proof}

Now, we determine the structure of the group $ E_{6,\sH} $.

\begin{theorem}\label{theorem 5.2.3}
The group $ E_{6,\sH} $ is isomorphic to the group $ SU(6)/\Z_2 ${\rm :} $ E_{6,\sH} \cong SU(6)/\Z_2 $.
\end{theorem}
\begin{proof}
Let the group $ E_{6,\sH} $ as the group $ ((E_{6,\sH})^C)^{\tau\lambda} $. The we define a mapping $ g_{6,\tau\lambda}:SU(6) \to ((E_{6,\sH})^C)^{\tau\lambda} $ by
\begin{align*}
g_{6,\tau\lambda}(A)X=k^{-1}(\phi(A)(k X)\phi(A)^*), \;\;X \in \mathfrak{J}(3,\H^C).
\end{align*}
Note that this mapping is the restriction of the mapping $ g_{6,\sH^C} $.  First, we will prove that $ g_{6,\tau\lambda} $ is well-defined and a homomorphism. Since $ g_{6,\tau\lambda} $ is the restriction of the mapping $ g_{6,\sH^C} $, it is clear that $ g_{6,\tau\lambda}(A) \in (E_{6,\sH})^C $ and $ g_{6,\tau\lambda} $ is a homomorphism. Moreover, from Lemma \ref{lemma 5.2.2}, we have $ \tau {}^t\!(g_{6,\tau\lambda}(A))^{-1}\tau=g_{6,\tau\lambda}(A) $, so that $ g_{6,\tau\lambda}(A) \in ((E_{6,\sH})^C)^{\tau\lambda} $.

Next, we will prove that $ g_{6,\tau\lambda} $ is surjective. Let $ \alpha \in E_{6,\sH}=((E_{6,\sH})^C)^{\tau\lambda} \subset (E_{6,\sH})^C $. Then there exists $ B \in SU(6,\C^C) $ such that $ \alpha=g_{6,\sH^C}(B) $ (as mentioned in the beginning of p.26). Moreover $ \alpha $ satisfies the condition $ \tau{}^t\!\alpha^{-1}\tau=\alpha $, that is, $ \tau {}^t\!(g_{6,\sH^C}(B))^{-1}\tau=g_{6,\sH^C}(B) $, so that from Lemma \ref{lemma 5.2.2} we have the following
\begin{align*}
\tau B=B \quad \text{or} \quad \tau B=-B.
\end{align*}
In the former case, we have $ B \in SU(6) $. In the latter case, $ B $ is of the form $ iB', B' \in M(6,\C) $. \vspace{1mm}
Since $ B^*B=E $, we have $ {B'}^*B'=-E $. Here, set $ B':=\left(  \begin{array}{c@{\,\,}c@{\,\,}c@{\,\,}c@{\,\,}c@{\,\,}c}
b_{11} & b_{12} & b_{13} & b_{14} & b_{15} & b_{16}
\\
b_{21} & b_{22} &\ldots &&& b_{26}
\\
b_{31} & \ldots & &&& b_{36}
\\
b_{41} & \ldots && && b_{46}
\\
b_{51} & \ldots &&& & b_{56}
\\
b_{61} & \ldots &\ldots&\ldots&\ldots& b_{66}
\end{array}\right),b_{ij} \in \C   $, then we have
\begin{align*}
-E={B'}^*B'=\left( \begin{array}{c@{\,\,}c@{\,\,}c@{\,\,}c@{\,\,}c@{\,\,}c}
c_{11} &&&&&
\\
& c_{22} &&&*&
\\
&& c_{33} &&&
\\
&&& c_{44} &&
\\
&*&&& c_{55} &
\\
&&&&& c_{66}
\end{array}
\right), \;\text{where}\;\;
\begin{array}{l}
c_{ii}=\ov{b_{1i}}b_{1i}+\ov{b_{2i}}b_{2i}+\ov{b_{3i}}b_{3i}
\vspace{1mm}\\
\quad\;\;\,+\ov{b_{4i}}b_{4i}+\ov{b_{5i}}b_{5i}+\ov{b_{6i}}b_{6i}\geq 0 ,
\vspace{1mm}\\
\quad\;\;\, i=1,2,\ldots,6.
\end{array}
\end{align*}
By comparing the diagonal entries of both sides, the contradiction is followed, so that this case is impossible. With above, the proof of surjective is completed.

Finally, we will determine $ \Ker\,g_{6,\tau\lambda} $. However, since $ g_{6,\tau\lambda} $ the restriction of the mapping $ g_{6,\sH^C} $, it is easy to obtain $ \Ker\,g_{6,\tau\lambda}=\Ker\,g_{6,\sH^C}=\{E,-E\} \cong \Z_2 $.

Therefore we have the required isomorphism
\begin{align*}
E_{6,\sH} \cong SU(6)/\Z_2.
\end{align*}

\end{proof}

\subsection{The group $ E_{6,\sH'} $}

We consider a subalgebra $ (\H^C)_{\tau\gamma_{{}_{\scalebox{0.8}{\sC}}}} $ of $ \H^C $:
\begin{align*}
(\H^C)_{\tau\gamma_{{}_{\scalebox{0.8}{\sC}}}}:&=\left\lbrace p \in \H^C \relmiddle{|} \tau\gamma_{{}_{\scalebox{0.8}{\sC}}}p=p \right\rbrace
\\
&=\left\lbrace p=(x+yi\i)+(z+i\i w)\j \relmiddle{|} x,y,z,w \in \R \right\rbrace
\\
&=(\C^C)_{\tau\gamma_{{}_{\scalebox{0.8}{\sC}}}}\oplus (\C^C)_{\tau\gamma_{{}_{\scalebox{0.8}{\sC}}}}\j.
\end{align*}

Let $ \H' $ be the algebra of split quaternion numbers: $ \H':=\C' \oplus \C'\j, \j^2=-1 $. Then, using the mapping $ l $ defined in the beginning of Subsection 4.3, the correspondence
\begin{align*}
(\H^C)_{\tau\gamma_{{}_{\scalebox{0.8}{\sC}}}} \ni (x+yi\i)+(z+wi\i)\j \underset{l}{\mapsto} (x+y\i')+(z+w\i')\j \in \H'
\end{align*}
gives an isomorphism $ (\H^C)_{\tau\gamma_{{}_{\scalebox{0.8}{\sC}}}} \simeq \H' $ as algebras.

We define a subalgebra $ (\mathfrak{J}(3,\H^C))_{\tau\gamma_{{}_{\scalebox{0.8}{\sC}}}} $ of $ \mathfrak{J}(3,\H^C) $ by
\begin{align*}
(\mathfrak{J}(3,\H^C))_{\tau\gamma_{{}_{\scalebox{0.8}{\sC}}}}:&=\left\lbrace X \in \mathfrak{J}(3,\H^C) \relmiddle{|} {\tau\gamma_{{}_{\scalebox{0.8}{\sC}}}} X=X \right\rbrace
\\
&=\left\lbrace
X=\begin{pmatrix}
\xi & x_3 & \ov{x}_2 \\
\ov{x}_3 & \xi_2 & x_1 \\
x_2 & \ov{x}_1 & \xi_3
\end{pmatrix} \relmiddle{|} \xi_i \in \R, \x_i \in (\H^C)_{\tau\gamma_{{}_{\scalebox{0.8}{\sC}}}} \right\rbrace.
\end{align*}
Let $ \mathfrak{J}(3,\H') $ be the split Jordan algebra. Then the correspondence
\begin{align*}
(\mathfrak{J}(3,\H^C))_{\tau\gamma_{{}_{\scalebox{0.8}{\sC}}}}\ni
\begin{pmatrix}
\xi & x_3 & \ov{x}_2 \\
\ov{x}_3 & \xi_2 & x_1 \\
x_2 & \ov{x}_1 & \xi_3
\end{pmatrix} \underset{f}{\mapsto}
\begin{pmatrix}
\xi & l(x_3) & \ov{l(x_2)} \\
\ov{l(x_3)} & \xi_2 & l(x_1) \\
l(x_2) & \ov{l(x_1)} & \xi_3
\end{pmatrix} \in \mathfrak{J}(3,\H')
\end{align*}
gives an isomorphism $ (\mathfrak{J}(3,\C^C))_{\tau\gamma_{{}_{\scalebox{0.8}{\sC}}}} \simeq  \mathfrak{J}(3,\C') $ as algebras with the properties $ \det(fX)=\det\,X $.

We define a group $ E_{6(6),(\sH^C)_{\tau\gamma_{{}_{\scalebox{0.8}{\sC}}}}} $ by
\begin{align*}
E_{6(6),(\sH^C)_{\tau\gamma_{{}_{\scalebox{0.8}{\sC}}}}}:=\left\lbrace \alpha \in \Iso_{\sR}((\mathfrak{J}(3,\H^C))_{\tau\gamma_{{}_{\scalebox{0.8}{\sC}}}}) \relmiddle{|}\det(\alpha X)=\det\, X\right\rbrace.
\end{align*}

Then we have the following proposition.

\begin{proposition}\label{proposition 5.3.1}
The group $ E_{6(6),(\sH^C)_{\tau\gamma_{{}_{\scalebox{0.8}{\sC}}}} } $ is isomorphic to the group $ E_{6,\sH'} ${\rm :} $ E_{6(6),(\sH^C)_{\tau\gamma_{{}_{\scalebox{0.8}{\sC}}}} } \cong E_{6(6),\sH'} $.
\end{proposition}
\begin{proof}
We can prove this proposition by replacing $ \C^C, \C' $ with $ \H^C, \H' $ in the proof of Proposition \ref{proposition 4.3.1}, respectively.
\end{proof}

Since the group $ (E_{6,\sH})^C $ has an involutive automorphism $ \tilde{\tau\gamma_{{}_{\scalebox{0.8}{\sC}}}} $, we can define a subgroup $ ((E_{6,\sH})^C)^{\tau\gamma_{{}_{\scalebox{0.8}{\sC}}}} $ of $ (E_{6,\sH})^C $ by
\begin{align*}
((E_{6,\sH})^C)^{\tau\gamma_{{}_{\scalebox{0.8}{\sC}}}}:=\left\lbrace\alpha \in (E_{6,\sH})^C \relmiddle{|}\tilde{\tau\gamma_{{}_{\scalebox{0.8}{\sC}}}}(\alpha)=\alpha \right\rbrace.
\end{align*}

Then we prove the following theorem.

\begin{theorem}\label{theorem 5.3.2}
The group $ ((E_{6,\sH})^C)^{\tau\gamma_{{}_{\scalebox{0.8}{\sC}}}} $ coincides with the group $ E_{6(6),(\sH^C)_{\tau\gamma_{\scalebox{0.8}{\sC}}} } ${\rm :} $ ((E_{6,\sH})^C)^{\tau\gamma_{{}_{\scalebox{0.8}{\sC}}}}=E_{6(6),(\sH^C)_{\tau\gamma_{\scalebox{0.8}{\sC}}} } $.

In particular, we have the isomorphism $ ((E_{6,\sH})^C)^{\tau\gamma_{{}_{\scalebox{0.8}{\sC}}}} \cong E_{6(6),\sH'} $.
\end{theorem}
\begin{proof}
This proposition can be proved in almost the same way
by replacing $ \C^C, \C' $ with $ \H^C, \H' $ in the proof of Theorem \ref{theorem 4.3.2}, respectively, however we will rewrite the proof as detailed as possible.

Let $ \alpha \in ((E_{6,\sH})^C)^{\tau\gamma_{{}_{\scalebox{0.8}{\sC}}}} $. Since $ (\tau\gamma_{{}_{\scalebox{0.8}{\sC}}})\alpha=\alpha(\tau\gamma_{{}_{\scalebox{0.8}{\sC}}}) $, $ \alpha $ induces an $ \R $-linear isomorphism of $ (\mathfrak{J}(3,\H^C))_{\tau\gamma_{{}_{\scalebox{0.8}{\sC}}}} $. Moreover, since it is clear that $ \det(\alpha X)=\det\,X, X \in (\mathfrak{J}(3,\H^C))_{\tau\gamma_{{}_{\scalebox{0.8}{\sC}}}} $, we have $ \alpha \in E_{6(6),(\sH^C)_{\tau\gamma_{{}_{\scalebox{0.8}{\sC}}}} } $. Conversely, let $ \beta \in  E_{6(6),(\sH^C)_{\tau\gamma_{{}_{\scalebox{0.8}{\sC}}}} } $.
Since $ \mathfrak{J}(3,\H^C) $ is decomposed as $ (\mathfrak{J}(3,\H^C))_{\tau\gamma_{{}_{\scalebox{0.8}{\sC}}}} \oplus i(\mathfrak{J}(3,\H^C))_{\tau\gamma_{{}_{\scalebox{0.8}{\sC}}}} $: $ \mathfrak{J}(3,\H^C)=(\mathfrak{J}(3,\H^C))_{\tau\gamma_{{}_{\scalebox{0.8}{\sC}}}} \oplus i(\mathfrak{J}(3,\H^C))_{\tau\gamma_{{}_{\scalebox{0.8}{\sC}}}} $, that is, $ \mathfrak{J}(3,\H^C) $ is the complexification of $ (\mathfrak{J}(3,\H^C))_{\tau\gamma_{{}_{\scalebox{0.8}{\sC}}}} $, we can define an action to $ \mathfrak{J}(3,\H^C) $ of the group $  E_{6(6),(\sH^C)_{\tau\gamma_{{}_{\scalebox{0.8}{\sC}}}} } $ by
\begin{align*}
\beta X=\beta (X_1+iX_2)=\beta X_1+i\beta X_2,\;\; X:=X_1+iX_2 \in \mathfrak{J}(3,\H^C), X_i \in (\mathfrak{J}(3,\H^C))_{\tau\gamma_{{}_{\scalebox{0.8}{\sC}}}}.
\end{align*}
Then it follows from
\begin{align*}
\beta X \times \beta Y&=\beta (X_1+iX_2) \times \beta (Y_1+iY_2)=(\beta X_1+i\beta X_2) \times (\beta Y_1+i\beta Y_2)
\\
&=(\beta X_1 \times \beta Y_1 - \beta X_2 \times \beta Y_2)+i(\beta X_1 \times \beta Y_2+\beta X_2 \times \beta Y_1)
\\
&=({}^t\!\beta^{-1}(X_1 \times Y_1)-{}^t\!\beta^{-1}(X_2 \times Y_2))+i({}^t\!\beta^{-1}(X_1 \times Y_2)+{}^t\!\beta^{-1}(X_2 \times Y_1))
\\
&={}^t\!\beta^{-1}((X_1 \times Y_1)-(X_2 \times Y_2))+i((X_1 \times Y_2)+(X_2 \times Y_1))
\\
&={}^t\!\beta^{-1}((X_1+iX_2)\times (Y_1+iY_2))
\\
&={}^t\!\beta^{-1}(X \times Y)
\end{align*}
that $ \beta \in (E_{6,\sH})^C $. Moreover, it is easy to verify that $ (\tau\gamma_{{}_{\scalebox{0.8}{\sC}}})\beta=\beta(\gamma_{{}_{\scalebox{0.8}{\sC}}}\tau) $. Indeed, for $ X:=X_1+iX_2 \in (\mathfrak{J}(3,\H^C))_{\tau\gamma_{{}_{\scalebox{0.8}{\sC}}}} \oplus i(\mathfrak{J}(3,\H^C))_{\tau\gamma_{{}_{\scalebox{0.8}{\sC}}}}=\mathfrak{J}(3,\H^C) $, it follows that
\begin{align*}
(\tau\gamma_{{}_{\scalebox{0.8}{\sC}}})\beta(\gamma_{{}_{\scalebox{0.8}{\sC}}}\tau)X&=(\tau\gamma_{{}_{\scalebox{0.8}{\sC}}})\beta(\gamma_{{}_{\scalebox{0.8}{\sC}}}\tau)(X_1+iX_2)=(\tau\gamma_{{}_{\scalebox{0.8}{\sC}}})\beta(X_1-iX_2)
\\
&=(\tau\gamma_{{}_{\scalebox{0.8}{\sC}}})(\beta X_1-i\beta X_2)=\beta X_1+i\beta X_2=\beta(X_1+iX_2)
\\
&=\beta X,
\end{align*}
that is, $ (\tau\gamma_{{}_{\scalebox{0.8}{\sC}}})\beta(\gamma_{{}_{\scalebox{0.8}{\sC}}}\tau)=\beta $.
With above, we have $ \beta \in ((E_{6,\sH})^C)^{\tau\gamma_{{}_{\scalebox{0.8}{\sC}}}} $.

Thus we have the required result $ ((E_{6,\sC})^C)^{\tau\gamma_{{}_{\scalebox{0.8}{\sC}}}}=E_{6(6),(\sC^C)_{\tau\gamma_{{}_{\scalebox{0.8}{\sC}}}} } $.

Therefore, together with Proposition \ref{proposition 5.3.1}, we have the isomorphism
\begin{align*}
((E_{6,\sH})^C)^{\tau\gamma_{{}_{\scalebox{0.8}{\sC}}}} \cong E_{6(6),\sH'}.
\end{align*}
\end{proof}

Let the mapping $ \psi:SL(6,C) \to SU^*(6,\C^C) $ defined by $ \psi(B)=\iota B-\ov{\iota}JBJ,\iota:=(1/2)(1+i\i) $, then this mapping gives an isomorphism $ SL(6,C) \cong SU^*(6,\C^C) $ (\cite[in the proof of Theorem 3.5.9]{iy7}). Hence the composition mapping $ f_{6,\sH^C}\psi $ of $ \psi $ and $ f_{6,\sH^C} $ induces the isomorphism $ (E_{6,\sH})^C \cong SL(6,C)/\Z_2,\Z_2=\{E,-E\} $:
\begin{align*}
SL(6,C) \overset{\psi}{\longrightarrow} SU^*(6,\C^C) \overset{f_{6,{\scalebox{0.8}{\sH}^C}}}{\longrightarrow} (E_{6,\sH})^C.
\end{align*}
We denote the composition mapping $ f_{6,\sH^C}\psi $ by $ h_{6,\sH^C} $: $ h_{6,\sH^C}:=f_{6,\sH^C}\psi $.

Here, we prove the lemma needed in the proof of theorem below.

\begin{lemma}\label{lemma 5.3.3}
For $ B \in SL(6,C) $, the mapping $ h_{6,\sH^C} $ satisfies $ \tau \gamma_{{}_{\scalebox{0.8}{\sC}}}h_{6,\sH^C}(B)\gamma_{{}_{\scalebox{0.8}{\sC}}}\tau=h_{6,\sH^C}(\tau B) $.
\end{lemma}
\begin{proof}
First we have $ J\psi(B)J=\psi(JBJ) $. Indeed, it follows that
\begin{align*}
J\psi(B)J&=J(\iota B-\ov{\iota}JBJ)J=\iota JBJ-\ov{\iota}B=\iota JBJ-\ov{\iota}J(JBJ)J
\\
&=\psi(JBJ).
\end{align*}
Moreover, since $ \tau\psi(B)=\psi(-J(\tau B)J) $, we have $ \tau(J\psi(B)J)=\psi(-\tau B) $. Indeed, it follows that
\begin{align*}
\tau(J\psi(B)J)=\tau\psi(JBJ)=\phi(-\tau B).
\end{align*}
Hence, note that $ f_{6,\sH^C}(-E)=1 $, it follows from Lemma \ref{lemma 5.1.2} (2), (3) and the result above that
\begin{align*}
\tau \gamma_{{}_{\scalebox{0.8}{\sC}}}h_{6,\sH^C}(B)\gamma_{{}_{\scalebox{0.8}{\sC}}}\tau&=\tau \gamma_{{}_{\scalebox{0.8}{\sC}}}f_{6,\sH^C}(\psi(B))\gamma_{{}_{\scalebox{0.8}{\sC}}}\tau=\tau f_{6,\sH^C}(J\psi(B)J)\tau
\\
&=f_{6,\sH^C}(\tau(J\psi(B)J))=f_{6,\sH^C}(\psi(-\tau B))
\\
&=f_{6,\sH^C}((-E)\phi(\tau B))=f_{6,\sH^C}(-E)f_{6,\sH^C}(\phi(\tau B))
\\
&=f_{6,\sH^C}(\psi(\tau B))=f_{6,\sH^C}(\psi(\tau\gamma_{{}_{\scalebox{0.8}{\sC}}} B))
\\
&=h_{6,\sH^C}(\tau B).
\end{align*}

With above, this lemma is proved.
\end{proof}

Let $ I_i:=\diag(-i,i,i,i,i,i) \in SL(6,C) $. We define a $ C $-linear isomorphism $ \nu $ of $ \mathfrak{J}(3,\H^C) $ by $ \nu:=h_{6,\sH^C}(I_i) $:
\begin{align*}
\nu X=h_{6,\sH^C}(I_i)X=k^{-1}(\psi(I_i)(k X)\psi(I_i)^*),\;\; X \in \mathfrak{J}(3,\H^C).
\end{align*}
Then we have $ \nu \in (E_{6,\sH})^C $ and $ \nu^2=1 $.

We consider a discrete group $ \Z_2:=\{1,\nu \} $. This group acts on $ SL(6,\R) $ by
\begin{align*}
1 A=A, \quad \nu A=I_i A {I_i}^{-1}\,(\in SL(6,\R)),
\end{align*}
and then $ SL(6,\R) \rtimes \Z_2 $ be the semi-direct product of $ SL(6,\R) $ and $ \Z_2=\{1,\nu\} $ with the multiplication
\begin{align*}
&(A_1,1)(A_2,1)=(A_1A_2,1), \quad (A_1,1)(A_2,\nu)=(A_1A_2,\nu),
\\
&(A_1,\nu)(A_2,1)=(A_1(I_iA_2{I_i}^{-1}),\nu), \quad (A_1,\nu)(A_2,\nu)=(A_1(I_iA_2{I_i}^{-1}),1).
\end{align*}

Now, we determine the structure of the group $ E_{6(6),\sH'} $.

\begin{theorem}\label{theorem 5.3.5}
The group $ E_{6(6),\sH'} $ is isomorphic to the group $ SL(6,\R)/\Z_2 \rtimes \Z_2,\Z_2=\{E,-E \}, \Z_2=\{1,\nu \} ${\rm :} $E_{6(6),\sH'} \cong SL(6,\R)/\Z_2 \rtimes \Z_2 $.
\end{theorem}
\begin{proof}
Let the group $ E_{6(6),\sH'} $ as the group $ ((E_{6,\sH})^C)^{\tau\gamma_{{}_{\scalebox{0.8}{\sC}}}} $ (Theorem \ref{theorem 5.3.2}).
Then we define a mapping $ h_{6,\tau\gamma_{{}_{\scalebox{0.8}{\sC}}}}:SL(6,\R) \rtimes \{1,\nu \} \to ((E_{6,\sH})^C)^{\tau\gamma_{{}_{\scalebox{0.8}{\sC}}}} $ by
\begin{align*}
h_{6,\tau\gamma_{{}_{\scalebox{0.8}{\sC}}}}(A,1)X&=k^{-1}(\psi(A)(k X)\psi(A)^*),
\\
h_{6,\tau\gamma_{{}_{\scalebox{0.8}{\sC}}}}(A,\nu)X&=k^{-1}(\psi(AI_i)(k X)\psi(AI_i)^*)\;\;X \in \mathfrak{J}(3,\H^C).
\end{align*}
First, we will prove that $ h_{6,\tau\gamma_{{}_{\scalebox{0.8}{\sC}}}} $ is well-defined. Since the mapping $ h_{6,\tau\gamma_{{}_{\scalebox{0.8}{\sC}}}} $ of the former case is the restriction of the mapping $ h_{6,\sH^C} $ and together with Lemma \ref{lemma 5.3.3}, it is clear that $ h_{6,\tau\gamma_{{}_{\scalebox{0.8}{\sC}}}}(A,1) \in ((E_{6,\sH})^C)^{\tau\gamma_{{}_{\scalebox{0.8}{\sC}}}} $. Since the mapping $ h_{6,\tau\gamma_{{}_{\scalebox{0.8}{\sC}}}} $ is a homomorphism, which is shown later, we see $ h_{6,\tau\gamma_{{}_{\scalebox{0.8}{\sC}}}}(A,\nu) =h_{6,\tau\gamma_{{}_{\scalebox{0.8}{\sC}}}}(A,1)\allowbreak h_{6,\tau\gamma_{{}_{\scalebox{0.8}{\sC}}}}(E,\nu)$. Hence we have to show $ h_{6,\tau\gamma_{{}_{\scalebox{0.8}{\sC}}}}(E,\nu) \in ((E_{6,\sH})^C)^{\tau\gamma_{{}_{\scalebox{0.8}{\sC}}}} $. It follows from Lemma \ref{lemma 5.3.3} that
\begin{align*}
\tau\gamma_{{}_{\scalebox{0.8}{\sC}}}h_{6,\tau\gamma_{{}_{\scalebox{0.8}{\sC}}}}(E,\nu)\gamma_{{}_{\scalebox{0.8}{\sC}}}\tau
&=\tau\gamma_{{}_{\scalebox{0.8}{\sC}}}h_{6,\sH^C}(I_i)\gamma_{{}_{\scalebox{0.8}{\sC}}}\tau=h_{6,\sH^C}(\tau I_i)=h_{6,\sH^C}(-I_i)
\\
&=h_{6,\sH^C}(-E)h_{6,\sH^C}(I_i)=h_{6,\sH^C}(I_i)
\\
&=h_{6,\tau\gamma_{{}_{\scalebox{0.8}{\sC}}}}(E,\nu).
\end{align*}
Hence we have $ h_{6,\tau\gamma_{{}_{\scalebox{0.8}{\sC}}}}(E,\nu) \in ((E_{6,\sH})^C)^{\tau\gamma_{{}_{\scalebox{0.8}{\sC}}}} $, so that $ h_{6,\tau\gamma_{{}_{\scalebox{0.8}{\sC}}}}(A,\nu) \in ((E_{6,\sH})^C)^{\tau\gamma_{{}_{\scalebox{0.8}{\sC}}}} $. With above, the proof of well-defined is proved. Subsequently, we will prove that $ h_{6,\tau\gamma_{{}_{\scalebox{0.8}{\sC}}}} $ is a homomorphism. That is shown as follows.
\begin{align*}
h_{6,\tau\gamma_{{}_{\scalebox{0.8}{\sC}}}}(A_1,1)h_{6,\tau\gamma_{{}_{\scalebox{0.8}{\sC}}}}(A_2,1)X
&=h_{6,\tau\gamma_{{}_{\scalebox{0.8}{\sC}}}}(A_1,1)(k^{-1}(\psi(A_2)(k X)\psi(A_2)^*))
\\
&=k^{-1}(\psi(A_1)k(k^{-1}(\psi(A_2)(k X)\psi(A_2)^*))\psi(A_1)^*)
\\
&=k^{-1}(\psi(A_1)\psi(A_2)(k X)\psi(A_2)^*\psi(A_1)^*)
\\
&=k^{-1}(\psi(A_1A_2)(k X)\psi(A_1A_2)^*)
\\
&=h_{6,\tau\gamma_{{}_{\scalebox{0.8}{\sC}}}}(A_1A_2,1)X
\\
&=h_{6,\tau\gamma_{{}_{\scalebox{0.8}{\sC}}}}((A_1,1)(A_2,1))X,
\\[1mm]
h_{6,\tau\gamma_{{}_{\scalebox{0.8}{\sC}}}}(A_1,1)h_{6,\tau\gamma_{{}_{\scalebox{0.8}{\sC}}}}(A_2,\nu)X
&=h_{6,\tau\gamma_{{}_{\scalebox{0.8}{\sC}}}}(A_1,1)(k^{-1}(\psi(A_2I_i)(k X)\psi(A_2I_i)^*))
\\
&=k^{-1}(\psi(A_1)k(k^{-1}(\psi(A_2)(k X)\psi(A_2)^*))\psi(A_1)^*)
\\
&=k^{-1}(\psi(A_1)\psi(A_2I_i)(k X)\psi(A_2I_i)^*\psi(A_1)^*)
\\
&=k^{-1}(\psi(A_1A_2I_i)(k X)\psi(A_1A_2I_i)^*)
\\
&=h_{6,\tau\gamma_{{}_{\scalebox{0.8}{\sC}}}}(A_1A_2,\nu)X
\\
&=h_{6,\tau\gamma_{{}_{\scalebox{0.8}{\sC}}}}((A_1,1)(A_2,\nu))X,
\\[1mm]
h_{6,\tau\gamma_{{}_{\scalebox{0.8}{\sC}}}}(A_1,\nu)h_{6,\tau\gamma_{{}_{\scalebox{0.8}{\sC}}}}(A_2,1)X
&=h_{6,\tau\gamma_{{}_{\scalebox{0.8}{\sC}}}}(A_1,\nu)(k^{-1}(\psi(A_2)(k X)\psi(A_2)^*))
\\
&=k^{-1}(\psi(A_1I_i)k(k^{-1}(\psi(A_2)(k X)\psi(A_2)^*))\psi(A_1I_i)^*)
\\
&=k^{-1}(\psi(A_1I_i)\psi(A_2)(k X)\psi(A_2)^*\psi(A_1I_i)^*)
\\
&=k^{-1}(\psi(A_1I_iA_2)(k X)\psi(A_1I_iA_2)^*)
\\
&=k^{-1}(\psi(A_1I_iA_2{I_i}^{-1}I_i)(k X)\psi(A_1I_iA_2{I_i}^{-1}I_i)^*)
\\
&=k^{-1}(\psi(A_1(I_iA_2{I_i}^{-1})I_i)(k X)\psi(A_1(I_iA_2{I_i}^{-1})I_i)^*)X
\\
&=h_{6,\tau\gamma_{{}_{\scalebox{0.8}{\sC}}}}(A_1(I_iA_2{I_i}^{-1}),\nu)
\\
&=h_{6,\tau\gamma_{{}_{\scalebox{0.8}{\sC}}}}((A_1,\nu)(A_2,1))X,
\\[1mm]
h_{6,\tau\gamma_{{}_{\scalebox{0.8}{\sC}}}}(A_1,\nu)h_{6,\tau\gamma_{{}_{\scalebox{0.8}{\sC}}}}(A_2,\nu)X
&=h_{6,\tau\gamma_{{}_{\scalebox{0.8}{\sC}}}}(A_1,\nu)(k^{-1}(\psi(A_2I_i)(k X)\psi(A_2I_i)^*))
\\
&=k^{-1}(\psi(A_1I_i)k(k^{-1}(\psi(A_2I_i)(k X)\psi(A_2I_i)^*))\psi(A_1I_i)^*)
\\
&=k^{-1}(\psi(A_1I_i)\psi(A_2I_i)(k X)\psi(A_2I_i)^*\psi(A_1I_i)^*)
\\
&=k^{-1}(\psi(A_1I_iA_2I_i)(k X)\psi(A_1I_iA_2I_i)^*)
\\
&=k^{-1}(\psi(A_1(I_iA_2{I_i}^{-1})(-E))(k X)\psi(A_1(I_iA_2{I_i}^{-1})(-E))^*)
\\
&=k^{-1}(\psi(A_1(I_iA_2{I_i}^{-1}))(k X)\psi(A_1(I_iA_2{I_i}^{-1}))^*)X
\\
&=h_{6,\tau\gamma_{{}_{\scalebox{0.8}{\sC}}}}(A_1(I_iA_2{I_i}^{-1}),1)
\\
&=h_{6,\tau\gamma_{{}_{\scalebox{0.8}{\sC}}}}((A_1,\nu)(A_2,\nu))X.
\end{align*}

Next, we will prove that $ h_{6,\tau\gamma_{{}_{\scalebox{0.8}{\sC}}}} $  is surjective. Let $ \alpha \in E_{6,\sH'}=((E_{6,\sH})^C)^{\tau\gamma_{{}_{\scalebox{0.8}{\sC}}}} \subset (E_{6,\sH})^C $. Then there exists $ B \in SL(6,C) $ such that $ \alpha=h_{6,\sH^C}(B) $ (as mentioned in the beginning of this page). Moreover $ \alpha $ satisfies the condition $ \tau\gamma_{{}_{\scalebox{0.8}{\sC}}}\alpha\gamma_{{}_{\scalebox{0.8}{\sC}}}\tau=\alpha $, that is, $ \tau\gamma_{{}_{\scalebox{0.8}{\sC}}}h_{6,\sH^C}(B)\gamma_{{}_{\scalebox{0.8}{\sC}}}\tau=h_{6,\sH^C}(\tau B) $, so that from Lemma \ref{lemma 5.3.3} we have the following
\begin{align*}
\tau B=B \quad \text{or} \quad \tau B=-B.
\end{align*}
In the former case, we have $ B \in SL(6,\R)$. Hence there exists $ (A,1) \in SL(6,\R) \rtimes \{1,\nu \} $ such that $ \alpha=h_{6,\sH^C}(A)=h_{6,\tau\gamma_{{}_{\scalebox{0.8}{\sC}}}}(A,1) $.
In the latter case, $ B $ is of the form $ iB', B' \in M(6,\R) $. Since $ \det\,B=1 $, we have $ \det\,B'=-1 $, and moreover $ B $ can be modified the form $ I_i(I_1B'),I_1:=\diag(-1,1,1,1,1,1) $: $ B=I_i(I_1B') $. Then we have $ I_1B' \in SL(6,\R) $. Here, set $ A:=I_1B' $, then there exists $ (A,\nu) \in SL(6,\R) \rtimes \{1,\nu \} $ such that $ \alpha=h_{6,\sH^C}(I_iA)=h_{6,\tau\gamma_{{}_{\scalebox{0.8}{\sC}}}}(A,\nu) $. With above, the proof of surjective is completed.

Finally, we will determine $ \Ker\,h_{6,\tau\gamma_{{}_{\scalebox{0.8}{\sC}}}} $. It follows from the definition of kernel that
\begin{align*}
\Ker\,h_{6,\tau\gamma_{{}_{\scalebox{0.8}{\sC}}}}&=\left\lbrace (A,1) \in  SL(6,\R) \rtimes \Z_2 \relmiddle{|} h_{6,\tau\gamma_{{}_{\scalebox{0.8}{\sC}}}}(A,1)=1 \right\rbrace
\\
&\cup \left\lbrace (A,\nu) \in  SL(6,\R) \rtimes \Z_2 \relmiddle{|} h_{6,\tau\gamma_{{}_{\scalebox{0.8}{\sC}}}}(A,\nu)=1 \right\rbrace.
\end{align*}
In the former case, since $ h_{6,\tau\gamma_{{}_{\scalebox{0.8}{\sC}}}}(A,1)=h_{6,\sH^C}(A) $, together with the result of $ \Ker\,h_{6,\sH^C} $, we have
\begin{align*}
\left\lbrace (A,1) \in  SL(6,\R) \rtimes \Z_2 \relmiddle{|} h_{6,\tau\gamma_{{}_{\scalebox{0.8}{\sC}}}}(A,1)=1 \right\rbrace&=\left\lbrace (E,1), (-E,1) \right\rbrace.
\end{align*}
In the latter case, let $ (A,\nu) \in \Ker\,h_{6,\tau\gamma_{{}_{\scalebox{0.8}{\sC}}}} $. Then $ (A,\nu) $ satisfies $ h_{6,\tau\gamma_{{}_{\scalebox{0.8}{\sC}}}}(A,\nu)=1 $, that is, $ h_{6,\tau\gamma_{{}_{\scalebox{0.8}{\sC}}}}(A,1)=h_{6,\tau\gamma_{{}_{\scalebox{0.8}{\sC}}}}(E,\nu) $. However there exists no $ (A.\nu) \in SL(6,\R) \rtimes \Z_2 $. Hence we see
\begin{align*}
\left\lbrace (A,\nu) \in  SL(6,\R) \rtimes \Z_2 \relmiddle{|} h_{6,\tau\gamma_{{}_{\scalebox{0.8}{\sC}}}}(A,\nu)=1 \right\rbrace =\emptyset.
\end{align*}

Therefore, from Theorem \ref{theorem 5.3.2}, we have the required isomorphism
\begin{align*}
E_{6(6),\sH'} \cong SL(6,\R)/\Z_2 \rtimes \Z_2.
\end{align*}
\end{proof}

\subsection{The group $ E_{6(-14),\sH} $}

As mentioned in the surface of this section, since the group $ (E_{6,\sH})^C $ has an involutive automorphism $ \tilde{\tau\lambda\sigma} $, we can consider the subgroup $ ((E_{6,\sH})^C)^{\tau\lambda\sigma} $ of $ (E_{6,\sH})^C $:
\begin{align*}
((E_{6,\sH})^C)^{\tau\lambda\sigma}:=\left\lbrace \alpha \in (E_{6,\sH})^C \relmiddle{|} \tilde{\tau\lambda\sigma}(\alpha)=\alpha \right\rbrace.
\end{align*}

Then we have the following theorem.

\begin{theorem}\label{theorem 5.4.1}
The group $ ((E_{6,\sH})^C)^{\tau\lambda\sigma} $ coincides with to the group $ E_{6(-14),\sH} ${\rm :} $ ((E_{6,\sH})^C)^{\tau\lambda\sigma}=E_{6(-14),\sH} $.
\end{theorem}
\begin{proof}
Let $ \alpha \in ((E_{6,\sH})^C)^{\tau\lambda\sigma} $. Then it follows from $ (\tau\sigma) {}^t\!\alpha^{-1}(\sigma\tau)=\alpha $ that
\begin{align*}
\langle \alpha X, \alpha Y \rangle_\sigma &=(\tau \sigma \alpha X,\alpha Y)=({}^t\!\alpha^{-1}(\sigma\tau)X,\alpha Y)=(\sigma\tau X,\alpha^{-1}\alpha Y)=(\tau\sigma X,Y)
\\
&=\langle X, Y \rangle_\sigma.
\end{align*}
Hence we see $ \alpha \in E_{6(-14),\sH} $. Conversely, let $ \beta \in E_{6(-14),\sH} $. Then it follows from $ \langle \beta X, \beta Y \rangle_\sigma=\langle X ,Y \rangle_\sigma $ that
\begin{align*}
(\tau\sigma X,Y)=\langle X ,Y \rangle_\sigma=\langle \beta X, \beta Y \rangle_\sigma=(\tau\sigma\beta X,\beta Y)=({}^t\!\beta\tau\sigma \alpha X,Y)
\end{align*}
that $ \tau\sigma={}^t\!\beta\tau\sigma \beta $, that is, $ (\tau\sigma){}^t\!\beta^{-1}(\sigma\tau)=\beta $, so that $ \beta \in ((E_{6,\sH})^C)^{\tau\lambda\sigma} $.

With above, the proof of this theorem is completed.
\end{proof}

Let the mapping $ f:SU(2,4,\C^C) \to SU(6,\C^C) $ defined by $ f(B)=\varGamma_2B{\varGamma_2}^{-1} $, where $ \varGamma_2:=\diag(-i,-i,1,1,1,1) $ (\cite[in the proof of Theorem 3.5.11 (2)]{iy7}). Then this mapping $ f $ gives the isomorphism $ SU(2,4,\C^C) \cong SU(6,\C^C) $ (\cite[in the proof of Theorem 3.5.11 (2)]{iy7}). Using the mapping $ f $ and $ \phi $ defined in page 26, we can define an isomorphism $ \zeta:SU(2,4,\C^C) \to SU^*(6,\C^C)  $ by the composition mapping of $ f $ and $ \phi $: $ \zeta:=f \circ \phi $, and the explicit form of $ \zeta $ is give by $ \zeta(B)=\iota(\varGamma_2 B{\varGamma_2}^{-1})-\ov{\iota}J\ov{(\varGamma_2 B{\varGamma_2}^{-1}})J, \iota=(1/2)(1+i\i) $.
Hence the composition mapping $ f_{6,\sH^C}\zeta $ of $ \zeta $ and $ f_{6,\sH^C} $ induces the isomorphism $ (E_{6,\sH})^C \cong SU(2,4,\C^C)/\Z_2,\Z_2=\{E,-E\} $:
\begin{align*}
SU(2,4,\C^C) \overset{\zeta}{\longrightarrow} SU^*(6,\C^C) \overset{f_{6,{\scalebox{0.8}{\sH}^C}}}{\longrightarrow} (E_{6,\sH})^C.
\end{align*}
We denote the composition mapping $ f_{6,\sH^C}\zeta $ by $ l_{6,\sH^C} $: $ l_{6,\sH^C}:=f_{6,\sH^C}\zeta $.

Here, we prove the lemma needed in the proof of theorem below.

\begin{lemma}\label{lemma 5.4.2}
For $ B \in SU(2,4,\C^C) $, the mapping $ l_{6,\sH^C} $ satisfies $ (\tau\sigma) {}^t\!(l_{6,\sH^C}(B))^{-1}(\sigma\tau)=l_{6,\sH^C}\allowbreak (\tau B) $.
\end{lemma}
\begin{proof}
First, we have $ I_2\zeta(B)I_2=\zeta(I_2BI_2) $. Indeed, it follows from
\begin{align*}
I_2\zeta(B)I_2&=I_2(\iota(\varGamma_2 B{\varGamma_2}^{-1})-\ov{\iota}J\ov{(\varGamma_2 B{\varGamma_2}^{-1}})J)I_2=\iota(\varGamma_2 (I_2BI_2){\varGamma_2}^{-1})-\ov{\iota}J\ov{(\varGamma_2 (I_2BI_2){\varGamma_2}^{-1}})J
\\
&=\zeta(I_2BI_2).
\end{align*}
Hence it follows from Lemma \ref{lemma 5.1.2} (4) that
\begin{align*}
\sigma l_{6,\sH^C}(B)\sigma&=\sigma f_{6,\sH^C}(\zeta(B))\sigma=f_{6,\sH^C}(I_2\zeta(B)I_2)=f_{6,\sH^C}(\zeta(I_2BI_2))
\\
&=l_{6,\sH^C}(I_2BI_2),
\end{align*}
that is, $ \sigma l_{6,\sH^C}(B)\sigma=l_{6,\sH^C}(I_2BI_2) $.
Moreover, it follows that
\begin{align*}
\tau f(B)=\tau(\varGamma_2B{\varGamma_2}^{-1})={\varGamma_2}^{-1}(\tau B)\varGamma_2=\varGamma_2(I_2(\tau B)I_2){\varGamma_2}^{-1}=f(I_2(\tau B)I_2),
\end{align*}
that is, $ \tau f(B)=f(I_2(\tau B)I_2) $. Hence it follows from Lemma \ref{lemma 5.2.2} that
\begin{align*}
\tau{}^t\!l_{6,\sH^C}(B)^{-1}\tau&=\tau{}^t\!g_{6,\sH^C}(f(B))^{-1}\tau=g_{6,\sH^C}(\tau f(B))=g_{6,\sH^C}(f(I_2(\tau B)I_2))
\\
&=l_{6,\sH^C}(I_2(\tau B)I_2),
\end{align*}
that is, $ \tau{}^t\!l_{6,\sH^C}(B)^{-1}\tau=l_{6,\sH^C}(I_2(\tau B)I_2) $.

Thus it follows that
\begin{align*}
(\tau\sigma){}^t\!l_{6,\sH^C}(B)^{-1}(\sigma\tau)&=l_{6,\sH^C}(I_2(\tau (I_2BI_2))I_2)=l_{6,\sH^C}({I_2}^2(\tau B){I_2}^2)
\\
&=l_{6,\sH^C}(\tau B),
\end{align*}
that is, $ (\tau\sigma) {}^t\!(l_{6,\sH^C}(B))^{-1}(\sigma\tau)=l_{6,\sH^C}\allowbreak (\tau B) $.
\end{proof}

Now, we determine the structure of the group $ E_{6(-14),\sH} $.

\begin{theorem}\label{theorem 5.4.3}
The group $ E_{6(-14),\sH} $ is isomorphic to the group $ SU(2,4)/\Z_2 ${\rm :} $ E_{6(-14),\sH} \cong SU(2,4)/\Z_2 $.
\end{theorem}
\begin{proof}
Let the group $ E_{6(-14),\sH} $ as the group $ ((E_{6,\sH})^C)^{\tau\lambda\sigma} $. The we define a mapping $ l_{6,\tau\lambda\sigma}:SU(2,4) \to ((E_{6,\sH})^C)^{\tau\lambda\sigma} $ by
\begin{align*}
l_{6,\tau\lambda\sigma}(A)X=k^{-1}(\zeta(A)(k X)\zeta(A)^*), \;\;X \in \mathfrak{J}(3,\H^C).
\end{align*}
Note that this mapping is the restriction of the mapping $ l_{6,\sH^C} $.  First, we will prove that $ l_{6,\tau\lambda\sigma} $ is well-defined and a homomorphism. Since $ l_{6,\tau\lambda\sigma} $ is the restriction of the mapping $ l_{6,\sH^C} $, it is clear that $ g_{6,\tau\lambda}(A) \in (E_{6,\sH})^C $ and $ l_{6,\tau\lambda\sigma} $ is a homomorphism. Moreover, from Lemma \ref{lemma 5.4.2}, we have $ (\tau\sigma){}^t\!(l_{6,\tau\lambda\sigma}(A))^{-1}(\sigma\tau)=l_{6,\tau\lambda\sigma}(A) $, so that $ l_{6,\tau\lambda\sigma}(A) \in ((E_{6,\sH})^C)^{\tau\lambda\sigma} $.

Next, we will prove that $ l_{6,\tau\lambda\sigma} $ is surjective. Let $ \alpha \in E_{6(-14),\sH}=((E_{6,\sH})^C)^{\tau\lambda\sigma} \subset (E_{6,\sH})^C $. Then there exists $ B \in SU(2,4,\C^C) $ such that $ \alpha=l_{6,\sH^C}(B) $ (as mentioned in the beginning of page 32). Moreover $ \alpha $ satisfies the condition $ \tau{}^t\!\alpha^{-1}\tau=\alpha $, that is, $ (\tau\sigma){}^t\!(l_{6,\sH^C}(B))^{-1}(\sigma\tau)=l_{6,\sH^C}(B) $, so that from Lemma \ref{lemma 5.4.2} we have the following
\begin{align*}
\tau B=B \quad \text{or} \quad \tau B=-B.
\end{align*}
In the former case, we have $ B \in SU(2,4) $. Hence there exists $ A \in SU(2,4) $ such that $ \alpha=l_{6,\sH^C}(A)=l_{6,\tau\lambda\sigma}(A) $.
In the latter case, $ B $ is of the form
$ iB', B' \in M(6,\C) $, so that since $ B^*I_2B=I_2 $ and $ \det\,B=1 $, we have $ {B'}^*I_2B'=-I_2 $ and $ \det\,B'=-1 $. However, these leads to the following
\begin{align*}
-1&=\det\,(-I_2)=\det\,({B'}^*I_2B')=(\det\,{B'}^*)(\det\,I_2)(\det\,B')
\\
&=(\ov{\det\,B'})(\det\,I_2)(\det\,B')=(-1)1(-1)=1.
\end{align*}
Hence this is contradiction, so that this case is impossible.
 With above, the proof of surjective is completed.

Finally, we will determine $ \Ker\,l_{6,\tau\lambda\sigma} $. However, since $ l_{6,\tau\lambda\sigma} $ the restriction of the mapping $ l_{6,\sH^C} $, it is easy to obtain $ \Ker\,l_{6,\tau\lambda}=\Ker\,l_{6,\sH^C}=\{E,-E\} \cong \Z_2 $.

Therefore, from Theorem \ref{theorem 5.4.1}, we have the required isomorphism
\begin{align*}
E_{6(-13),\sH} \cong SU(2,4)/\Z_2.
\end{align*}

\end{proof}

\subsection{The group $ E_{6(-26),\sH} $}

As mentioned in the surface of this section, since the group $ (E_{6,\sH})^C $ has an involutive automorphism $ \tilde{\tau} $, we can consider the subgroup $ ((E_{6,\sH})^C)^{\tau} $ of $ (E_{6,\sH})^C $:
\begin{align*}
((E_{6,\sH})^C)^{\tau}:=\left\lbrace \alpha \in (E_{6,\sH})^C \relmiddle{|} \tilde{\tau}(\alpha)=\alpha \right\rbrace.
\end{align*}

Then we have the following theorem.

\begin{theorem}\label{theorem 5.5.1}
The group $ ((E_{6,\sH})^C)^{\tau} $ coincides with to the group $ E_{6(-26),\sH} ${\rm :} $ ((E_{6,\sH})^C)^{\tau}=E_{6(-26),\sH} $.
\end{theorem}
\begin{proof}
Let $ \alpha \in ((E_{6,\sH})^C)^{\tau} $. Then it follows from $ \tau\alpha=\alpha\tau $ that $ \alpha X=\alpha (\tau X)=\tau (\alpha X), X \in \mathfrak{J}(3,\H) $, that is, $ \alpha X \in \mathfrak{J}(3,\H) $, so that $ \alpha $ indices an $ \R $-linear isomorphism of $ \mathfrak{J}(3,\H) $. Hence we see $ \alpha \in E_{6(-26),\sH} $. Conversely, let $ \beta \in E_{6(-26),\sH} $. Then we define an action of $ \beta $ to $ \mathfrak{J}(3,\H^C) $ by
\begin{align*}
\beta X=\beta(X_1+iX_2)=\beta X_1+ i\beta X_2,\;\; X \in \mathfrak{J}(3,\H^C), X_i \in \mathfrak{J}(3,\H).
\end{align*}
Hence $ \beta $ induces a $ C $-linear isomorphism of $ \mathfrak{J}(3,\H^C) $. Moreover, it follows that
\begin{align*}
\beta X \times \beta Y&=\beta(X_1+iX_2) \times \beta(Y_1+iY_2)=(\beta X_1+i\beta X_2) \times (\beta Y_1+i\beta Y_2)
\\
&=(\beta X_1 \times \beta Y_1-\beta X_2 \times \beta Y_2)+i(\beta X_1 \times \beta Y_2+\beta X_2 \times \beta Y_1)
\\
&=({}^t\!\beta^{-1}(X_1 \times Y_1)-{}^t\!\beta^{-1}(X_2 \times Y_2))+i({}^t\!\beta^{-1}(X_1 \times Y_2)+{}^t\!\beta^{-1}(X_2 \times Y_1))
\\
&={}^t\!\beta^{-1}((X_1 \times Y_1-X_2 \times Y_2)+i(X_1 \times Y_2+X_2 \times Y_1))
\\
&={}^t\!\beta^{-1}((X_1+iX_2)\times (Y_1+iY_2))
\\
&={}^t\!\beta^{-1}(X \times Y),
\end{align*}
so that $ \beta \in (E_{6,\sH})^C $. In addition, we have
\begin{align*}
\tau\beta X&=\tau\beta(X_1+iX_2)=\tau(\beta X_1+i\beta X_2)=\beta X_1-i\beta X_2=\beta(X_1-iX_2)=\beta\tau(X_1+iX_2)
\\
&=\beta\tau X,
\end{align*}
that is, $ \tau\beta=\beta\tau $. Hence we see $ \beta \in ((E_{6,\sH})^C)^{\tau} $.

With above, the proof of this theorem is completed.
\end{proof}

\begin{theorem}\label{theorem 5.5.2}
{\rm (1)} The Lie algebra $ \mathfrak{e}_{6(-26),\sH} $ of the group $ E_{6(-26),\sH} $ is given by
\begin{align*}
  \mathfrak{e}_{6(-26),\sH}=\left\lbrace \phi \in \Hom_{\sR}(\mathfrak{J}(3,\H))\relmiddle{|}(\phi X, X, X)=0 \right\rbrace.
\end{align*}

{\rm (2)} Any element $ \phi \in \mathfrak{e}_{6(-26),\sH} $ is uniquely expressed by the form
\begin{align*}
    \phi=\delta+\tilde{T},\;\; \delta \in \mathfrak{f}_{4,\sH}, T \in \mathfrak{J}_{\sH}, \tr(T)=0,
\end{align*}
where $ \mathfrak{f}_{4,\sH}:=\left\lbrace \phi \in \mathfrak{e}_{6(-26),\sH} \relmiddle{|} \phi E=0 \right\rbrace  $.

In particular, we have $ \dim(\mathfrak{e}_{6(-26),\sH})=35 $.
\end{theorem}
\begin{proof}
(1) The proof is evident (cf. \cite[Lemma 2.3.1]{iy0}).

(2) First, we have to prove the following
\begin{align*}
   \mathfrak{f}_{4,\sH}=\left\lbrace \phi \in \mathfrak{e}_{6(-26),\sH} \relmiddle{|} \phi E=0 \right\rbrace.
\end{align*}
Indeed, we can easily prove this as follows:
\begin{align*}
    \mathfrak{f}_{4,\sH}:&=\left\lbrace \phi \in \Iso_{\sR}(\mathfrak{J}(3,\H))\relmiddle{|}\phi(X \circ Y)=\phi X \circ Y +X \circ \delta Y \right\rbrace
    \\
    &=\left\lbrace \phi \in \Iso_{\sR}(\mathfrak{J}(3,\H))\relmiddle{|}
    (\phi X, X, X)=0, \phi E=0 \right\rbrace
    \\
    &=\left\lbrace \phi \in \mathfrak{e}_{6(-26),\sH} \relmiddle{|}\phi E=0
    \right\rbrace.
\end{align*}
Here, set $ T:=\phi E $, then we have
\begin{align*}
    T \in \mathfrak{J}(3,\H),\;\; \tr(T)=0.
\end{align*}
Indeed, it is trivial that $ T \in \mathfrak{J}(3,\H) $, and it follows that $ \tr(T)=(T,E,E)=(\phi E,E,E)=0 $. Moreover, set $ \delta:=\phi-\tilde{T} $, then we have
\begin{align*}
    \delta \in \mathfrak{f}_{4,\sH}.
\end{align*}
Indeed, since it follows from \cite[Lemma 2.3.5]{iy0} that
\begin{align*}
    (\tilde{T}X, X, X)&=(T \circ X, X ,X)=((1/2)(TX+XT),X \times X)
    \\
    &=((1/2)TX,X \times X)+((1/2)XT,X \times X)
    \\
    &=(T,(1/2)X(X \times X))+(T,(1/2)(X\times X)X)
    \\
    &=(T,(1/2)(X(X \times X)+(X \times X)X))
    \\
    &=(T, X \circ (X \times X))=(T,(\det\,X)E)
    \\
    &=(\det\,X)(T,E)=(\det\,X) \tr(T)
    \\
    &=0,
\end{align*}
we see $ \tilde{T} \in \mathfrak{e}_{6(-26),\sH} $, that is, $ \delta \in \mathfrak{e}_{6(-26),\sH} $. In addition, we have $ \delta E=(\phi-\tilde{T})E=\phi E-\tilde{T}E=T-T=0 $, so that $ \delta \in \mathfrak{f}_{4,\sH} $. Thus we have the following
\begin{align*}
    \phi=\delta + \tilde{T},\;\; \delta \in \mathfrak{f}_{4,\sH} \cong \mathfrak{sp}(3), T \in \mathfrak{J}(3,\H), \tr(T)=0.
\end{align*}
Note that $ \mathfrak{f}_{4,\sH} \cong \mathfrak{sp}(3) $ is a direct result of the isomorphism $ F_{4,\sH} \cong Sp(3)/\Z_2 $ (\cite[Proposition 2.11.1]{iy0}).

Finally, we will prove the uniqueness of its expression. In order to prove this, it is sufficient to show that $ \delta+\tilde{T}=0 $ implies $ \delta=0 $ and $ T=0 $. Certainly, let apply it on $ E $, then we have $ T=0 $, so that $ \delta=0 $.

Therefore we have $ \dim(\mathfrak{e}_{6(-26),\sH})=21+(2+4 \times 3)=35 $.
\end{proof}

In order to give a polar decomposition of the group $ E_{6(-26),\sH} $, we prove the following lemma.

\begin{lemma}\label{lemma 5.5.3}
 The group $  E_{6(-26),\sH} $ is an algebraic subgroup of the general linear group $ GL(15,\R)=\Iso_{\sR}(\mathfrak{J}(3,\H)) $ and satisfies the condition that $ \alpha \in E_{6(-26),\sH} $ implies $ {}^t\alpha \in E_{6(-26),\sH} $, where $ {}^t\alpha $ in the transpose of $ \alpha $ with respect to the inner product $ (X,Y) ${\rm :} $ (\alpha X, Y)=(X,{}^t\alpha Y) $.
\end{lemma}
\begin{proof}
 We use the identity formula: $ (Z \times Z)\times (Z \times Z)=(\det \,Z)Z,\;\; Z \in \mathfrak{J}(3,\H) $. For $ \alpha \in E_{6(-26),\sH} $ and $ Y \in \mathfrak{J}(3,\H) $, we have
 \begin{align*}
      {}^t\alpha^{-1}(Y \times Y)\times {}^t\alpha^{-1}(Y \times Y)&=(\alpha Y \times \alpha Y)\times (\alpha Y \times \alpha Y)
      =(\det\,\alpha Y)(\alpha Y)
      \\
      &=(\det\, Y)(\alpha Y)=\alpha((\det\, Y)Y)
      \\
      &=\alpha((Y\times Y)\times (Y\times Y)).
    \end{align*}
Here, set $ Y:=X \times X $ for any $ X \in \mathfrak{J}(3,\H) $ in the formula above, then we have
\begin{align*}
    {}^t\alpha^{-1}((\det\,X)X)\times {}^t\alpha^{-1}((\det\,X)X)=\alpha((\det\,X)X\times (\det\,X)X),
\end{align*}
that is, $ (\det\,X) ({}^t\alpha^{-1} X\times {}^t\alpha^{-1} X)=(\det\,X)\alpha(X \times  X) \cdots (*) $.
\vspace{1mm}

(i) Case where $ \det\,X \not=0 $. Then we have $ {}^t\alpha^{-1} X\times {}^t\alpha^{-1} X=\alpha (X\times X) $, that is, $ \det\,{}^t\alpha^{-1}X=\det\,X $. Thus, by considering $ \alpha^{-1} $ instead of $ \alpha $, we have
\begin{align*}
    \det\,{}^t\alpha X=\det\,X.
\end{align*}

(ii) Case where $ \det\,X=0 $. Then we have the same result above. Indeed, assume $ \det\,{}^t\alpha X \not=0 $, then by replacing $ X $ with $ {}^t\alpha X=:X' $ in the formula $ (*) $ above, we have $ (\det\,X') ({}^t\alpha^{-1} X' \times {}^t\alpha^{-1} X')=(\det\,X')\alpha(X' \times  X'), \det\,X' \not=0 $. Hence, as in the case (i), we have
\begin{align*}
    \det\,{}^t\alpha^{-1}X'=\det\,X',
\end{align*}
so that it follows that
\begin{align*}
    0=\det\,X=\det\,{}^t\alpha^{-1}({}^t\alpha X)\underset{*}{=}\det\,
    ({}^t\alpha X).\;\;(\mbox{\boldmath{$\underset{*}{=}\gets \det\,{}
    ^t\alpha^{-1}X'=\det\,X'$}})
\end{align*}
This is contradiction. Thus we obtain $ \det\,{}^t\alpha X =0=\det\,X
$, that is, $ {}^t\alpha \in E_{6(-26),\sH} $.

Finally, since the group $ E_{6(-26),\sH} $ is defined by the algebraic relation $ \det\,\alpha X=\det\,X $, it is clear that $ E_{6(-26),\sH} $ is real algebraic.
\end{proof}

Let $ O(\mathfrak{J}(3,\H)) $ be the orthogonal subgroup of $ GL(15,\R)=\Iso_{\sR}(\mathfrak{J}(3,\H)) $:
\begin{align*}
    O(15)=O(\mathfrak{J}(3,\H)):=\left\lbrace \alpha \in \Iso_{\sR}(\mathfrak{J}(3,\H))  \relmiddle{|} (\alpha X, \alpha Y)=(X,Y)\right\rbrace.
\end{align*}
Then, from \cite[Proposition 2.11.1]{iy0} we have the following
\begin{align*}
    E_{6(-26),\sH} \cap O(\mathfrak{J}(3,\H))&=\left\lbrace \alpha \in \Iso_{\sR}(\mathfrak{J}(3,\H))  \relmiddle{|} \det\,\alpha X=\det\,X,(\alpha X, \alpha Y)=(X,Y)\right\rbrace
    \\
    &=F_{4,\sH}
    \\
    &\cong Sp(3)/\Z_2.
\end{align*}

Using Chevalley's lemma (\cite[Lemma 2]{che}), we have a homeomorphism
\begin{align*}
    E_{6(-26),\sH} &\simeq (E_{6(-26),\sH} \cap O(\mathfrak{J}(3,\H))) \times R^d
    \\
    &\simeq Sp(3)/\Z_2 \times \R^d,
\end{align*}
where the dimension $ d $ of the Euclidean part is computed by Theorem
\ref{theorem 5.5.2} as follows:
\begin{align*}
    d=\dim(E_{6(-26),\sH})-\dim(Sp(3))=35-21=14.
\end{align*}

With above, we have the following theorem.

\begin{theorem}\label{theorem 5.5.4}
  The group $ E_{6(-26),\sH} $ is homeomorphic to the topological product of the group $ Sp(3)/\Z_2 $ and a $ 14 $-dimensional Euclidean space $ \R^{14} $:
  \begin{align*}
      E_{6(-26),\sH} \simeq Sp(3)/\Z_2 \times \R^{14}.
  \end{align*}

In particular, $ E_{6(-26),\sH} $  is a connected Lie group.
\end{theorem}

Now, we determine the structure of the group $ E_{6(-26),\sH} $.

\begin{theorem}\label{theorem 5.5.5}
The group $ E_{6(-26),\sH} $ is isomorphism to the group $ SU^*(6)/\Z_2 , \Z_2=\{E,-E\}${\rm :} $ E_{6(-26),\sH} \cong SU^*(6)/\Z_2 $.
\end{theorem}
\begin{proof}
Let the group $ E_{6(-26),\sH} $ as the group $ ((E_{6,\sH})^C)^{\tau} $ (Theorem \ref{theorem 5.5.1}).
Then we define a mapping $ f_{6,\tau}:SU^*(6) \to ((E_{6,\sH})^C)^{\tau} $ by
\begin{align*}
f_{6,\tau}(A)X=k^{-1}(A(k X)A^*),\;\; X \in \mathfrak{J}(3,\H^C).
\end{align*}
Note that this mapping $ f_{6,\tau} $ is the restriction of the mapping $ f_{6,\sH^C} $.
First, we will prove that $ f_{6,\tau} $ is well-defined a homomorphism. Since the mapping $ f_{6,\tau} $ is the restriction of the mapping $ f_{6,\sH^C} $, we easily see that $ f_{6,\tau} $ ia a homomorphism,
and together with Lemma \ref{lemma 5.1.2} (2), we have $ f_{6,\tau}(A) \in ((E_{6,\sH})^C)^{\tau} $.

Next, we will determine $ \Ker\,f_{6,\tau} $. Since the mapping $ f_{6,\tau} $ is the restriction of the mapping $ f_{6,\sH^C} $, we easily obtain $ \Ker\,f_{6,\tau}=\Ker\,f_{6,\sH^C}=\{E,-E\} \cong \Z_2. $

Finally, we will prove that $ f_{6,\tau} $ is surjective. Since the group $ ((E_{6,\sH})^C)^\tau=E_{6,\sH} $ is connected (Theorems \ref{theorem 5.5.4}, \ref{theorem 5.5.1}) and $ \Ker\,f_{6,\tau} $ is discrete, together with $ \dim(((\mathfrak{e}_{6,\sH})^C)^\tau=\mathfrak{e}_{6,\sH})=35=\dim(\mathfrak{su}^*(6)) $ (Lemma \ref{lemma 5.2.2} (2)), we see that $ f_{6,\tau} $ is surjective.

Therefore, from Theorem \ref{theorem 5.5.1}, we have the required isomorphism
\begin{align*}
E_{6(-26),\sH} \cong SU^*(6)/\Z_2.
\end{align*}
\end{proof}


\end{document}